\documentclass[a4paper,11pt,DIV=11,abstract=on%
]{scrartcl}
\usepackage[utf8]{inputenc}
\usepackage{fontenc}
\usepackage[english]{babel}
\usepackage{csquotes}
\usepackage{amsmath, amsthm, amssymb,mathtools}
\usepackage{graphicx}
\usepackage{mathrsfs}
\usepackage{subcaption}
\usepackage{xargs}
\usepackage{algorithm}
\usepackage[noend]{algpseudocode}
\usepackage{booktabs}
\usepackage{enumerate}
\usepackage{listings}
\usepackage{environ,ifthen,xcolor}
\usepackage{hyperref}
\usepackage{animate}
\usepackage{pgfplots,tikz}
\usetikzlibrary{spy,calc,external}
\usepackage{bbm}
\mathtoolsset{showonlyrefs}
\definecolor{tuklblue}{RGB}{0,95,140}

\usepackage[bordercolor=tuklblue,linecolor=tuklblue,backgroundcolor=tuklblue!20]{todonotes}
\DeclareCaptionLabelSeparator{periodspace}{.\ }
\setkomafont{sectioning}{\rmfamily\bfseries}
\setkomafont{title}{\rmfamily}
\newcommand{\e}{{\,\mathrm{e}}}

\newcommand{\E}{\mathbb{E}}
\newcommand{\N}{\mathbb{N}}
\renewcommand{\P}{\mathbb{P}}

\newcommand{\R}{\mathbb{R}}

\newcommand{\DD}{\mathcal{D}}

\newcommand{\GG}{\mathcal{G}}
\newcommand{\HH}{\mathcal{H}}

\newcommand{\LL}{\mathcal{L}}

\newcommand{\NN}{\mathcal{N}}

\renewcommand{\SS}{\mathcal{S}}

\newcommand{\dx}[1][x]{\,\mathrm{d}#1}

\newcommandx{\abs}[2][1=\@empty]{#1\lvert #2 #1\rvert}
\newcommandx{\norm}[3][1=\@empty,3=\@empty]{#1\lVert #2 #1\rVert_{#3}}

\DeclareMathOperator*{\argmin}{arg\,min} 
\DeclareMathOperator*{\argmax}{arg\,max} 

\DeclareMathOperator{\divv}{div}

\DeclareMathOperator{\MSE}{MSE}
\DeclareMathOperator{\TV}{TV}
\newtheorem{theorem}{Theorem}[section]

\newtheorem{proposition}[theorem]{Proposition}
\newtheorem{lemma}[theorem]{Lemma}
\newtheorem{remark}[theorem]{Remark}

\newtheorem{corollary}[theorem]{Corollary}

\begin{document}
\title{Nonlocal Myriad Filters\\ for Cauchy Noise Removal}
\author{Friederike Laus\footnotemark[1] \and Fabien Pierre\footnotemark[3] \and Gabriele Steidl\footnotemark[1] \footnotemark[2]}

\date{\today}

	\maketitle
\footnotetext[1]{Department of Mathematics,
	Technische Universität Kaiserslautern,
	Paul-Ehrlich-Straße~31, D-67663 Kaiserslautern, Germany,
	\{friederike.laus,steidl\}@mathematik.uni-kl.de.}
\footnotetext[2]{Fraunhofer ITWM, Fraunhofer-Platz 1,
	D-67663 Kaiserslautern, Germany} 
\footnotetext[3]{
	Laboratoire Lorrain de Recherche en Informatique et ses Applications, UMR CNRS 7503, Université de Lorraine, INRIA projet Magrit, France, fabien.pierre@univ-lorraine.fr} 
\renewcommand{\proofname}{\textbf{\emph{Proof:}}}
\allowdisplaybreaks[4]
\newlength\figureheight 
\newlength\figurewidth

\noindent	
	
	\begin{abstract}
		The contribution of this paper is two-fold.
		First, we introduce a generalized myriad filter, which is a method to compute  the joint maximum likelihood estimator 
		of the location and the scale parameter of the Cauchy distribution. 
		Estimating only the location parameter is known as myriad filter. 
		We propose an efficient algorithm to compute the generalized myriad filter and prove its convergence. 
		Special cases of this algorithm result in the classical myriad filtering
		and an algorithm for estimating only the scale parameter. 
		Based on an asymptotic analysis, we develop a second, even faster generalized myriad filtering technique.
		
		Second, we use our new approaches within a nonlocal, fully unsupervised  method to denoise images corrupted by Cauchy noise. 
		Special attention is paid to the determination of similar patches in noisy images.
		Numerical examples demonstrate the excellent performance of our algorithms which have moreover the advantage to be robust with respect to the parameter choice.
	\end{abstract}

	\section{Introduction}
	%
	Myriad filters (MF) as introduced in~\cite{GA96}  form a large class of nonlinear filters for robust non-Gaussian signal processing. 
	Analogously as  mean and median filters can be derived from Gaussian and Laplacian distributions respectively, 
	the myriad filter arises from 
	the Cauchy distribution, a special kind of \emph{$\alpha$-stable distributions}, which is due to its heavy tails often used to model an impulsive  behavior. 
	Examples of  such data can be found in low-frequency atmospheric
	signals~\cite{SG74}, underwater acoustic signals~\cite{BA13}, 
	radar clutter~\cite{KFR98}, and multiple-access interference in wireless communication systems~\cite{Mid77}, to mention only a few. 
	MFs have been successfully employed in robust signal and image processing~\cite{Arc05}. 
	Applications in image processing include myriad filtering 
	to denoise images corrupted by Gaussian plus $\alpha$-stable noise~\cite{ZGA96}, whereas in~\cite{HK01}, 
	convex combinations of local mean and median filters or local mean and myriad filters 
	have been used to denoise images corrupted by Salt-and-Pepper and mixed Gaussian-Laplacian noise. 
	
	The class of MFs can be derived from the \emph{sample myriad}, which is the maximum likelihood (ML) estimator of the location parameter of the Cauchy distribution. 
	More general, the MF belongs to the class of so called M-\emph{estimators}, which are estimators obtained as  minima of sums of functions of given data, that means
	\begin{equation*}
	\hat{\theta} \in \argmin_{\theta}\sum_{i=1}^n \rho(x_i;\theta).
	\end{equation*}
	Here, the function $\rho$ is the \emph{cost function} of the M-estimator and choosing e.g.\  $\rho(x;\theta) = -\log\bigl(f(x|\theta)\bigr)$ 
	for a probability density function (pdf) corresponds to classical ML estimator.
	Taking the density functions of the Gaussian and Laplacian
	distribution  yields the cost functions $\rho(x) = x^2$ for the sample mean 
	and $\rho(x) = |x|$ for the sample median, and the sample myriad is obtained using the density of the Cauchy distribution, 
	resulting in the cost function $\rho(x) = \log(x^2 + \gamma^2)$, where the so-called \emph{scale parameter} $\gamma>0$ controls the robustness of the estimator.

	In general, estimating the parameter(s) of the Cauchy distribution is a difficult task, 
	since they are not related to any moments of the distribution or transformations thereof; 
	in fact, the Cauchy distribution has no finite moments. 
	Besides ML estimators~\cite{HBA70}, there exist several  
	other approaches 
	in the literature, 
	for instance  order statistics~\cite{Bar66,Blo66,RFT64}, sample quantiles~\cite{Can74}, window estimates~\cite{HT77}, empirical
	characteristic function~\cite{BM01,Kou82,MT05}, Bayesian estimators~\cite{HW88}, L-estimators~\cite{Zha10}, Pitman estimator~\cite{Fre07}
	or linear rank estimators~\cite{BBT03}. 
	As always in statistical estimation, one has to find the trade off between  the ease of computation and the 
	efficiency, robustness and consistency of the estimator. 
	ML estimation possesses  a number of desirable limiting properties such as consistency, asymptotic normality and efficiency. 
	A requirement for this are sufficiently large sample sizes, which is usually the case in image processing applications.
	For the Cauchy distribution, MFs have shown to be optimal in \cite{GA01}.
	However, the MFs do not admit a closed-form expression, which makes its computation a nontrivial task. 
	Different computation methods have been proposed, ranging from fixed point algorithms~\cite{KA00} over branch-and-bound search~\cite{NRGA08} 
	to polynomial and trimming approaches~\cite{Pan10,Pan16}. 
	
	In this paper, we propose to combine MFs with nonlocal methods in image processing.
	Nonlocal, patch-based methods have shown to provide state-of-the-art results 
	in many image restoration tasks, and the concept of non-locality is in particular central to
	most of the recent denoising techniques. These include the nonlocal means algorithm~\cite{BCM05} and its generalizations~\cite{GO08,Sal10,WPCMB07,YSM12}, 
	BM3D~\cite{DFKE08} and BM3D-SAPCA~\cite{DFKE09}, patch-ordering based wavelet methods~\cite{REC14},
	and the nonlocal Bayes algorithm~\cite{LBM13b,LBM13}.
	For a recent review of the denoising problem and the different denoising principles 
	we refer to~\cite{LCBM12} and for a generalization to manifold-valued images to \cite{LNPS17}.
	The standard noise model is that of additive Gaussian white noise and  the quality of the denoised images 
	has become excellent for moderate noise levels, as summarized e.g.\ in \cite{CM10,LN11}.
	
	However, in many situation the image acquisition process suggests other noise models.
	Recently, the Cauchy distribution attracted attention in image processing.  In~\cite{MDHY16,SDZ15}, the authors proposed
	a variational method  for removing Cauchy noise. Their model can also be used in the context of inverse problems, e.g.\ for deblurring. 
	Following a Bayesian approach, the model consists of the data term 
	$\mathcal{D}(u;f)\coloneqq  \int_{\Omega} \log( (u-f)^2 + \gamma^2 ) \dx $
	which resembles the noise statistics and a total variation regularization term.
	The first model of Sciacchitano, Dong and Zeng \cite{SDZ15} contains additionally a quadratic term 
	$\|u - u_0\|^2$  relating the image $u$ to the median filtered version $u_0$
	of the noisy image $f$ and making the whole model convex.
	In \cite{MDHY16}, the artificial quadratic term is skipped and the nonconvex minimization problem is
	solved by a nonconvex ADMM version of Wang, Yin and Zeng \cite{WYZ15}. 
	This model shows a better performance than the convexified one.
	For both methods, the scale parameter $\gamma$ of the Cauchy distribution has to be known.
	
	In this paper, we propose a generalized myriad filter  (GMF) for estimating both the location and the scale parameters of the Cauchy distribution. 
	The classical MF for estimating only the location parameter~$a$ as well as an algorithm for estimating only the scale parameter~$\gamma$ 
	can be obtained as special cases of our algorithm. Additionally, we analyze the asymptotic behavior of the algorithms if the sample size approaches infinity.
	To this end, we use the interesting fact that applying the functions appearing in the gradient of the log-likelihood function
	to a Cauchy distributed random variable results in a random variable which possesses first and second order moments.
	Our considerations result in a further, even faster GMF algorithm. For all our algorithms we provide convergence results.
	
\begin{figure*}[thb]
		\centering
		\begin{subfigure}[t]{0.24\textwidth}
			\centering
			\includegraphics[width=.98\textwidth]{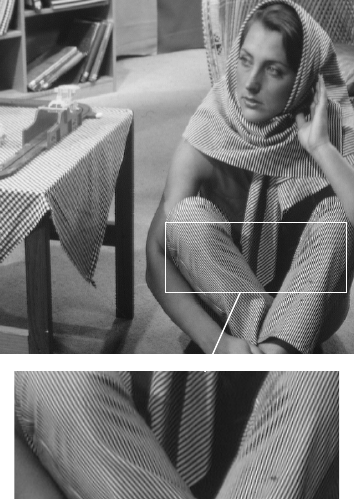}
		\end{subfigure}
		\begin{subfigure}[t]{0.24\textwidth}
			\centering
			\includegraphics[width=.98\textwidth]{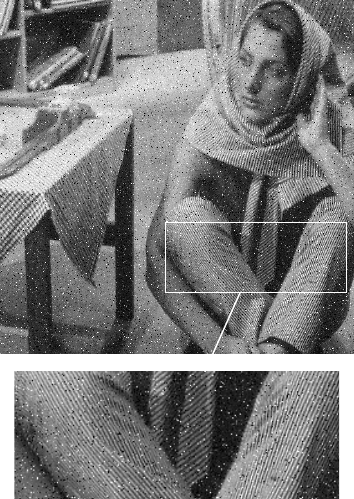}
		\end{subfigure}
		\begin{subfigure}[t]{0.24\textwidth}
			\centering
			\includegraphics[width=.98\textwidth]{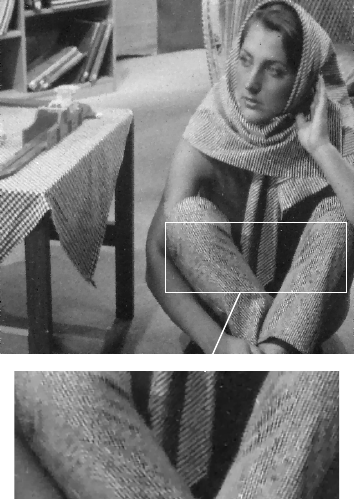}
		\end{subfigure}
		\begin{subfigure}[t]{0.24\textwidth}
			\centering
			\includegraphics[width=.98\textwidth]{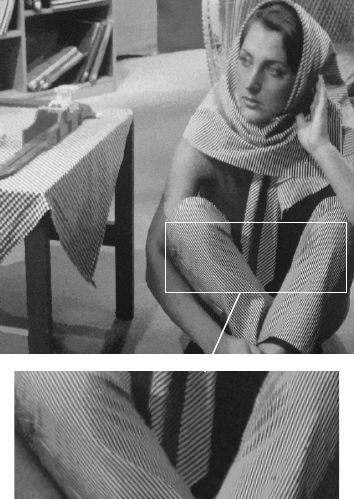}
		\end{subfigure}
		
		\vspace{0.15cm}
		
		\begin{subfigure}[t]{0.24\textwidth}
			\centering
			\includegraphics[width=.98\textwidth]{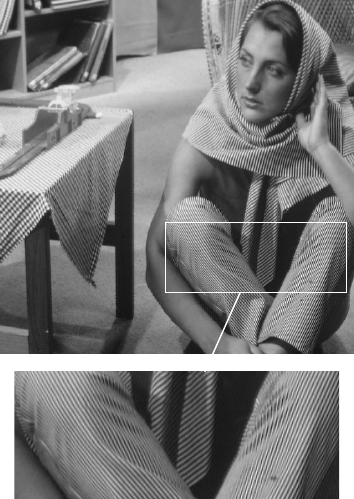}
		\end{subfigure}
		\begin{subfigure}[t]{0.24\textwidth}
			\centering
			\includegraphics[width=.98\textwidth]{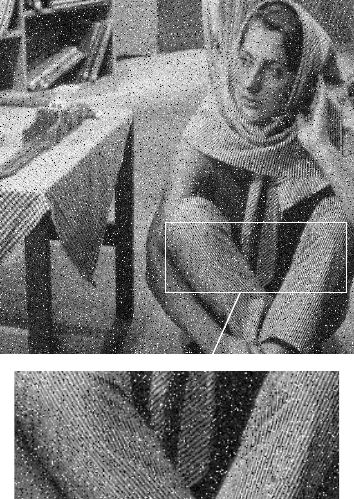}
		\end{subfigure}
		\begin{subfigure}[t]{0.24\textwidth}
			\centering
			\includegraphics[width=.98\textwidth]{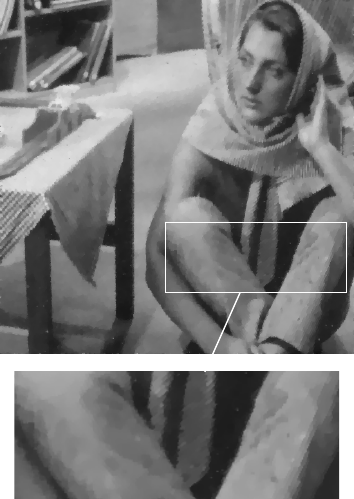}
		\end{subfigure}
		\begin{subfigure}[t]{0.24\textwidth}
			\centering
			\includegraphics[width=.98\textwidth]{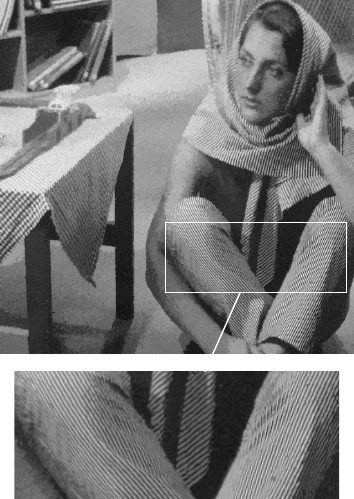}
		\end{subfigure}		
		\caption{Results  of our proposed nonlocal generalized myriad filtering  
			in comparison with the variational method~\cite{MDHY16} 
			for different noise levels $\gamma = 5$ (top) and $\gamma = 10$ (bottom). 
			From left to right: Original image, noisy image, result of~\cite{MDHY16}, our result.  }
		\label{Fig:ADMM_vs_NGMF}
	\end{figure*}

	In a second step, we apply our GMF algorithm  to design a nonlocal myriad filter, called~\emph{nonlocal generalized  myriad filter} (N-GMF). 
	At this point, it is interesting to notice that our  generalized nonlocal myriad filter which estimates both the location and the scale parameter 
	does not only perform better than a local myriad filter, which is to be expected, but also better than a nonlocal myriad filter.
	Further, an important issue in our nonlocal approach is the selection of similar patches, which will serve as samples in the myriad filter. 
	Here, the adaptation of the similarity measure to the noise
	distribution is essential for a robust similarity evaluation, see~\cite{DDT12}. We propose a similarity measure based on 
	a likelihood ratio statistical test under the hypothesis of Cauchy noise corruption. 
	Figure~\ref{Fig:ADMM_vs_NGMF} indicates the very good performance of our new N-GMF in comparison with the method in~\cite{MDHY16}.
	
	{\textbf{Outline of the paper}}: 
	In Section \ref{sec:cauchy} we briefly introduce the Cauchy distribution. 
	Section~\ref{sec:ML_prop} deals with the ML approach to estimate the parameters of a Cauchy distribution, 
	and, in particular, with properties of the (log-)likelihood function.
	Some of these properties are known, however, we incorporated the section for the following  reasons: 
	Beside we wish to make the paper self-contained, most of the results are needed in the following sections, 
	and further, there are some gaps in proofs  in the literature.
	In Section~\ref{sec:alg}, we propose a new algorithm (GMF) for maximizing the log-likelihood function in both 
	the  location and scale parameter. 
	Based on this algorithm, we additionally derive two algorithms for single parameter estimation.
	We provide convergence proofs for the various algorithms. Moreover, 
	we set our results for maximizing the log-likelihood function in the location parameter 
	in relation to an algorithm from~\cite{KA00}. In addition to~\cite{KA00} we take also 
	saddle points of the objective function into consideration 
	and show that we have indeed convergence to a local maximum with probability one. 
	Section~\ref{sec:asymp} provides an asymptotic analysis of the proposed algorithm if the sample size approaches infinity.
	Based on this analysis we develop an improved GMF algorithm.
	Next, in Section~\ref{sec:myr} we apply the developed parameter estimation algorithms within a nonlocal denoising method method.  
	In particular, we explain how to estimate the noise level and how to find similar patches  in images corrupted by Cauchy noise.
	Several numerical examples and a comparison with the variational method from~\cite{MDHY16}
	are given in Section \ref{sec:numerics}. Finally, Section \ref{sec:conclusions} contains  conclusions and an outline for future work. The proofs of the results of Sections~\ref{sec:ML_prop}-\ref{sec:myr} can be found in Appendices~\ref{app:ML_prop}-\ref{app:myr}, respectively.

	\section{Cauchy Distribution}\label{sec:cauchy}
	
	The Cauchy distribution $C(a,\gamma)$ 
	belongs to the class of \emph{$\alpha$-stable distributions}. 
	Recall that for $\alpha\in (0,2]$, a random variable $X$ is said to have an $\alpha$-stable distribution, if for all $n\in \N$ and i.i.d.\ random variables $X_1,\ldots,X_n,X$ it holds
	$
	X_1 + \ldots + X_n\sim n^{\frac{1}{\alpha}}X
	$.
	In general, there is no analytic solution for the density of the resulting distribution. 
	However, in three special cases the resulting distribution can be further specified, namely  $\alpha = 2$ results in the normal distribution, 
	$ \alpha = 1$ in the Cauchy distribution and $\alpha=\frac{1}{2}$ in the L\'{e}vy distribution.
	
\begin{figure*} [thb]
		\centering  
		\centering  
		{\includegraphics[width=0.4\textwidth]{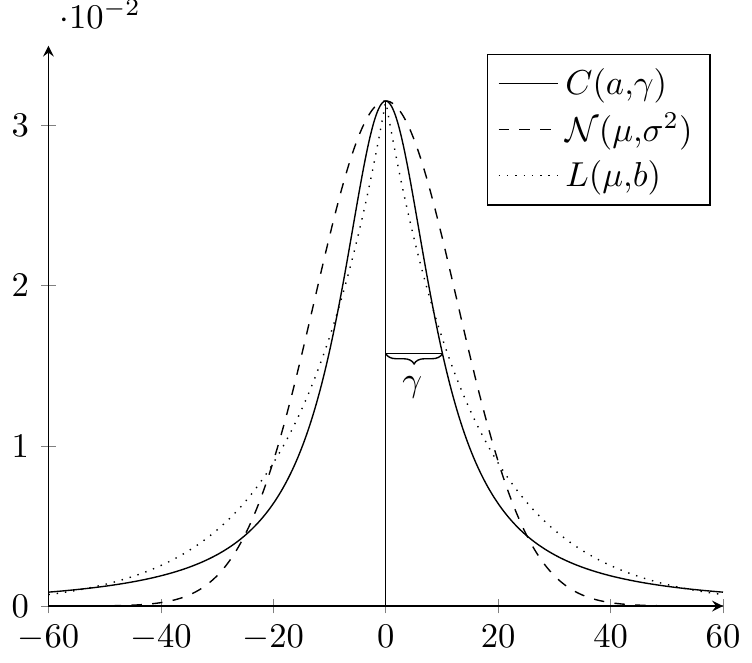}}
		\caption{Comparison of the Gaussian, Laplacian and Cauchy distribution for $\mu = a = 0$,   $\sigma^2 = 160$, $b = \sqrt{\frac{\pi}{2}\sigma^2}\approx 15.8533 $ and $\gamma^2 = \frac{2\sigma^2}{\pi}\approx 10.0925$.
		}\label{Fig:Cauchy_tails}
	\end{figure*}

	The Cauchy distribution depends on two parameters,
	the \emph{location} parameter $a\in \R$ and the \emph{scale parameter} $\gamma >0$. 
	Its probability density function (pdf) and cumulative density function (cdf) are given by
	\begin{align}
	p(x|a,\gamma) &= \frac{1}{\pi \gamma}  \frac{\gamma^2}{(x-a)^2 + \gamma^2} = \frac{1}{\pi \gamma}\frac{1}{\left(\frac{x-a}{\gamma} \right)^2+1},\\
	F(x|a,\gamma) & = \frac{1}{\pi} \arctan\left(\frac{x-a}{\gamma} \right) + \frac{1}{2}\label{cdf}. 
	\end{align}
	The parameter $a$ is at the same time median and mode of the distribution, 
	whereas $\gamma$ specifies the half-width at half-maximum (HWHM), alternatively $2\gamma$ is full width at half maximum (FWHM), see Figure~\ref{Fig:Cauchy_tails} . In contrast to the normal distribution, no moments of the Cauchy distribution exist, in particular its mean and variance are undefined.
	
	Because of the fact that the parameters of the Cauchy distribution are not related to any moment, 
	attempting to estimate the parameters of the Cauchy distribution  using  sample mean, sample variance or any transformation thereof does not work.  
	The fact that the Cauchy distribution has no finite mean and variance is closely related to its heavy-tailedness. A distribution is said to be right \emph{heavy-tailed}, if 
	\begin{equation*}
	\lim_{x\to \infty} \e^{\lambda x} \P(X>x) = + \infty \qquad \text{for all }\lambda >0,
	\end{equation*}
	and similarly for heavy left tails.
	
	The heavy-tails of the Cauchy distribution are illustrated in Figure~\ref{Fig:Cauchy_tails}, where the Cauchy distribution $C(0,\gamma)$
	is compared to the normal distribution $\NN(0,\sigma^2)$ for $\sigma^2 = 160$ and to the Laplacian distribution $L(0,b)$, where $b = \sqrt{\frac{\pi}{2}\sigma^2}\approx 15.8533  $ and  $\gamma = \sqrt{\frac{2\sigma^2}{\pi}}\approx 10.0925$. The parameters $b$ and $\gamma$ are chosen in such a way that the densities of the two distributions coincide in zero. Due to the heavy tails, the probability of outliers is, in case of the Cauchy distribution, much higher than the Gaussian and also different from the Laplace distribution, see also noisy images in Figure~\ref{Fig:ADMM_vs_NGMF}.

	In analogy to the normal distribution, the distribution $C(0,1)$ is called \emph{standard Cauchy distribution}, and using e.g.\ characteristic functions
	one easily verifies the following result.
	
	\begin{proposition} \label{Prop:properties_Cauchy}
		If $X\sim C(a,\gamma)$ and $Y = \alpha X + \beta$ for some $\alpha, \beta\in \R$, then it holds $Y\sim C(\alpha a + \beta, |\alpha|\gamma)$. 
		Further, for independent random variables $X_1,\ldots,X_n$, $X_i\sim C(a_i,\gamma_i)$,  
		the relation $\sum\limits_{i=1}^n X_i \sim C\left(\sum\limits_{i=1}^{n}a_i,\sum\limits_{i=1}^{n}\gamma_i\right)$
		is fulfilled.
	\end{proposition}

	Sampling of the Cauchy distribution is straightforward and can be done for instance using the inversion method. 
	Indeed, from~\eqref{cdf} one sees that if $U\sim \mathcal{U}(0,1)$ is uniformly distributed on $[0,1]$, then 
	\begin{equation*}
	X = a+\gamma \tan\left(\pi(U-\tfrac{1}{2})\right)\sim C(a,\gamma). 
	\end{equation*}
	Alternatively, one might use the fact that $\frac{X}{Y}\sim  C(0,1)$ for i.i.d.\ $X,Y\sim \NN(0,1)$ in order to sample from $C(0,1)$ and then use Proposition~\ref{Prop:properties_Cauchy} to obtain samples from $C(a,\gamma)$.

	\section{Properties of the ML Function}\label{sec:ML_prop}
	
	In this section, we (re)consider properties of the (log-) likelihood function of the Cauchy distribution.
	We start with the analysis of the joint likelihood function. It is known that in this case, the maximizer is unique, 
	and most authors refer to the paper of Copas~\cite{Cop75} for this result. However, as already pointed out by Gabrielsen~\cite{Gab82}, 
	the proof of Copas is incomplete since he only showed that all critical points  are maxima. Here, 
	a further step is required to conclude that this leads to the desired uniqueness result.
	Gabrielsen mentioned that Morse theory can be used to show a slightly more general result 
	and that another ,,proof can be obtained from the author on request''. In \cite{Scho96} an argument is given for Gabrielsen's result using properties 
	of the solutions of ordinary differential equations.
	We apply in the second part of the proof of Theorem \ref{theo_both} the fact that simple roots of polynomials vary  
	smoothly with the coefficients of the polynomial  to show the final result. 
	Then, we provide in Lemmata~\ref{lem_1} to~\ref{lem_2} properties of the likelihood function if one of the parameters is fixed. In particular, 
	Lemma~\ref{Lemma:no_saddle_point} deals with the occurrence of saddle points of the likelihood function for fixed scale parameter. 
	It  makes use of the resolvent of polynomials and plays an important role for the analysis of the classical myriad filter. The proofs of this section can be found in Appendix~\ref{app:ML_prop}.
	
	Let $x_1,\ldots,x_n$ be i.i.d.\ realizations of a Cauchy random variable $X \sim C(a,\gamma)$. In what follows we aim at estimating the parameter $\theta = (a,\gamma)$ using the ML approach. 
	For the Cauchy distribution, the likelihood function reads 
	\begin{align*}
	{\mathcal L}(a, \gamma|x_1,\ldots,x_n) &= \prod_{i=1}^n p(x_i|a,\gamma) = \left(\frac{\gamma}{\pi}\right)^n \prod_{i=1}^n \frac{1}{(x_i-a)^2 + \gamma^2},
	\end{align*}
	and the log-likelihood function is given by
	\begin{equation*}
	\log {\mathcal L} (a,\gamma|x_1,\ldots,x_n) = n \log(\gamma) - n \log(\pi) - \sum_{i=1}^n \log\bigl((x_i-a)^2 + \gamma^2\bigr).
	\end{equation*}
	Maximizing $\log {\cal L}$ it is equivalent to minimizing $-\log {\cal L}$, on which we focus in the following. Here and in the following, 
	the notation ${\mathcal L}(a, \gamma|x_1,\ldots,x_n)$ does not refer to conditional probabilities, but emphasizes that $\LL$ depends on the samples $x_1,\ldots,x_n$.
	
	More generally, we might also allow for different weighting of the summands,
	that is, we introduce positive weights $w_1,\ldots,w_n$, $w_i> 0$, $\sum\limits_{i=1}^n w_i = 1$, and consider the function
	\begin{align*}
	L(a,\gamma)&\coloneqq L(a,\gamma|x_1,\ldots,x_n)=  \sum_{i=1}^n w_i  \log\bigl((x_i-a)^2 + \gamma^2\bigr)-\log(\gamma).
	\end{align*}
	Note that by doing so we might assume
	\begin{equation*}
	x_1 < \ldots < x_n, \qquad n \ge 2,	
	\end{equation*}
	since multiple samples can be handled by updating the weight accordingly. 
	We  skip the dependence of $L$ on the samples $x_1,\ldots,x_n$ since it will always be  clear from the context.
	
	We will need the gradient 
	$\nabla L = \begin{pmatrix} \frac{\partial L}{\partial a} \\[0.5ex] \frac{\partial L}{\partial \gamma} \end{pmatrix}$ 
	and the Hessian 
	$\nabla^2 L = 
	\begin{pmatrix} 
	\frac{\partial^2 L}{\partial a^2} & \frac{\partial^2 L}{\partial a \partial \gamma }\\[0.5ex]
	\frac{\partial^2 L}{\partial a \partial \gamma} & \frac{\partial^2 L}{\partial \gamma^2} 
	\end{pmatrix}
	$
	of $L$.
	The partial derivatives of $L$ are given by
		\begin{align}
		\frac{\partial L}{\partial a}(a,\gamma) &= 2\sum_{i=1}^n w_i\frac{a-x_i}{(x_i-a)^2 + \gamma^2} \label{der_a},\\
		\frac{\partial L}{\partial \gamma}(a,\gamma) &= 2\sum_{i=1}^n w_i \frac{\gamma}{(x_i-a)^2 + \gamma^2}-\frac{1}{\gamma} \label{der_gamma},\\
		\frac{\partial^2 L}{\partial a^2} (a,\gamma) &= - 2 \sum\limits_{i=1}^n w_i \frac{ (x_i-a)^2- \gamma^2}{\bigl((x_i-a)^2 + \gamma^2\bigr)^2} \label{der_a2},\\
		\frac{\partial^2 L}{\partial a \partial \gamma} (a,\gamma)&= 4 \sum\limits_{i=1}^n w_i\frac{\gamma(x_i-a)}{\bigl((x_i-a)^2 + \gamma^2\bigr)^2} \label{der_agamma},\\
		\frac{\partial^2 L}{\partial \gamma^2}(a,\gamma) &=  2 \sum\limits_{i=1}^n w_i \frac{(x_i-a)^2-\gamma^2}{\bigl((x_i-a)^2 + \gamma^2\bigr)^2}+\frac{1}{\gamma^2}\label{der_gamma2}.
		\end{align}
	Then, $(\hat a,\hat \gamma)$ is a critical point of $L$ if and only if it solves the system of equations 
	\begin{align}
	\sum_{i=1}^n & w_i\frac{x_i-\hat a}{(x_i- \hat a)^2 + \hat \gamma^2} = 0,\label{cond_a}\\
	\sum_{i=1}^n& w_i \frac{\hat \gamma^2}{(x_i-\hat a)^2 + \hat \gamma^2} = \frac{1}{2} \label{cond_gamma}.
	\end{align}
	These equations are in general not solvable in closed form, except $n\leq 4$. 
	For $n=3,4$,    closed form solutions are given in~\cite{Fer78}.
	For later usage we introduce the notation
	\begin{align}
	S_0(a, \gamma)& \coloneqq \sum_{i=1}^n w_i \frac{\gamma^2}{(x_i-a)^2 + \gamma^2} =\frac{1}{2} +  \frac{\gamma}{2}\frac{\partial L}{\partial \gamma}(a,\gamma)  \label{cond_gamma_1},\\
	S_1(a, \gamma) &\coloneqq \sum_{i=1}^n w_i\frac{\gamma (x_i-a)}{(x_i- a)^2 + \gamma^2}=-\frac{\gamma}{2}\frac{\partial L}{\partial a}(a,\gamma). \label{cond_a_1}\\
	\end{align}

	The next theorem recalls that the likelihood function of the Cauchy distribution has a unique maximizer. We  add the proof for two reasons: 
	The first part of the proof following the arguments of Copas~\cite{Cop75} is given to make the paper
	self-contained. The second part contains a new proof that the conclusions from the first part 
	indeed yield the desired result.
	
	\begin{theorem}\label{theo_both}
		Let $n \ge 3$ and $w_i < \frac12$, $i=1,\ldots,n$. 
		Then $L$ has exactly one critical point on $\R\times \R_{>0}$. This point is a minimizer of $L$. 
	\end{theorem}

	Next, we fix one of the parameters and consider the functions 
	$ L(\cdot,\gamma)$ and $L(a,\cdot)$.

	\begin{lemma} \label{lem_1}
		Let $n \ge 2$.
		Then, for fixed $\gamma > 0$, the function 
		$L(\cdot,\gamma)$
		has at least one and at most $2n-1$ critical points. All critical points lie in the interval $(x_1,x_n)$.
	\end{lemma}

	In dependence on the samples $x_i$, $i=1,\ldots,n$ we can further characterize the critical points of $L(\cdot,\gamma)$.
	
	\begin{lemma}\label{Lemma:no_saddle_point} 
		For fixed $\gamma > 0$, the function $L(\cdot,\gamma|x_1,\ldots,x_n)$ has $\lambda^n $-a.s. no saddle points, i.e.
		\begin{equation}
		\lambda^n\Bigl(\bigl\{ 
		(x_i)_{i=1}^n : \ \exists\, \hat a \in \R \colon 
		\frac{\partial L}{\partial a}(\hat a,\gamma|x_1,\ldots,x_n) 
		= \frac{\partial^2 L}{\partial a^2}(\hat a,\gamma|x_1,\ldots,x_n) = 0 \bigr\}\Bigr)=0,\label{saddle_point}
		\end{equation}	
		where $\lambda^n$ denotes the $n$-dimensional Lebesgue measure.
	\end{lemma}

	Finally, we consider the likelihood function $L(a,\cdot)$ for a fixed location parameter $a$. Here, we have the following result:
	
	\begin{lemma} \label{lem_2}
		Let $n \ge 3$ and $w_i < \frac12$, $i=1,\ldots,n$.
		Then, for fixed $a \in (x_1,x_n)$, the function 
		$L(a,\cdot)$
		has exactly one critical point. This critical point is a minimizer of $L(a,\cdot)$ and lies
		in $(d \epsilon,x_n - x_1)$, where $d \coloneqq \min_i (x_{i+1} - x_i)$, $\epsilon \coloneqq  (\frac12 - w_{\max})^\frac12$ and $w_{\max} \coloneqq \max\limits_{i=1,\ldots,n}\{ w_i\}$.
	\end{lemma}

	\section{Algorithms for Parameter Estimation}\label{sec:alg}
	In this section, we propose a new algorithm for minimizing
	$L(\cdot,\cdot)$. Then, fixing the parameter $\gamma$, resp. $a$ in this algorithm, provides efficient algorithms for
	computing a minimizer of $L(\cdot, \gamma)$ and the minimizer of $L(a,\cdot)$.
	The first algorithm coincides with those suggested by Kalluri and Arce in~\cite{KA00}.
	Unfortunately, the convergence proof in~\cite{KA00} is incomplete. We provide a simpler proof, 
	that takes in particular saddle points of $L(\cdot, \gamma)$ into account. The proofs of this section can be found in Appendix~\ref{app:alg}.

	\subsection{Minimization of $L(\cdot,\cdot)$}
	We start with the joint minimization of $L(\cdot,\cdot)$.
	Reformulating~\eqref{cond_a} yields
	\begin{align*}
	\sum_{i=1}^n w_i \frac{\hat{a}}{(x_i-\hat{a})^2 + \hat{\gamma}^2} &= 	\sum_{i=1}^n w_i \frac{x_i}{(x_i-\hat{a})^2 + \hat{\gamma}^2},
	\end{align*}  
	which gives rise to the semi-implicit iteration
	\begin{align}
	a_{r+1}
	&= \frac{\sum\limits_{i=1}^n w_i\frac{x_i}{( x_i-a_r )^2 + \gamma_r^2}}{\sum\limits_{i=1}^n w_i\frac{1}{( x_i-a_r )^2 + \gamma_r^2} }\label{update_a_convex}\\
	&= a_r + \gamma_r \frac{S_{1}(a_r,\gamma_r)}{S_{0}(a_r,\gamma_r)},\label{fp_a}\\
	&=  a_r  - \frac{\gamma_r^2}{2S_0(a_r,\gamma_r)} \frac{\partial L}{\partial a}(a_r,\gamma_r).
	\end{align}
	Similarly, we obtain by \eqref{cond_gamma}
	\begin{align}
	\gamma^2_{r+1}  &=\gamma_r^2  \frac{1-S_0(a_r,\gamma_r)}{S_0(a_r,\gamma_r)} , \label{fp_gamma}	\\
	&=\gamma_r^2 - \frac{\gamma_r^3}{S_0(a_r,\gamma_r)} \frac{\partial L}{\partial \gamma}(a_r,\gamma_r).
	\end{align}
	It turns out that the above balancing of $S_{0}(a_r,\gamma_r)$ and $1-S_{0}(a_r,\gamma_r)$ leads to a fast convergence. More precisely, the convergence is monotone, which is stated in Theorem~\ref{Theo:convergence_gamma}. Combining~\eqref{fp_a} and~\eqref{fp_gamma} results in Algorithm \ref{Alg:myriad_general}. 
	It  would also be possible to use the new $a_{r+1}$ in the update of~$\gamma_r$, that means, to replace the update rule for $\gamma$ by
	\begin{equation*}
	\gamma_{r+1}^2 \coloneqq \gamma_r^2 \, \frac{1-S_0(a_{r+1},\gamma_r)}{S_0(a_{r+1},\gamma_r)}.
	\end{equation*}
	However, in our numerical experiments, this iteration scheme was slower than Algorithm \ref{Alg:myriad_general}. 
	
	\begin{algorithm}[!ht]
		\caption{Minimization of $L(\cdot,\cdot)$ (Generalized Myriad Filter, GMF) }\label{Alg:myriad_general}
		\begin{algorithmic}
			\State \textbf{Input:} $x_1< \ldots < x_n$, $n \ge 3$, $0< w_i < \frac12$, $i=1,\ldots,n$, $\sum\limits_{i=1}^n w_i = 1$,
			\State \textbf{Initialization:} $a_0 \in (x_1,x_n)$, $\gamma_0 >0$
			\For{$r=0,\ldots$}
			\begin{align} \label{update_ag}		
			a_{r+1} &\coloneqq  a_r + \gamma_r \frac{S_1(a_r,\gamma_r)}{S_0(a_r,\gamma_r)}  \\
			\gamma_{r+1}^2 &\coloneqq \gamma_r^2 \, \frac{1-S_0(a_r,\gamma_r)}{S_0(a_r,\gamma_r)}
			\end{align}
			\EndFor
		\end{algorithmic}
	\end{algorithm}
	
	The following theorem establishes the convergence of the proposed algorithm.
	
	\begin{theorem}\label{Theo:myriad_general}
		Let $n \ge 3$ and $w_i < \frac12$, $i=1,\ldots,n$.
		Then, for any starting point $\gamma_0 >0$, the sequence $\{ (a_r,\gamma_r) \}_{r\in \N}$ 	
		generated by  Algorithm \ref{Alg:myriad_general}
		converges to the minimizer
		$(\hat a,\hat \gamma)$ of $L$.		
	\end{theorem}

	\subsection{Minimization of $L(\cdot,\gamma)$}
	For a given $\gamma > 0$,  
	Algorithm~\ref{Alg:myriad_general} can be replaced by  Algorithm \ref{Alg:myriad_a} to obtain a minimizer 
	of $L(\cdot,\gamma)$. It turns out that in this situation our algorithm coincides with the one proposed in~\cite{KA00}.
	\begin{algorithm}[!ht]
		\caption{Minimization of $L(\cdot,\gamma)$ (Myriad Filter, MF)}\label{Alg:myriad_a}
		\begin{algorithmic}
			\State \textbf{Input:} $x_1< \ldots < x_n$, $n \ge 2$, $w_i > 0$, $i=1,\ldots,n$, $\sum\limits_{i=1}^n w_i = 1$,
			$\gamma > 0$
			\State \textbf{Initialization:} $a_0 \in (x_1,x_n)$
			\For{$r=0,\ldots$}
			\begin{align}		
			a_{r+1} &\coloneqq  a_r + \gamma \frac{S_1(a_r,\gamma)}{S_0(a_r,\gamma)}
			\end{align}
			\EndFor
		\end{algorithmic}
	\end{algorithm}
	

	%
	
	The next theorem guarantees the convergence of Algorithm~\ref{Alg:myriad_a}. 
	Slightly weaker results can be found in \cite{KA00}, where, however, the proof is much more complicated. Further, in \cite{KA00} among others the treatment of saddle points is missing. 
	
	We define the following function
	\begin{equation}
	Q(a) \coloneqq L(a,\gamma) + \log(\gamma) = \sum_{i=1}^n w_i \log \left((x_i-a)^2 + \gamma^2 \right)	
	\end{equation}
	and associate to Algorithm~\ref{Alg:myriad_a} the operator $T_1$ given by $a_{r+1} \coloneqq T_1(a_r)$. With these definitions we have the following result:
	\begin{theorem}\label{Theo:descent}
		Let $n \ge 2$. Then, for every starting point $a_0 \in (x_1,x_n)$, the sequence $\{ a_r\}_ {r\in \N}$  
		generated by Algorithm \ref{Alg:myriad_a}
		converges to a critical point of $Q$.
	\end{theorem}

	Theorem~\ref{Theo:descent} establishes the convergence to a critical point of $Q$, that is, a fixed point of $T_1$. 
	However, according to Lemma~\ref{lem_1} there might exist up to $2n-1$ fixed points and consequently 
	the starting point $a_0$ determines the computed fixed-point.
	Further, it is a  priori not guaranteed that we end up in a (local) minimum.
	That this is a.s.\ the case  is stated in the next theorem. Note that in \cite{KA00} the treatment of saddle points is missing.
	
	\begin{theorem}\label{Theorem:convergence_local_min}
		Let $n \geq 2$ and $a_0 \in (x_1,x_n)$ be an arbitrary starting point. 
		Then, $\lambda$-a.s.\ the sequence $\{ a_r\}_ {r\in \N}$ generated by
		Algorithm~\ref{Alg:myriad_a} converges to a local minimum of $Q$, i.e.\ 
		\begin{equation*}
		\lambda\Bigl(\bigl\{a_0\in (x_1,x_n) \colon  \hat a = \lim_{r\to \infty} a_r\text{ is not a local minimum of }Q     \bigr\}\Bigr)=0.
		\end{equation*}
	\end{theorem}

	\subsection{Minimization of $L(a,\cdot)$}
	Next, we consider the minimization of $L(a,\cdot)$ for fixed $a \in \mathbb R$. 
	In this case,  Algorithm \ref{Alg:myriad_general} simplifies to Algorithm~\ref{Alg:myriad_gamma}.
	
	
	\begin{algorithm}[!ht]
		\caption{Minimization of $L(a,\cdot)$ }\label{Alg:myriad_gamma}
		\begin{algorithmic}
			\State \textbf{Input:} $x_1< \ldots < x_n$, $n \ge 3$, $0< w_i < \frac12$, $i=1,\ldots,n$, $\sum\limits_{i=1}^n w_i = 1$,
			$a \in (x_1,x_n) > 0$
			\State \textbf{Initialization:} $\gamma_0 \in (0,x_n-x_1)$
			\For{$r=0,\ldots$}
			\begin{align} \label{update_gamma}
			\gamma_{r+1}^2 \coloneqq \gamma_r^2 \, \frac{1-S_0(a,\gamma_r)}{S_0(a,\gamma_r)}		
			\end{align}
			\EndFor
		\end{algorithmic}
	\end{algorithm}
	
	Here, we have the following convergence result.
	
	\begin{theorem}\label{Theo:convergence_gamma}
		Let $n \ge 3$ and $w_i < \frac12$, $i=1,\ldots,n$.
		Then, for any starting point $\gamma_0 \in (0, x_n - x_1)$, the sequence $\{ \gamma_r \}_{r\in \N}$ 	
		generated by Algorithm \ref{Alg:myriad_gamma} converges to the minimizer
		$\hat \gamma$ of $L(a,\cdot)$.	
		Furthermore, we have monotone convergence in the sense that one of the
		following relation is fulfilled for all $r\in \N$
		\begin{align} \label{mono}
		\gamma_r \leq \gamma_{r+1} \leq \hat \gamma \quad \mathrm{or} \quad
		\gamma_r \ge \gamma_{r+1} \ge \hat \gamma,
		\end{align}
		where $\gamma_r = \gamma_{r+1}$ if and only if $\gamma_r = \hat \gamma$.
	\end{theorem}


	\begin{remark}[Initialization of the Algorithms]
		Initialization is not an issue when estimating the scale parameter $\gamma$ or both parameters, as we have global convergence in these cases.  It is, however, crucial in the case of estimating only the location parameter $a$ in Algorithm~\ref{Alg:myriad_a}, since the result depends on the starting point of the iterative scheme. Further, a good initialization improves of course the speed of convergence, see also the simulation study in Subsection~\ref{Sec:Simulation}. One strategy is to initialize the algorithms with estimates that can be easily computed, for instance the sample median for $a$ and the Hodge-Lehman-estimator~\cite{KP12} for~$\gamma$. In our experiments, both perform very well.
		Concerning the initialization of Algorithm~\ref{Alg:myriad_a}, we further observed that  generally the global minimum is located near the mode of the samples $x_1,\ldots,x_n$, and in particular one of them is usually very close to it. Thus, 
		to initialize our algorithm, we choose the value $x_i$ for which the objective function becomes minimal, that is
		\begin{equation}\label{init_algo}
		a_0 = \argmin_{x_i\in \{ x_1, \dots, x_n \} } Q(x_i). 
		\end{equation}
		Although this does not guarantee the convergence to the global minimum, it turned out that in all our numerical experiments, Algorithm~\ref{Alg:myriad_a} initialized with~\eqref{init_algo} converges to the global minimum. 
	\end{remark}
	\section{Asymptotic Analysis and Speed of Convergence} \label{sec:asymp}
	The fact that the Cauchy distribution has no finite moments makes it hard to estimate the parameters 
	$a$ and~$\gamma$ compared to for instance the parameters of a normal distribution.
	Furthermore, asymptotic results such that the \emph{Law of Large Numbers} 
	or the \emph{Central Limit Theorem} are not applicable. 
	This is however different when considering transformed Cauchy random variables 
	as they appear in myriad objective function, and in the following, we aim at analyzing 
	the distribution of these transformed random variables in more details. 
	In particular, we prove an asymptotic convergence result including details 
	on the speed of convergence, which might serve as an explanation of   
	the good performance of our generalized  myriad filter. 
	Finally, the analysis leads to new   algorithm whose convergence is even faster than the convergence of the previous algorithms. The proofs of this section can be found in Appendix~\ref{app:asymp}.
	
	\subsection{Expectation of Transformed Cauchy Random Variables}
	To this aim, let $X_1,\ldots, X_n$ be i.i.d.\ random variables, $X_i\sim C(a,\gamma)$. 
	In order to simplify notation we define $X_{ir} = \frac{X_i-a_r}{\gamma_r}$ and consider as before
	\begin{align}
	S_0(a_r,\gamma_r) &\coloneqq \sum_{i=1}^n w_i \frac{\gamma_r^2}{( X_i-a_r)^2 + \gamma_r^2}\\
	&=  \sum_{i=1}^n w_i\frac{1}{1 + \left(\frac{ X_i-a_r}{\gamma_r}\right)^2} =  \sum_{i=1}^n w_i \frac{1}{1 + X_{ir}^2}, \\
	S_1(a_r,\gamma_r) &\coloneqq \sum_{i=1}^n w_i\frac{\gamma_r( X_i-a_r)}{( X_i-a_r)^2 + \gamma_r^2} =  \sum_{i=1}^n w_i \frac{X_{ir}}{1 + X_{ir}^2},
	\end{align}
	that is, we use the same notation as in \eqref{cond_gamma} and \eqref{cond_a_1}, but with random variables $X_i$ instead of samples $x_i$.
	The quantities $S_{0}$ and $S_{1}$ contain transformation of Cauchy random variables $X\sim C(a,\gamma)$ 
	of the form $\frac{1}{1 + X^2}$ and $\frac{X}{1 + X^2}$, whose expectation we  analyze in the next lemma.

	\begin{lemma}\label{Theo:expected_values}
		For $X\sim C(a,\gamma)$,  the random variables $Y = \frac{1}{1 + X^2}$ and $Z = \frac{X}{1 + X^2}$ fulfill
		\begin{equation} \label{E1}
		E(Y) = \begin{cases}
		\frac{\gamma(a^2 + \gamma^2-1) 
			+ a^2-\gamma^2 + 1}{(a^2 + \gamma^2 + 1)^2 - 4\gamma^2}  & \text{ for } a\neq 0,\\
		\frac{1}{1+\gamma} & \text{ for } a= 0,
		\end{cases}		
		\end{equation}
		and 
		\begin{equation}\label{E2}
		E(Z) = \begin{cases}
		\frac{a(a^2 + \gamma^2+1-2\gamma) }{(a^2 + \gamma^2 + 1)^2 - 4\gamma^2} & \text{ for } \; (a,\gamma) \not = (0,1),\\
		0 & \text{ for }  \; (a,\gamma)  = (0,1).
		\end{cases}		
		\end{equation}
	\end{lemma}
	
	As a direct consequence of Lemma~\ref{Theo:expected_values} we obtain the following corollary.
	
	\begin{corollary}\label{Coro:mean_Y_ir}
		Let $X\sim C(a,\gamma)$ and  $X_r \coloneqq \frac{X-a_r}{\gamma_r}$. Then, $X_r\sim C\left(\tfrac{a-a_r}{\gamma_r},\tfrac{\gamma}{\gamma_r}\right)$, and it holds
		\begin{align*}
		\E\bigl(g(X_r)\bigr) &= \frac{\gamma_r(\gamma + \gamma_r)}{(a-a_r)^2 + (\gamma + \gamma_r)^2}\qquad \text{and}\\ \E\bigl(h(X_r)\bigr) 
		&	= \frac{\gamma_r(a-a_r) }{(a-a_r)^2 + (\gamma + \gamma_r)^2}.
		\end{align*}
	\end{corollary}
	Corollary~\ref{Coro:mean_Y_ir} can be used to compute the expectation  of $S_i(a_r,\gamma_r)$, $i=0,1$. Indeed, since the random variables $X_i$ are i.i.d.\ and $\sum\limits_{i=1}^n w_i =1$,  we have
		\begin{align*}
		m_{0}\bigl(S_0(a_r,\gamma_r)\bigr)  \coloneqq \E\bigl(S_0(a_r,\gamma_r)\bigr)&=  \frac{\gamma_r(\gamma+\gamma_r)}{(a - a_r)^2 + (\gamma+\gamma_r)^2},\\
		m_{1} \bigl(S_1(a_r,\gamma_r)\bigr) : = \E\bigl(S_1(a_r,\gamma_r)\bigr) & = \frac{\gamma_r(a - a_r)}{(a - a_r)^2 + (\gamma+\gamma_r)^2}.
		\end{align*}

	%
	%
	
	\subsection{Asymptotic Analysis}
	To emphasize the dependence of $S_0$, $S_1$ on the sample size $n$, 
	we write $S^{(n)}_{0}$, $S^{(n)}_{1}$ in the following. 
	According to the Strong Law of Large Numbers, see, e.g.~\cite{CB02}, it holds
	\begin{align*}
	S^{(n)}_{0}(a_r,\gamma_r) &\overset{{a.s.}}{\longrightarrow} m_{0}(a_r,\gamma_r),\\
	S^{(n)}_{1}(a_r,\gamma_r)&\overset{{a.s.}}{\longrightarrow} m_{1}(a_r,\gamma_r),
	\end{align*}
	as $n\to \infty$
	and since $m_{0}(a_r,\gamma_r)> 0$ by the Continuous Mapping Theorem further
	\begin{align*}
	\frac{S^{(n)}_{1}(a_r,\gamma_r)}{S^{(n)}_{0}(a_r,\gamma_r)}&\overset{a.s.}{\longrightarrow}  \frac{m_{1}(a_r,\gamma_r)}{m_{0}(a_r,\gamma_r)},\\
	\frac{1}{S^{(n)}_{0}(a_r,\gamma_r)}&\overset{{a.s.}}{\longrightarrow} \frac{1}{m_{0}(a_r,\gamma_r)}
	\end{align*}
	as $n\to \infty$.
	As a consequence, for each $r\in \N$ we have
		\begin{align*}
		a_{r+1} & = a_r + \gamma_r \frac{S^{(n)}_{1}(a_r,\gamma_r)}{S^{(n)}_{0}(a_r,\gamma_r)}\overset{a.s.}{\longrightarrow}a_r + \gamma_r \frac{m_{1}(a_r,\gamma_r)}{m_{0}(a_r,\gamma_r)},\\
		\gamma^2_{r+1} & = \gamma^2_r \left( \frac{1}{S^{(n)}_{0}(a_r,\gamma_r)}-1\right)\overset{{a.s.}}{\longrightarrow} \gamma^2_r \left( \frac{1}{m_{0}(a_r,\gamma_r)}-1\right)
		\end{align*}	
	as $n\to \infty$.
	Thus, for a sufficiently large sample size, this suggests to replace $S_{0r}$ and $S_{1r}$ in Algorithm \ref{Alg:myriad_general}
	by their expected values $m_{0r}$ and $m_{1r}$, respectively. Setting $m_{ir}=m_{i}(a_r,\gamma_r)$, $i=1,2$, this results in sequences $\{\tilde{a}_r\}_{r\in \N}$ and $\{\tilde{\gamma}_r\}_{r\in \N}$ given by
	\begin{align}
	\tilde{a}_{r+1} &= \tilde{a}_r + \tilde{\gamma}_r \frac{m_{1r}}{m_{0r}}\label{seq1}\\
	& = \tilde{a}_r + \tilde{\gamma}_r \frac{\frac{\tilde{\gamma}_r(a-\tilde{a}_r)}{(a-\tilde{a}_r)^2 
			+ (\gamma + \tilde{\gamma}_r)^2}}{\frac{\tilde{\gamma}_r(\gamma+\tilde{\gamma}_r)}{(a - \tilde{a}_r)^2 + (\gamma+\tilde{\gamma}_r)^2}}\nonumber\\
	& = \tilde{a}_r + \tilde{\gamma}_r \frac{a-\tilde{a}_r}{\gamma + \tilde{\gamma}_r}\label{a_rek_mean}
	\end{align}
	and
	\begin{align}
	\tilde{\gamma}^2_{r+1} 
	& = \tilde{\gamma}^2_r \left( \frac{1}{m_{0r}}-1\right) \label{seq2}\\
	& = \tilde{\gamma}_r^2 \left(\frac{(a-\tilde{a}_r)^2 + (\gamma + \tilde{\gamma}_r)^2}{\tilde{\gamma}_r (\gamma + \tilde{\gamma}_r)}-1\right)\\
	&= \tilde{\gamma}_r \frac{(a-\tilde{a}_r)^2 + \gamma^2 + \gamma \tilde{\gamma}_r}{ \gamma + \tilde{\gamma}_r} \nonumber\\
	& = \tilde{\gamma}_r \left(\gamma + \frac{(a-\tilde{a}_r)^2}{ \gamma + \tilde{\gamma}_r}\right). \label{gamma_rek_mean}
	\end{align}
	Below we collect some properties of the generated sequences.
	
	\begin{theorem} \label{conv_komisch}
		Let $\tilde{a}_0\in \R$ and $\tilde{\gamma}_0>0$ be arbitrary. 
		Then the sequences $\{(\tilde{a}_r,\tilde{\gamma}_r)\}_{r\in \N}$  generated by \eqref{seq1}
		and \eqref{seq2} have the following properties:  
		\begin{enumerate}[\normalfont (i)]
			\item It holds $\tilde{\gamma}^2_{r+1}\geq \min\{\tilde{\gamma}^2_r,\gamma^2\}\geq \min\{\tilde{\gamma}^2_0,\gamma^2\}$.
			\item If $\tilde{a}_r\neq a$, then \begin{equation*}
			\tilde{a}_{r+1}\begin{cases}
			> \tilde{a}_r & \text{  if } \; \tilde{a}_r< a,\\
			< \tilde{a}_r  & \text{  if } \; \tilde{a}_r> a.
			\end{cases}
			\end{equation*}		
			\item The sequence converges linearly, $\lim\limits_{r\to \infty} (\tilde{a}_r, \tilde{\gamma}_r) = (a,\gamma)$ with rate $q = \max\left\{\frac{1}{2},\frac{\gamma}{\gamma + \tilde{\gamma}_0}\right\}$.
		\end{enumerate}
	\end{theorem}

	Solving ,,right-hand side of~\eqref{seq1}=\eqref{a_rek_mean}'' 
	and ,,right-hand side of~\eqref{seq2}=\eqref{gamma_rek_mean}'' for $a$ and $\gamma$, yields
	\begin{align*}
	a& = \tilde{a}_r + \tilde{\gamma}_r\frac{m_{1}(\tilde a_r,\tilde \gamma_r)}{m_{0}^2(\tilde a_r,\tilde \gamma_r) + m_{1}^2(\tilde a_r,\tilde \gamma_r)},\\
	\gamma & =  \tilde{\gamma}_r\left(\frac{m_{0}(\tilde a_r,\tilde \gamma_r)}{m_{0}^2(\tilde a_r,\tilde \gamma_r) + m_{1}^2(\tilde a_r,\tilde \gamma_r)}-1\right).\\
	\end{align*}
	In practice, the quantities $m_{0r}$ and $m_{1r}$ are unknown. However, for $n$ large enough they might be accurately approximated by $S_{0r}$ and  $S_{1r}$, 
	and substituting them we obtain Algorithm~\ref{Alg:myriad_general_fast}. 
	
	\begin{algorithm}[!ht]
		\caption{Minimization of $L(\cdot,\cdot)$ (Fast Generalized Myriad Filter, fast GMF) }\label{Alg:myriad_general_fast}
		\begin{algorithmic}
			\State \textbf{Input:} $x_1< \ldots < x_n$, $n \ge 3$, $0< w_i < \frac12$, $i=1,\ldots,n$, $\sum\limits_{i=1}^n w_i = 1$,
			\State \textbf{Initialization:} $a_0 \in (x_1,x_n)$, $\gamma_0 >0$
			\For{$r=0,\ldots$}
			\begin{align}
			a_{r+1}& = {a}_r + {\gamma}_r\frac{S_{1}(a_r,\gamma_r)}{S_{0}^2(a_r,\gamma_r) + S_{1}^2(a_r,\gamma_r)},\\
			\gamma_{r+1} & =  {\gamma}_r\left(\frac{S_{0}(a_r,\gamma_r)}{S_{0}^2(a_r,\gamma_r) + S_{1}^2(a_r,\gamma_r)}-1\right).\\
			\end{align} 	
			\EndFor
		\end{algorithmic}
	\end{algorithm}
	
	Interestingly,  this algorithm converges without any further assumption on the sample size~$n$, as the following theorem shows.
	
	\begin{theorem}\label{Theo:alg_fast}
		Let $n \geq 3$ and $w_i < \frac12$, $i=1,\ldots,n$.
		Then, for any starting point $a_0\in (x_1,x_n)$,  $\gamma_0 >0$, the sequence $\{ (a_r,\gamma_r) \}_{r\in \N}$ 	
		generated by  Algorithm \ref{Alg:myriad_general_fast}
		converges to the minimizer
		$(\hat{a},\hat{\gamma})$ of $L$.		
	\end{theorem}

	\begin{remark}
		\begin{enumerate}[\normalfont(i)]
			\item If $m_{0r}$ and $m_{1r}$ were known, Algorithm~\ref{Alg:myriad_general_fast} would converge in one step. 
			In practice we observed that the Algorithm converges very fast, in particular for large sample sizes for which the approximation of $m_{0r}$ and $m_{1r}$ by $S_{0r}$ and $S_{1r}$ is good. 
			Also that we have equality $\Upsilon=1$ in the proof of Theorem~\ref{Theo:alg_fast} indicates that the algorithm exhausts the possible step size in each iteration.  
			\item Of course, also for Algorithm~\ref{Alg:myriad_general_fast} one might consider the cases when one of the parameters is known and fix it in the iterations, 
			and our experiments revealed that the resulting algorithms are faster than their counterparts Algorithm~\ref{Alg:myriad_a} and Algorithm~\ref{Alg:myriad_gamma}. 
			However, convergence cannot be guaranteed in these cases; in fact, for small sample sizes it is possible that the algorithm for $a$ ($\gamma$ fixed) starts cycling, while $\gamma$ ($a$ fixed) might converge to zero. 
			This might be explained by the coupling between $a_r$ and $\gamma_r$ introduced in the update step. 
		\end{enumerate}
	\end{remark}
	
	\subsection{Comparison of the GMFs in Algorithm \ref{Alg:myriad_general} and \ref{Alg:myriad_general_fast}}\label{Sec:Simulation}
	In order to evaluate the numerical performance, in particular the speed of convergence, of the two proposed algorithms 
	we did the following Monte Carlo simulation: an i.i.d.\ sample of size $n$ from a $C(a,\gamma)$ distribution is drawn 
	and the Algorithms~\ref{Alg:myriad_general} and~\ref{Alg:myriad_general_fast} are run to compute the joint ML-estimate $(\hat{a},\hat{\gamma})$. Both algorithms are initialized with the median of the samples for $a$ and the Hodge-Lehmann estimator~\cite{KP12} for $\gamma$ and we used tolerance $\frac{\lVert (a_{r+1},\gamma_{r+1}) - (a_{r},\gamma_{r})\rVert}{\lVert  (a_{r},\gamma_{r})\rVert} < 10^{-6}$ as stopping criterion. This experiment is repeated $N=10000$ times and afterwards, we calculated the average number of iterations  $\overline{\text{iter}}_1$ and $\overline{\text{iter}}_4$ needed to reach the tolerance criterion together with their standard deviations, the averages  $\bar{{a}}_N $ and $\bar{{\gamma}}_N $ of the values $\hat{a}$ and $\hat{\gamma}$ and their standard deviations $\sigma(\bar{{a}}_N) $ and $ \sigma(\bar{{\gamma}}_N) $, i.e.
	\begin{align*}
	\bar{{\theta}}_N &= \frac{1}{N}\sum_{k=1}^N \hat{\theta}_k, \\
	\sigma(\bar{\theta}_N) &= \sqrt{\frac{1}{N-1} \sum_{k=1}^N (\hat{\theta}_k-\bar{\theta}_N)^2}, \qquad \hat \theta\in \{\hat a,\hat\gamma\}.
	\end{align*}     
	where $\hat{\theta}_k$ denotes the obtained estimate in the $k$-th experiment. Further, we computed the mean squared error
	\begin{equation*}
	\operatorname{MSE}(\theta) = \frac{1}{N} \sum_{k=1}^{N}(\hat{\theta}_k- \theta)^2.
	\end{equation*}
	The results are given in Table~\ref{Tab:sampling_experiment}, where we chose $n\in \{10,50,100\}$, $a=0$ and $\gamma\in\{0.1,1,5,10\}$. 
	First, we notice that the average number of iterations is approximately three times higher for Algorithm~\ref{Alg:myriad_general}, and further, it does merely not depend on $(a,\gamma)$, but rather only on the number of samples $n$. Here, larger sample sizes result in less iterations. As to be expected, the estimated parameters $\hat{a}$ and $\hat{\gamma}$ become more and more accurate for increasing sample size $n$, and their standard deviations are of the same order. The standard deviations and also the MSE become higher for larger values of $\gamma$, which is reasonable since the samples inherit a greater variability in this case. The results are qualitatively similar in case of estimating only one parameter while the other one is known and fixed.

	\begin{table*}[]
		\centering
		\footnotesize
		\caption{Comparison of Algorithm~\ref{Alg:myriad_general} and Algorithm~\ref{Alg:myriad_general_fast}.}
		\label{Tab:sampling_experiment}
		\begin{tabular}{cl|llllll}
			$\gamma$&	$n$	        & $\overline{\text{iter}}_1\pm \sigma(\overline{\text{iter}}_1)$&  $\overline{\text{iter}}_4\pm \sigma(\overline{\text{iter}}_4)$          &$\bar{\hat{a}}_N \pm \sigma(\bar{\hat{a}}_N) $   & $\MSE(\hat{a}_N)$ & $\bar{\hat{\gamma}}_N \pm \sigma(\bar{\hat{\gamma}}_N) $   
			& $\MSE(\hat{\gamma}_N)$ \\ \cline{1-8}
			&	10	& $26.5283\pm 9.8730$         & $11.5150\pm 5.7998$ &  $ -0.006\pm 0.0531$   &$ 0.0028$ & $0.0999\pm 0.0531$ & $0.0028$   \\
			$0.1$ &	50	& $18.6562\pm 2.5454$  & $6.7790\pm 1.4836$  &  $-0.0002\pm 0.0207$ &$0.0004$ & $0.1001\pm 0.0207 $ & $0.0004 $ \\
			&	100	& $ 17.1849\pm 1.9423  $  &   $5.8667\pm 1.0870$                   &   $0.0002\pm 0.0143$  &$0.0002$ & $0.1000\pm 0.0143$ & $0.0002$  \\ \cline{1-8}
			&	10	&   $26.6976\pm 9.5614$ & $11.6328\pm 5.6372$  & $0.0016\pm 0.5272$  & $0.2779$ & $1.0010\pm 0.5435$ & $0.2953$ \\
			$ 1$&	50	& $18.6468\pm 2.5664$  & $6.7959\pm 1.4903$                       & $0.0001\pm 0.2061$ & $0.0425$ & $0.9981\pm 0.2063$ & $0.0426$ \\
			&	100	&   $17.1643\pm 1.9431$ & $5.8671\pm 1.0816$ & $ -0.0006\pm 0.1432$ & $0.0205$ & $1.0012\pm 0.1460$ & $0.0213$ \\ \cline{1-8}
			&	10	&  $26.5601\pm 9.6086$  & $11.5128\pm 5.9935$  & $0.0242\pm 2.6473$   & $7.0082$ & $5.0442\pm 2.7048$ & $7.3173$ \\
			$ 5$&	50	&  $18.6685\pm 2.5569$&  $6.7773\pm 1.4976$ & $0.0136\pm 1.0295$  & $1.0599$   & $5.0005\pm 1.0252$ &  $1.0508$  \\
			&	100	&  $17.1735\pm 1.9372$  & $5.8558\pm 1.0774$   & $0.0075\pm 0.7087$  & $0.5023$ & $5.0001\pm 0.7120$ & $0.5068$ \\ \cline{1-8}
			&	10	&   $26.6397\pm 10.8184$  & $11.6081\pm 6.2560$   & $-0.0763\pm 5.2767$   & $27.8461$ & $10.0181\pm 5.4761$ & $29.9846$ \\
			$ 10$&	50	&  $18.6471\pm 2.5439$  & $6.8004\pm 1.4772$  & $0.0304\pm 2.0516$ & $4.2095$ & $10.0036\pm 2.0704$ & $4.2861$ \\
			&	100	& $17.1711\pm 1.9391$ & $5.8545\pm 1.0802$  & $-0.0130\pm 1.4356$ & $2.0610$ & $10.0000\pm 1.4452$  & $2.0885$ 
		\end{tabular}
	\end{table*}
	
	\begin{remark}
		Note that we chose a very strong convergence criterion. 
		In most applications, already a tolerance of $10^{-3}$ yields satisfying results, 
		which is achieved after 8-10 iterations of Algorithm~\ref{Alg:myriad_general}, and 3-4 iterations of Algorithm~\ref{Alg:myriad_general_fast}.
	\end{remark}

	\section{Myriad Filtering and Image Denoising}\label{sec:myr}
	
	In this section we describe how the developed MF algorithms can be used to denoise images corrupted by additive Cauchy noise. To this aim, let $f\colon \GG\rightarrow \R$ be a noisy image, where $\GG = \{1,\ldots,n_1\}\times \{1,\ldots,n_2\}$ 
	denotes the image domain. 
	We assume that each pixel $i=(i_1,i_2)\in \GG$ is affected by the noise in an independent and identical way, so that we can model the image as 
	\begin{equation}
	f = u +\gamma \eta,\qquad \eta\sim C(0,1),\quad \gamma>0,
	\end{equation} 
	where $u$ is the noise-free image we wish to reconstruct and the scale parameter $\gamma>0$ determines the noise level. Together with the properties of the Cauchy distribution, see Proposition~\ref{Prop:properties_Cauchy}, this results in independent realizations $f_i$ of $C(u_i,\gamma)$ random variables, $i\in \GG$. Now, for each $i\in \GG$ we wish to estimate the underlying $u_i$ using a myriad filtering approach. 
	
	\subsection{Local and Nonlocal Filtering}
	The estimation of the noise-free image requires to select for each $i\in \GG$ a set of indices of samples $\SS(i)$ that are interpreted as i.i.d.\ realizations of $C(u_i,\gamma)$. 
	The strategies used to determine the set $\SS(i)$ can roughly be divided into two different approaches, namely local and nonlocal ones. The local approach assumes that the image does not significantly change in a small neighborhood of a pixel $i\in \GG$ and thus takes the indices of this local neighborhood as set of samples. In case of a squared $r\times r$ neighborhood around $i\in \GG$, this results in 
	\begin{equation*}
	\SS(i)= \bigl\{j\in \GG\colon \lVert i-j\rVert_\infty = \max\{|i_1-j_1|, |i_2-j_2|\}\leq r\bigr\}.
	\end{equation*}
	Here and in all subsequent cases we extend the image at the boundary using mirror boundary conditions. 
	The parameter $r\in \N$ determining the size of the neighborhood has to maintain the following trade off: On the one hand, it has to be sufficiently large to guarantee an appropriate sample size, while on the other hand, for a too large neighborhood the local similarity assumption is unlikely to be fulfilled. 
	
	The nonlocal approach is based on an image self-similarity assumption stating that small patches of an image can be found several  times in the image. Then, the set $\SS(i)$ is constituted of the indices of the centers of patches that are similar to the patch centered at $i\in \GG$. This approach requires the selection of the patch size and an appropriate similarity measure. As we detail later on, both need to be adapted to the noise statistic and the noise level. Based on the similarity measure, one possibility is to take as the set $\SS(i)$ the indices of the centers of the $K$ most similar patches. In order to avoid a computational overload one typically restricts the search zone for similar patches to a $w\times w$ search window around $i\in \GG$. 
	
	From a statistical point of view, the nonlocal approach is more reasonable than the local one. Although the selection of similar patches still introduces a bias, the resulting samples are closer to the i.i.d.\ assumption than in the local case, in particular in image regions with high contrast and sharp edges. 
	
	Having defined for each $i\in \GG$ a set of indices of samples $\SS(i)$, the noise-free image can be estimated as
	\begin{align*}
	\hat{u} &\in \argmin_{u,\gamma}	\sum_{i\in \GG} L\left(u_i,\gamma_i|\{f_j\}_{j\in \SS(i)}\right)\\
	&= \argmin_{u,\gamma}	\sum_{i\in \GG} \sum_{j \in {\SS}(i)} \log\left( (f_j - u_i)^2 + \gamma_i^2 \right) - \log(\gamma_i).
	\end{align*}
	At this point, one might either assume that the noise level $\gamma$ is known and constant, i.e.\ $\gamma_i\equiv \gamma$, $i\in\GG$, or that it is unknown and also needs to be estimated. As we will see later on, even if $\gamma$ is known, it is preferable to estimate $\gamma_i$ individually for each $i\in \GG$. This might be explained by the fact that the selection of the samples used in the myriad filter introduces a bias in the estimation. Indeed, in nearly constant or smooth areas, where the image does not vary significantly, both the local and the nonlocal approach will yield very similar samples. However, in regions with sharp edges or complex patterns, the variability of the samples will be much higher and thus, estimating not only $\hat{u}_i$ but also the local noise level $\hat{\gamma}_i$ might compensate these effects. 
	Further, from a theoretical point of view, it makes the minimizer unique, see Theorem~\ref{theo_both}. 
	
	In both cases, the resulting minimization problem  can be solved pixelwise, either using one of the GMF Algorithms~\ref{Alg:myriad_general} or~\ref{Alg:myriad_general_fast}, or the classical MF Algorithm~\ref{Alg:myriad_a}.
	We call the resulting methods nonlocal (generalized) myriad filtering N-(G)MF and local (generalized) myriad filtering L-(G)MF. 
	\begin{remark}[Robustness of Parameters in N-GMF]\label{Rem:robustness}
		In both the local as well as the nonlocal approach the algorithm depends on several parameters, namely on the one hand the size of the local neighborhood and on the other hand, in the nonlocal approach one has to choose the patch width $s$, the size of the search
		window $w$, and the number of similar patches $n$. In the local approach in all our experiments a rather small local neighborhood of size $3\times 3$ gave the best results. Concerning the nonlocal method, an extensive grid search revealed that for all experiments, nearly the same parameters
		were optimal (w.r.t.\ PSNR), so we set them as follows: the sample size is $n = 40$, the patch width is $s = 3$ for $\gamma = 5$ and $s=5$ for $\gamma = 10$ and the search window size is $w = 31$. Further, we used uniform weights $w_i = \frac{1}{n}$, $i = 1,\ldots,n$. Additionally, we examined choosing weights based on the similarity of patches, more details can be found in our last numerical experiment.
	\end{remark}
	
	%
	%
	
	\subsection{Estimation of the Noise Level $\gamma$}
	If the overall noise level determined by the scale parameter $\gamma>0$ is unknown 
	and should not be estimated locally during the myriad filtering, it can be 
	estimated in constant areas of an image where the signal to noise ratio is weak and differences between pixel values 
	are solely caused by the noise. 
	%
	%
	
	
	%
	In order to detect constant regions we proceed as follows: First, the image grid $\GG$ is partitioned into $K$ small, 
	non-overlapping regions $\GG= \bigcup_{k=1}^K R_k$, and for each region we consider the hypothesis testing problem 
	\begin{align}
	\HH_0&\colon R_k\text{ is constant}\qquad \text{vs.}\qquad 
	\HH_1\colon R_k\text{ is not constant}	\label{constant_test}.
	\end{align}
	
\begin{figure*} [thb]
		\centering
		\begin{subfigure}[t]{0.32\textwidth}
			\centering	
			\includegraphics[width = .98\textwidth]{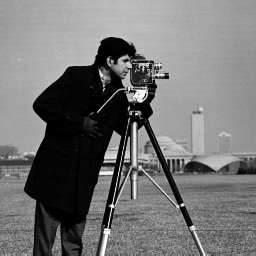}
		\end{subfigure}	
		\begin{subfigure}[t]{0.32\textwidth}
			\centering
			\includegraphics[width = .98\textwidth]{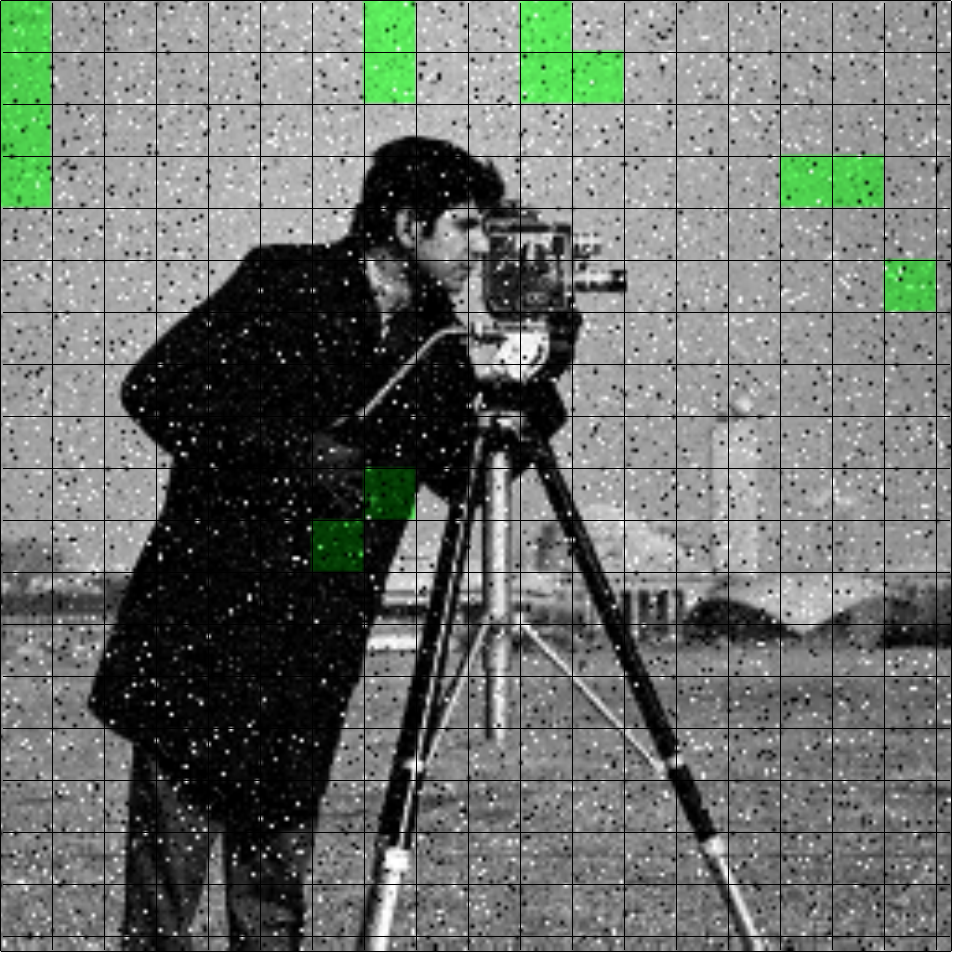}
		\end{subfigure}	
		\begin{subfigure}[t]{0.32\textwidth}
			\centering
			\includegraphics[width = .98\textwidth]{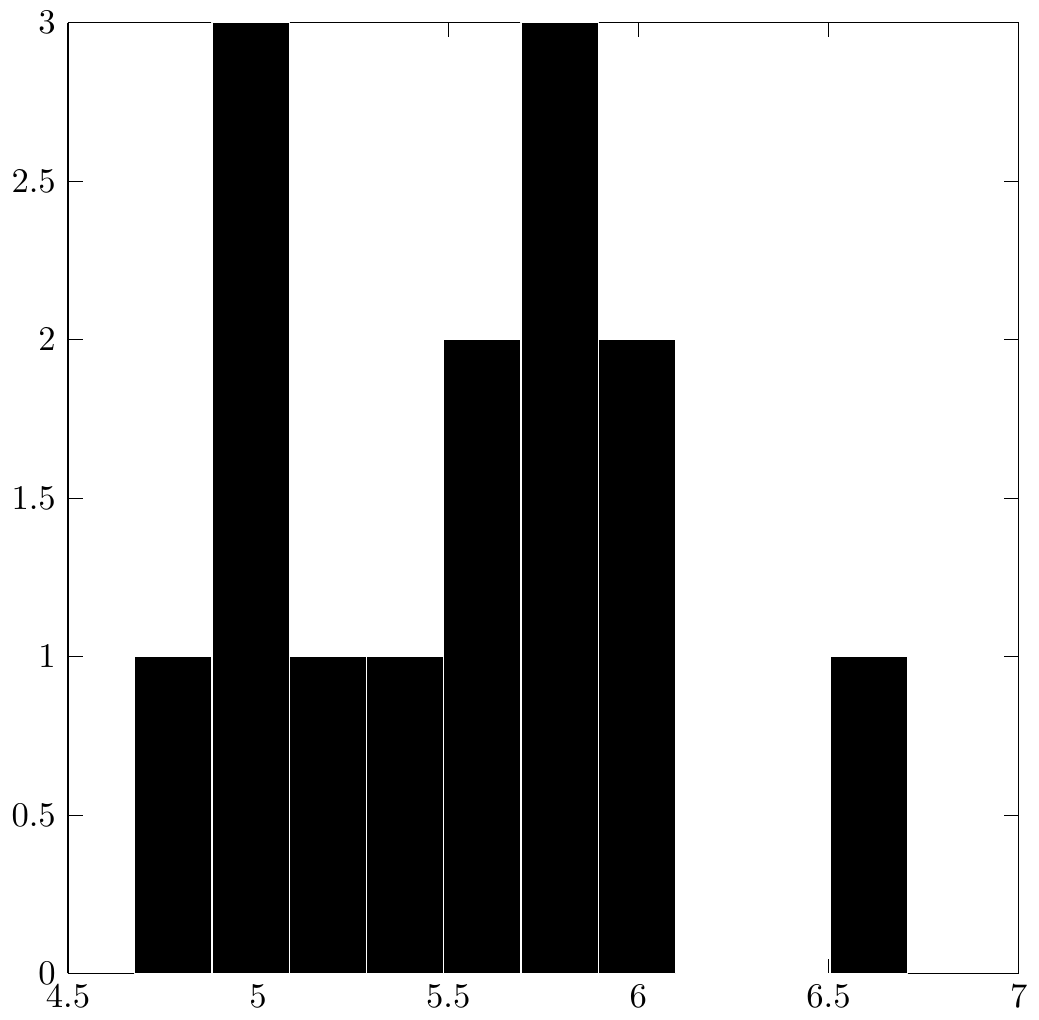}
		\end{subfigure}		
		\caption{Illustration of the estimation of the overall noise level $\gamma$ in constant areas.}\label{Fig:constant_area}
	\end{figure*}

	To decide whether to reject $\HH_0$ or not, we observe the following: Consider a fixed region $R_k$ and let $I, J\subseteq R_k$ be two disjoint subsets of $R_k$ with the same cardinality. Denote with $u_I$ and $u_J$ the vectors containing the values of $u$ at the positions indexed by $I$ and $J$. Then, under $\HH_0$, the vectors $u_I$ and $u_J$ are uncorrelated (in fact even independent) for all choices of $I, J\subseteq R_k$ with $I\cap J = \emptyset$ and $|I|=|J|$. As a consequence, the rejection of $\HH_0$ can be reformulated as the question whether we can find  $I,J$ such that $u_I$ and $u_J$ are significantly correlated, since in this case there has to be some structure in the image region $R_k$ and it cannot be constant. Now, a naive idea would be to use some test statistics based on empirical correlations between pixel values. However, this is not meaningful in case of the Cauchy distribution: Since the distribution does not have any finite moments, the variance of such an estimator will increase with increasing sample size, which is clearly undesirable. 
	As a remedy, we adopt an idea presented in~\cite{SDA15} and make use of Kendall's $\tau$-coefficient, which is a measure of rank correlation, and the associated $z$-score, see~\cite{K38,K45}. The key
	idea is to focus on the rank (i.e., on the relative order) of the values
	rather than on the values themselves. In this vein, a block is considered homogeneous if
	the ranking of the pixel values is uniformly distributed, regardless of the spatial arrangement of the pixels. For convenience, we briefly recall the main definitions and results. In the following, we assume that we have extracted two disjoint subsequences $x = u_I$ and $y = u_J$ from a region $R_k$ with $I$ and $J$ as above.
	Let $(x_i,y_i)$ and $(x_j,y_j)$ be two pairs of observations. Then, the pairs are said to be 
	\begin{equation*}
	\begin{cases}
	\text{concordant} & \text{if } x_i<x_j \text{ and } y_i<y_j\\& \text{or } 
	x_i>x_j \text{ and } y_i>y_j,\\
	\text{discordant} & \text{if } x_i<x_j \text{ and } y_i>y_j\\& \text{or } x_i>x_j \text{ and } y_i<y_j,\\
	\text{tied} & \text{if } x_i=x_j  \text{ or } y_i=y_j.
	\end{cases}
	\end{equation*}	 
	Next, let $x,y\in \R^n$ be two sequences without tied pairs and let $n_c$ and $n_d$ be the number of concordant and discordant pairs, respectively. Then, \emph{Kendall's $\tau$ coefficient}~\cite{Ken38}	is defined as $\tau\colon \R^n\times \R^n\to [-1,1]$, 
	\begin{equation*}
	\tau(x,y) = \frac{n_c - n_d}{\frac{n(n-1)}{2}}. 
	\end{equation*}
	From this definition we see that if the agreement between the two rankings is perfect, i.e.\ the two rankings are the same, then the coefficient attains its maximal value 1. On the other extreme, if the disagreement between the two rankings is perfect, that is, one ranking is the reverse of the other, then the coefficient has value -1.
	If the sequences $x$ and $y$ are uncorrelated,  we expect the coefficient to be approximately zero. Denoting with $X$ and $Y$ the underlying random variables that generated the sequences $x$ and $y$ we have the following result, whose proof can be found in~\cite{K38}. 
	\begin{theorem}\label{Theo:tau_asymptotic}
		Let $X$ and $Y$ be two arbitrary sequences under $\HH_0$ without tied pairs.
		Then, 	the random variable $\tau(X,Y)$ has an expected value of 0 and a variance of $\frac{2(2n+5)}{9n(n-1)}$. Moreover, for $n\to \infty$, the associated \emph{$z$-score} $z\colon \R^n\times \R^n\to \R$,
		\begin{align*}
		z(x,y) = \frac{3\sqrt{n(n-1)}}{\sqrt{2(2n+5)}}\tau(x,y)=\frac{3\sqrt{2}(n_c - n_d)}{\sqrt{n(n-1)(2n+5)}} 
		\end{align*}
		is asymptotically standard normal distributed,
		\begin{equation*}
		z(X,Y)\overset{n\to \infty}{\sim}\NN(0,1).
		\end{equation*}
	\end{theorem}
	With slight adaption, Kendall's $\tau$ coefficient can be generalized to sequences with tied pairs, see~\cite{Ken45}. As a consequence of Theorem~\ref{Theo:tau_asymptotic}, for a given significance level $\alpha\in (0,1)$, we can use the quantiles of the standard normal distribution to decide whether to reject $\HH_0$ or not. In practice, we cannot test any kind of region and any kind of disjoint sequences. As in~\cite{SDA15}, we restrict our attention to quadratic regions and pairwise comparisons of neighboring pixels. We use four kinds of neighboring relations (horizontal, vertical and two diagonal neighbors) thus perform in total four tests. We reject the hypothesis $\HH_0$ that the region is constant as soon as one of the four tests rejects it. Note that by doing so, the final significance level is smaller than the initially chosen one. We start with blocks of size $16\times 16$  whose side-length is incrementally decreased until enough constant areas are found. Then, in each constant region we use Algorithm~\ref{Alg:myriad_general} or Algorithm~\ref{Alg:myriad_general_fast} to estimate the parameters of the associated Cauchy distribution. The estimated location parameters  $a$ of the found regions are discarded, while the estimated scale parameters $\gamma$ are averaged to obtain the final estimate of the global noise level. Figure~\ref{Fig:constant_area} illustrates this procedure by means  of the \emph{cameraman} image. As reference, the original image is shown in the left, while the detected constant areas in the noisy image are depicted in the middle. On the right we give a histogram of the estimated values for $\gamma$. The final estimate for this example is $\hat{\gamma}=5.5283$, while the true parameter used to generate the noisy image was $\gamma = 5$.


	Finally, we would like to point out that our motivation for using this rather complicated way to determine the correlation differs from those of the authors of~\cite{SDA15}: In their work, the authors aim at determining the  noise characteristics, which is assumed to be unknown, so that they need an non-parametric test statistic that is independent of the underlying noise distribution. Therefore, simpler methods such as for instance empirical correlation functions cannot be applied. In our situation, however, the noise statistic is known. Nevertheless we cannot directly estimate the correlation, since any approach involving the empirical correlation of the data is no meaningful in case of the Cauchy distribution.

	\subsection{Selection of Similar Patches in Cauchy Noise Corrupted Images }\label{sec:patch}
	
	The selection of similar patches constitutes  a fundamental step in our nonlocal denoising approach.
	At this point, the question arises how to compare noisy patches and numerical 
	examples show that an adaptation of the similarity measure to the noise
	distribution is essential for a robust similarity evaluation.
	In~\cite{DDT12}, the authors formulated the similarity between  patches as a statistical hypothesis testing problem. Further, they collected and discussed several similarity measures, 
	among others see~\cite{DDT09}, for multiplicative noise~\cite{ST10,TL12} and the references therein. While they only considered Gaussian, 
	Poisson and Gamma noise, we extend some of their results to Cauchy noise.
	
	Modeling noisy images in a stochastic way allows to formulate the question whether two patches $p$ and $q$ are similar as a hypothesis test.
	Two noisy patches $p,q$ are   considered to be similar  if they are
	realizations  of  independent  random  variables $X\sim p_{\theta_1}$ and $Y\sim p_{\theta_2}$ 
	that follow the same parametric distribution $p_\theta$, $\theta\in \Theta$ 
	with a common parameter $\theta$ (corresponding to the underlying noise-free patch), i.e.\ $\theta_1 = \theta_2 \equiv \theta$. 
	Therewith, the evaluation of the similarity between noisy patches can be
	formulated as the following hypothesis test:
	\begin{align}
	\HH_0&\colon \theta_1 = \theta_2\qquad \text{vs.}\qquad 
	\HH_1\colon \theta_1 \not = \theta_2.	\label{similarity_test}
	\end{align}
	In this context, a similarity measure $S$ maps a pair of noisy patches $p,q$ to a real value $c\in \R$. The larger this value $c$ is, the more the patches are considered to be similar. 
	
	In the following, we describe how a similarity measure $S$ can be obtained based on a suitable test statistic for the hypothesis testing problem. 
	In general, according to the Neyman-Pearson Theorem, see, e.g.~\cite{CB02}, the optimal test statistic 
	(i.e.\ the one that maximizes the power for any given size $\alpha$) 
	for single-valued hypotheses of the form
	\begin{align*}
	\HH_0&\colon \theta= \theta_0\qquad \text{vs.}\qquad 
	\HH_1\colon \theta  = \theta_1	
	\end{align*} 
	is given by a likelihood ratio test. Note that single-valued testing problems correspond to a disjoint partition 
	of the parameter space of the form $\Theta = \Theta_0\dot{\cup} \Theta_1$, where $\Theta_i = \{\theta_i\}$, $i=0,1$.	
	
	Despite being a very strong theoretical result, the practical relevance of the Neyman-Pearson Theorem is limited due 
	to the fact that $\Theta_0$ and $\Theta_1$ are in most applications not single-valued. 
	Instead, the testing problem is a so called \emph{composite testing} problem, 
	meaning that $\Theta_0$ and/or $\Theta_1$ contain more than one element.
	It can be shown that for composite testing problems there does not exist a uniformly most powerful test. 
	Now, the idea to generalize the Neyman-Pearson test to composite testing problems 
	is to obtain first two candidates (or representatives) $\hat{\theta}_0$ and $\hat{\theta}_1$ of $\Theta_0$ and $\Theta_1$ respectively, e.g.\ by maximum-likelihood estimation, 
	and then to perform a Neyman-Pearson test using the computed candidates $\hat{\theta}_0$ and $\hat{\theta}_1$ 
	in the definition of the test statistic. 
	In case that an ML-estimation is used to determine $\hat{\theta}_0$ and $\hat{\theta}_1$, 
	the resulting test is called \emph{Likelihood Ratio Test} (LR test). Although there are in general no theoretical guarantees concerning the power 
	of LR tests, they usually perform very well in practice if the sample size used to estimate $\hat{\theta}_0$ and $\hat{\theta}_1$ is large enough.
	This  is due to the fact that ML estimators are asymptotically efficient.	Several classical tests, e.g.\ one and two-sided $t$-tests, 
	are either direct LR tests or equivalent to them.
	
	In the sequel we show how the above framework can be applied to our testing problem~\eqref{similarity_test}. 
	First, let $x_1$ and $y_1$ be two single pixels for which we want to test whether they are realizations of the same distribution with unknown common parameter.
	The LR statistic reads as
	\begin{align} \label{FR_stat}
	\lambda(x_1,y_1) 
	&= \frac{\sup\limits_{\theta\in \HH_0}  
		\bigl\{{\cal L}(\theta|x_1,y_1)\bigr\}}{\sup\limits_{\theta }\bigl\{{\cal L}(\theta|x_1,y_1)\bigr\}}\\
	& =\frac{\sup\limits_{\theta}\bigl\{{\cal L}(\theta|x_1){\cal L}(\theta|y_1)\bigr\}}{\sup\limits_{\theta}\bigl\{{\cal L}(\theta|x_1)\bigr\}
		\sup\limits_{\theta}\bigl\{{\cal L}(\theta|y_1)\bigr\}},
	\end{align}
	where in our situation ${\cal L} (\theta|x_1,y_1)$ denotes the likelihood function with respect to $a$ while $\gamma$ is assumed to be known, and the notation $\theta\in \HH_0$ means that the supremum is taken over those parameters $\theta$ fulfilling $\HH_0$.
	We use this statistic as similarity measure, i.e., 
	$S(x_1,y_1) \coloneqq \lambda(x_1,y_1)$.
	More generally, since we assume the noise to affect each pixel in an independent and identical way, 
	the similarity of two patches $p= (x_1,\ldots,x_t)$ and $q= (y_1,\ldots,y_t)$ is obtained as the product of the similarity of its pixels
	\begin{equation*}
	S(p,q) \coloneqq \prod_{i=1}^t S(x_i,y_i).
	\end{equation*}
	
	\begin{lemma}\label{Lem:Cauchy_similarity}
		For the Cauchy distribution, the LR statistics is given by 
		\begin{equation}
		\lambda(x_1,y_1) = \left(\left(\frac{x_1- y_1}{2\gamma}\right)^2 + 1\right)^{-2}.
		\end{equation}
	\end{lemma}
	The proof of Lemma~\ref{Lem:Cauchy_similarity} can be found in Appendix~\ref{app:myr}.
	Observe that this similarity measure requires the knowledge of the overall noise level $\gamma$, which is either assumed to be known or alternatively can be estimated in advance using the method described in the previous paragraph. Using Lemma~\ref{Lem:Cauchy_similarity}, we obtain the similarity measure between two patches $p,q$ 
	\begin{equation*}
	S(p,q) = \prod_{i=1}^{t} \left(\left(\frac{x_i-y_i}{2\gamma}\right)^2 + 1\right)^{-2}.
	\end{equation*}
	In practice, we take the logarithm of $S$ in order to avoid numerical instabilities.

	\section{Numerical Results} \label{sec:numerics}
	
\begin{figure*} [tp]
		\centering
		\begin{subfigure}[t]{0.16\textwidth}
			\centering
			\includegraphics[width = .98\textwidth]{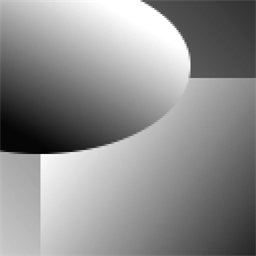}
		\end{subfigure}			
		\begin{subfigure}[t]{0.16\textwidth}
			\centering
			\includegraphics[width = .98\textwidth]{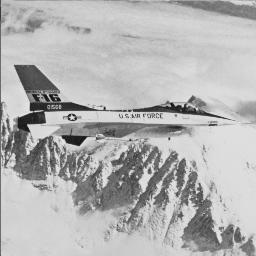}
		\end{subfigure}
		\begin{subfigure}[t]{0.16\textwidth}
			\centering
			\includegraphics[width = .98\textwidth]{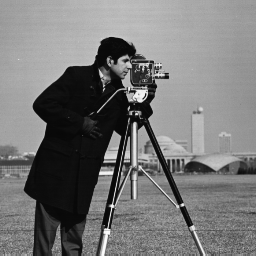}
		\end{subfigure}	
		\begin{subfigure}[t]{0.16\textwidth}
			\centering
			\includegraphics[width = .98\textwidth]{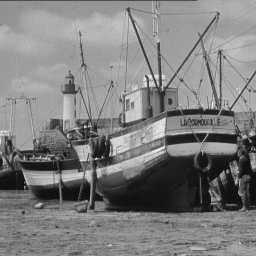}
		\end{subfigure}	
		\begin{subfigure}[t]{0.16\textwidth}
			\centering
			\includegraphics[width = .98\textwidth]{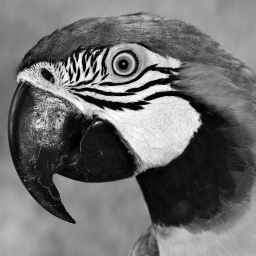}
		\end{subfigure}	
		\begin{subfigure}[t]{0.16\textwidth}
			\centering
			\includegraphics[width = .98\textwidth]{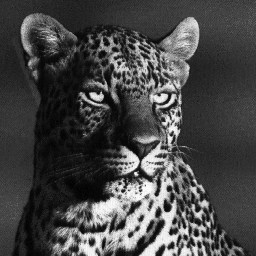}
		\end{subfigure}	
		
		\vspace{0.15cm}
		
		\begin{subfigure}[t]{0.16\textwidth}
			\centering
			\includegraphics[width = .98\textwidth]{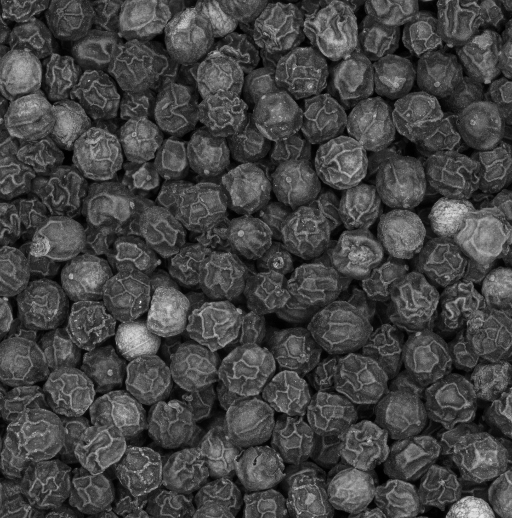}
		\end{subfigure}	
		\begin{subfigure}[t]{0.16\textwidth}
			\centering
			\includegraphics[width = .98\textwidth]{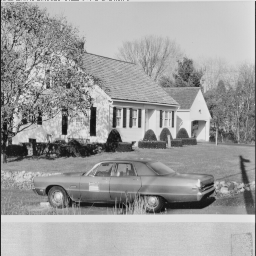}
		\end{subfigure}	
		\begin{subfigure}[t]{0.16\textwidth}
			\centering
			\includegraphics[width = .98\textwidth]{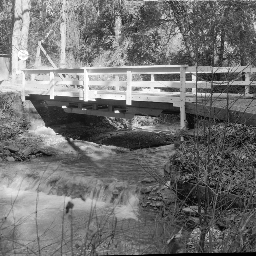}
		\end{subfigure}	
		\begin{subfigure}[t]{0.16\textwidth}
			\centering
			\includegraphics[width = .98\textwidth]{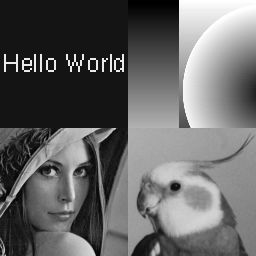}
		\end{subfigure}	
		\begin{subfigure}[t]{0.16\textwidth}
			\centering
			\includegraphics[width = .98\textwidth]{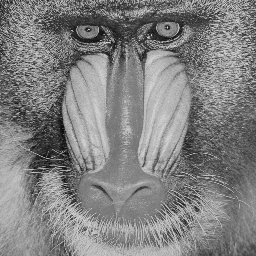}
		\end{subfigure}	
		\begin{subfigure}[t]{0.16\textwidth}
			\centering
			\includegraphics[width = .98\textwidth]{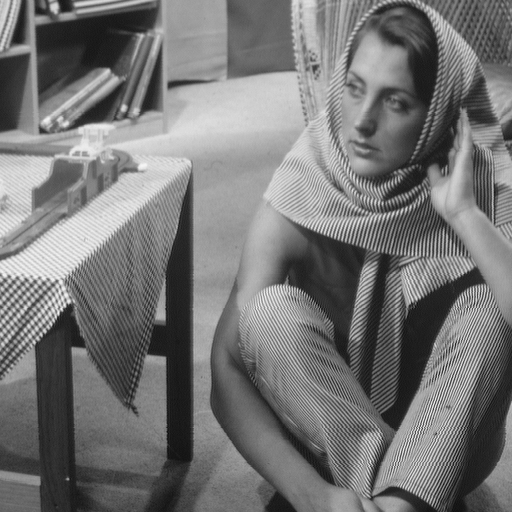}
		\end{subfigure}
		\caption[]{Test images used in our numerical experiments.}
		\label{Fig:test_images}	
	\end{figure*}

	Next, we provide numerical examples to illustrate the different denoising strategies presented in the previous section. 
	All our methods are implemented  in C and imported into Matlab
	using
	\texttt{mex}-interfaces.
	As test images we used the images shown in Figure~\ref{Fig:test_images}, which are standard test images representing the different types and structures one encounters in real world applications: On the one hand, images with sharp edges as well as smooth transitions and constant areas, such as the \emph{test} image or the \emph{cameraman}, and on the other hand images with fine structures and textures of different scale, such as the images \emph{leopard, peppers} or \emph{baboon} (and of course, combinations thereof). 
	
	As quality measures, we used the peak signal-to-noise ration (PSNR)
	\begin{equation}\label{psnr_formula}
	\operatorname{PSNR}(\hat{u},u) = 10\log_{10}\left(\frac{255^2}{\frac{1}{n_1 n_2}\lVert \hat{u}-u\rVert_2^2}\right), 
	\end{equation}
	as well as the structural similarity index (SSIM)~\cite{WBSS04}.
	
\begin{figure*} [tp]
		\centering
		\begin{subfigure}[t]{0.49\textwidth}
			\centering
			\includegraphics[width = .98\textwidth]{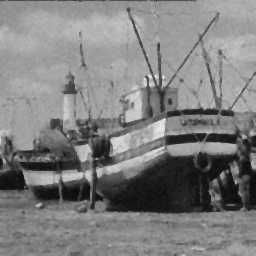}
		\end{subfigure}	
		\begin{subfigure}[t]{0.49\textwidth}
			\centering
			\includegraphics[width = .98\textwidth]{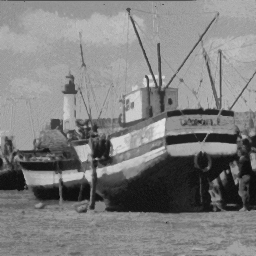}
		\end{subfigure}
		\caption[]{Local (left) versus nonlocal (right) generalized myriad filtering for a noise level of $\gamma = 5$.}
		\label{Fig:LGMF_vs_NGMF}	
	\end{figure*}

	\paragraph{Local vs.\ Nonlocal Generalized Myriad Filtering}
	In a first experiment, we compare the performance of the local versus the nonlocal GMF. 
	The left image in Figure~\ref{Fig:LGMF_vs_NGMF} gives the result of the L-GMF for a $3\times 3$ local neighborhood, while the right image depicts the result obtained using the N-GMF for the \emph{boat} image and a noise level  of $\gamma = 5$. As to be expected, the nonlocal approach yields better results, not only in terms of PSNR and SSIM (27.5307 and 0.7898 versus 28.9941 and 0.8350), but also the visual impression is better as the image is sharper and provides more details. 
	We show here only the results of the generalized versions of the algorithms since they perform in general better than the classical ones, see also our second example. 
	
\begin{figure*} [tp]
		\centering
		\begin{subfigure}[t]{0.49\textwidth}
			\raggedleft
			\includegraphics[width=\textwidth]{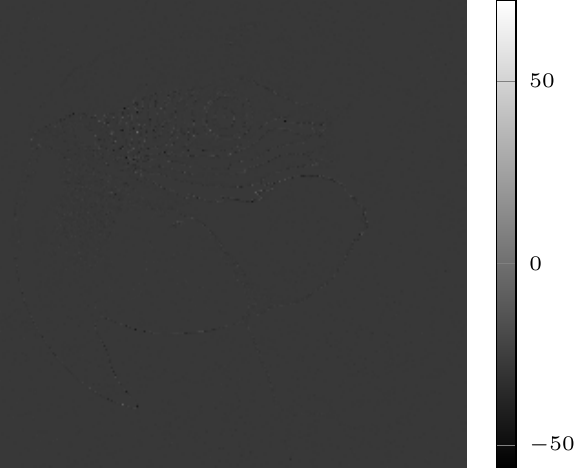}
		\end{subfigure}	
		\begin{subfigure}[t]{0.49\textwidth}
			\raggedright
			\includegraphics[width=.97\textwidth]{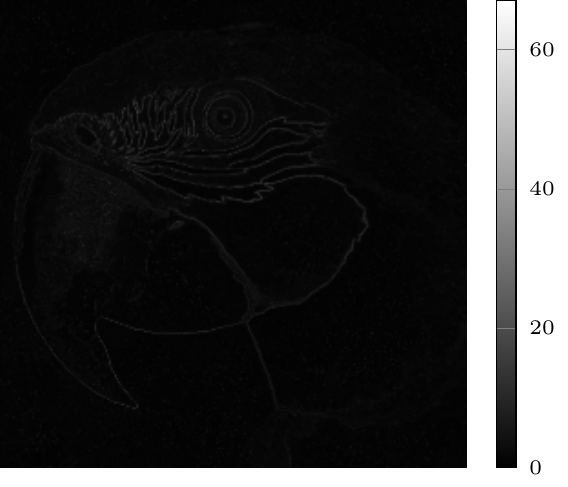}
		\end{subfigure}
		\caption[]{Differences between the generalized and the classical myriad filtering method. Left is shown the difference image between the two filtering methods while the right image shows the locally estimated $\gamma_i$, $i\in \GG$.}
		\label{Fig:GMF_vs_MF}	
	\end{figure*}

	\paragraph{Generalized vs.\ Classical Myriad Filtering}
	Our second experiment illustrates the difference between the generalized and the classical myriad filtering in case of the \emph{parrot} image. 
	We show the results only for the nonlocal approach,  the observations are similar in the local case. 
	Figure~\ref{Fig:GMF_vs_MF} left depicts the difference between the obtained images, while the right image contains in each pixel the locally estimated $\gamma_i$ (the global noise level was $\gamma = 5$). 
	
	While the differences in the images  are only sightly visible at edges, 
	there are significant differences at edges for the locally estimated $\gamma$.
	The difference between the locally estimated and the global value of the scale parameter $\gamma$ 
	can be explained as follows: Simple patches without much structure, for instance in homogeneous 
	or smoothly varying areas,  usually frequently occur in natural images, so that the selection 
	of similar patches results in a set of samples that is highly clustered around the pixel value to be estimated. 
	This naturally decreases the estimated value for $\gamma$. However, for more complicated and structured patches, 
	in particular at sharp edges, it is much more difficult to find similar patches, 
	and the resulting set of samples is very likely to posses a large variability, causing a higher value of $\gamma$. 
	This might be circumvented by allowing a varying number of patches instead of a fixed one, 
	i.e.\ only those whose similarity is larger than a certain threshold. 
	However, by doing so it can happen that one ends up with far to few samples 
	in case of structured patches, leading to wrong results in the estimation procedure. 
	Therefore, the use of the GMF is the better alternative 
	to cope with the bias introduced by the selection of similar patches as samples. 
	Our numerical experiments indeed indicate that one gains on average between 0.5dB and 2dB 
	in terms of PSNR compared to using the MF with fixed scale parameter $\gamma$. 
	This observation is independent of taking the true value of $\gamma$ 
	in the classical myriad filter or not; in fact a wrong $\gamma$ may lead to significantly worse results.

	\paragraph{Comparison with the Variational Method~\cite{MDHY16}}
	In our third example we compare our N-GMF to the variational denoising method 
	for Cauchy noise proposed in~\cite{MDHY16}. 
	Motivated by a \emph{maximum a posteriori} approach with a total variation prior the authors propose a continuous model of the form
	\begin{equation*}
	\mathcal{J}(u) = \DD(u;f) + \lambda \TV(u),\quad \lambda >0,
	\end{equation*}
	where the data term reads as
	\begin{equation*}
	\DD(u;f)=  \int_{\Omega} \log\bigl( (Ku-f)^2 + \gamma^2\bigr ) \dx
	\end{equation*}
	and the regularization term is given by 
	\begin{align*}
	\TV(u) = \int_{\Omega} |Du| = \sup \left\{\int_{\Omega} u(x) \divv\bigl(\phi(x)\bigr) \dx\colon
	 \phi \in C_c^1(\Omega,\R^2), \lVert \phi\rVert_{\infty}\leq 1 \right\}.
	\end{align*} 
	The model  incorporates a linear operator~$K$, which is the identity in the denoising case. 
	In~\cite{SDA15}, an additional quadratic penalty term $\alpha \lVert K u - u_0\rVert_2^2$, where $u_0$ is a median-filtered version of $f$, was added to the data term in order to make the functional convex. Then, the \emph{alternating method of multipliers} (ADMM) from convex optimization can be used to compute a (local) minimizer of the functional. 
	However, it turned out that the nonconvex model performs better than the convexified one, so that we compare our method only to the nonconvex model. For a fair comparison, we used the implementation, the noisy data  
	and the image-wise hand-tuned model- and algorithm-parameters provided by the authors of~\cite{MDHY16}. 
	At this point it is worth mentioning that their results highly depend on the initialization, 
	the parameters and the exact value of $\gamma$, which is not the case in our approach, 
	see Remark~\ref{Rem:robustness}. 
	
	Table~\ref{table_psnr} summarizes the comparison in terms of PSNR and SSIM 
	for two different noise levels $\gamma = 5$ and $\gamma = 10$. 
	Note that the authors of~\cite{MDHY16} used a slightly different definition of PSNR as 
	\begin{equation*}
	\operatorname{PSNR}(\hat{u},u) = 10\log_{10}\left(\frac{(\max{u}-\min{u})^2}{\frac{1}{n_1 n_2}\lVert \hat{u}-u\rVert_2^2}\right), 
	\end{equation*}
	which explains the difference compared to the values stated in~\cite{MDHY16}. 
	Table~\ref{table_psnr} shows that for some images our method yields slightly worse results,  
	while in most cases it performs significantly better in terms of PSNR and SSIM.
	
\begin{figure*} [tp]
		\centering	
		\begin{subfigure}[t]{0.32\textwidth}
			\centering
			\includegraphics[width=.98\textwidth]{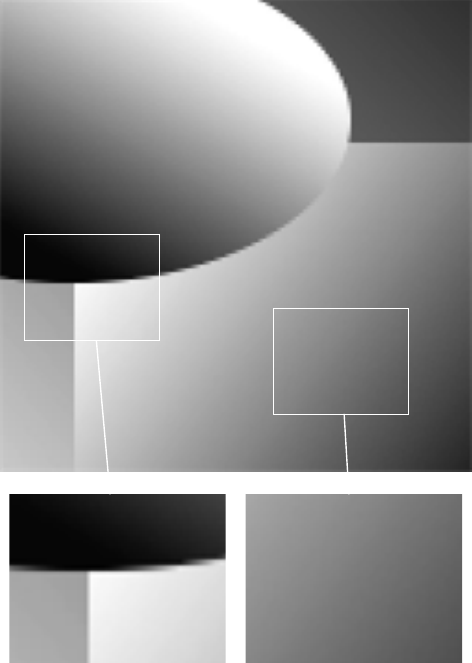}
		\end{subfigure}
		\begin{subfigure}[t]{0.32\textwidth}
			\centering
			\includegraphics[width=.98\textwidth]{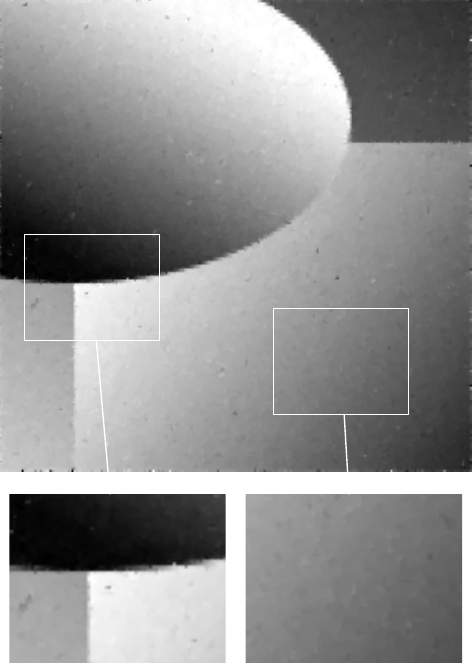}
		\end{subfigure}
		\begin{subfigure}[t]{0.32\textwidth}
			\centering
			\includegraphics[width=.98\textwidth]{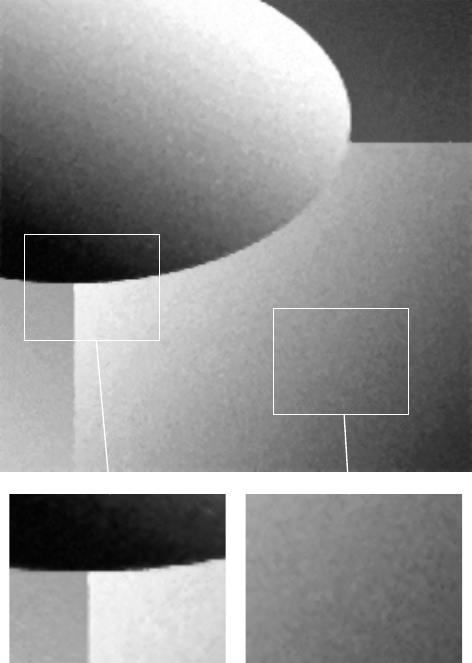}
		\end{subfigure}	
		
		\vspace{0.15cm}
		
		\begin{subfigure}[t]{0.32\textwidth}
			\centering
			\includegraphics[width=.98\textwidth]{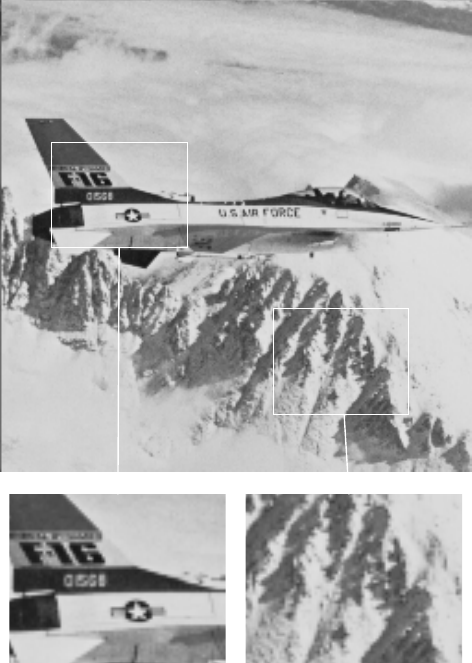}
		\end{subfigure}
		\begin{subfigure}[t]{0.32\textwidth}
			\centering
			\includegraphics[width=.98\textwidth]{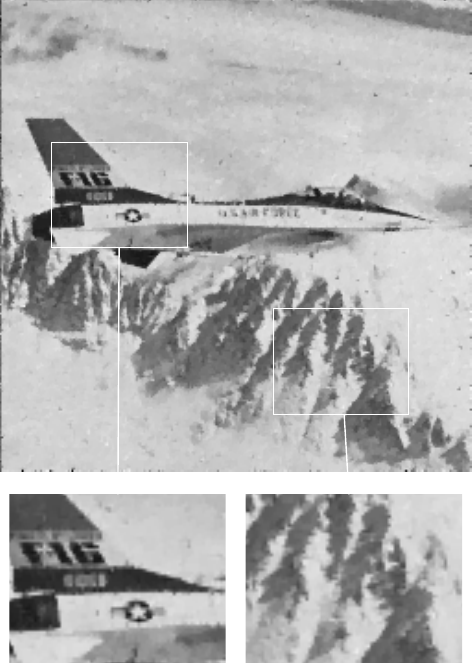}
		\end{subfigure}
		\begin{subfigure}[t]{0.32\textwidth}
			\centering
			\includegraphics[width=.98\textwidth]{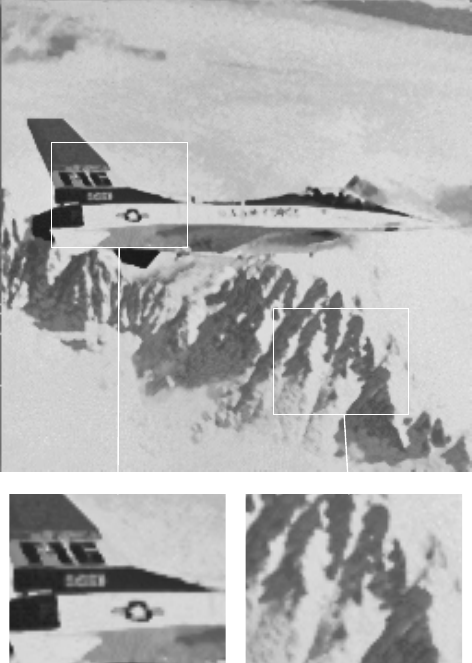}
		\end{subfigure}
		
		\vspace{0.15cm}		
		
		\begin{subfigure}[t]{0.32\textwidth}
			\centering
			\includegraphics[width=.98\textwidth]{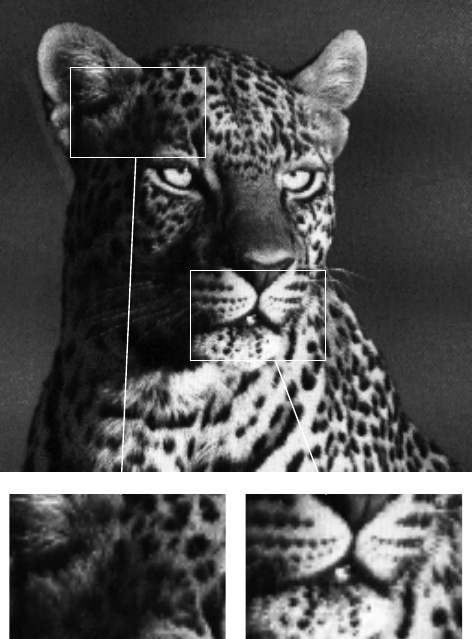}
		\end{subfigure}
		\begin{subfigure}[t]{0.32\textwidth}
			\centering
			\includegraphics[width=.98\textwidth]{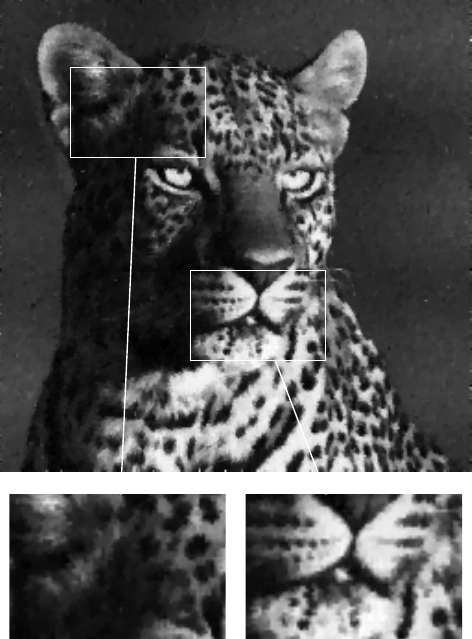}
		\end{subfigure}
		\begin{subfigure}[t]{0.32\textwidth}
			\centering
			\includegraphics[width=.98\textwidth]{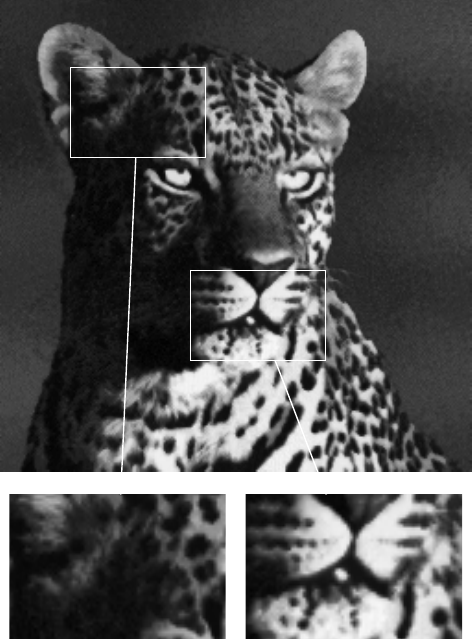}
		\end{subfigure}	
		\caption[]{Denoising of an image (left) for $\gamma =5$ using the method proposed in~\cite{MDHY16} (middle) and our N-GMF (right).}
		\label{Fig:examples_gamma_5}	
	\end{figure*}

\begin{figure*} [tp]
		\centering

		\begin{subfigure}[t]{0.32\textwidth}
			\centering
			\includegraphics[width=.98\textwidth]{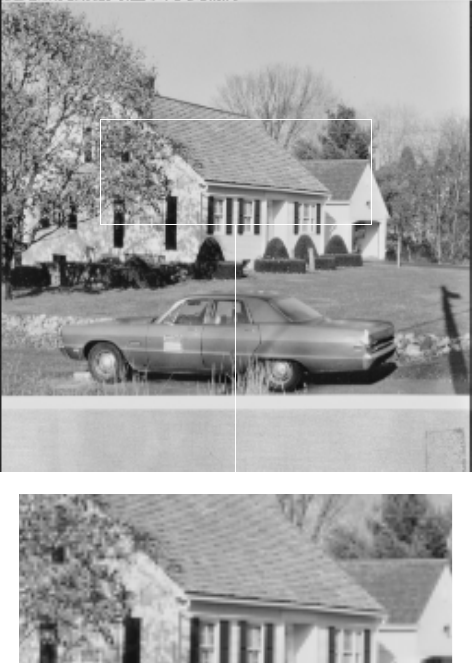}
		\end{subfigure}			
		\begin{subfigure}[t]{0.32\textwidth}
			\centering
			\includegraphics[width=.98\textwidth]{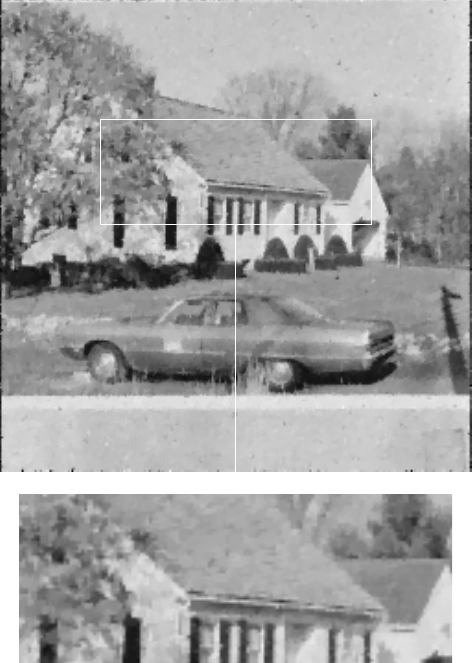}
		\end{subfigure}
		\begin{subfigure}[t]{0.32\textwidth}
			\centering
			\includegraphics[width=.98\textwidth]{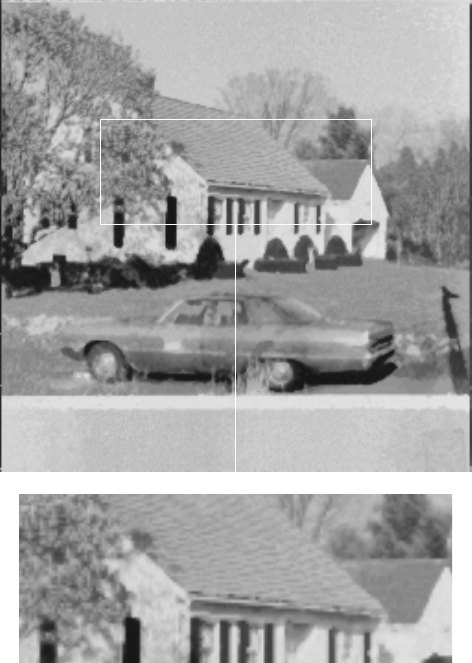}
		\end{subfigure}
		
		\vspace{0.15cm}	
		
		\begin{subfigure}[t]{0.32\textwidth}
			\centering
			\includegraphics[width=.98\textwidth]{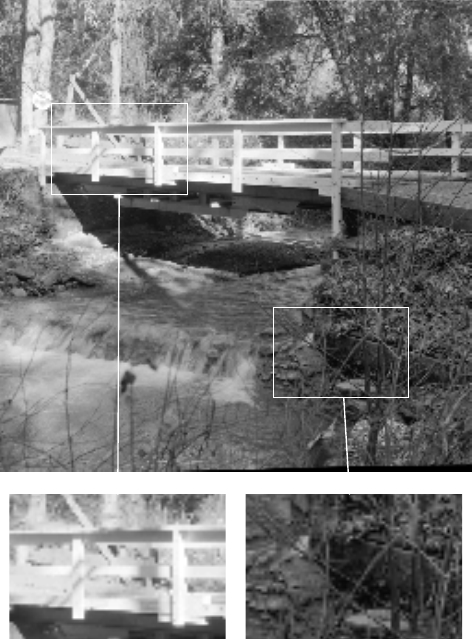}
		\end{subfigure}	
		\begin{subfigure}[t]{0.32\textwidth}
			\centering
			\includegraphics[width=.98\textwidth]{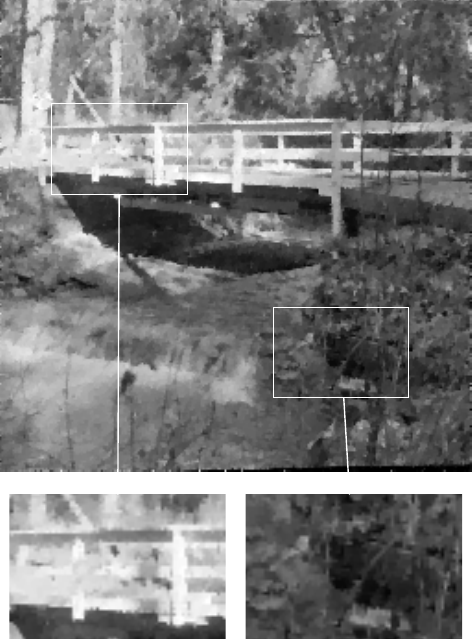}
		\end{subfigure}
		\begin{subfigure}[t]{0.32\textwidth}
			\centering
			\includegraphics[width=.98\textwidth]{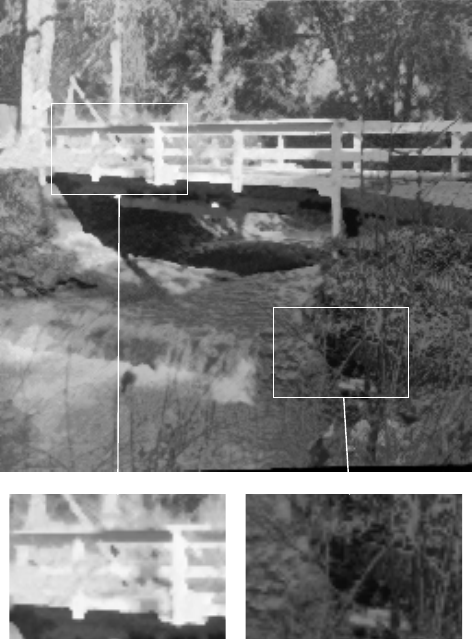}
		\end{subfigure}
		
		\vspace{0.15cm}	
		
		\begin{subfigure}[t]{0.32\textwidth}
			\centering
			\includegraphics[width=.98\textwidth]{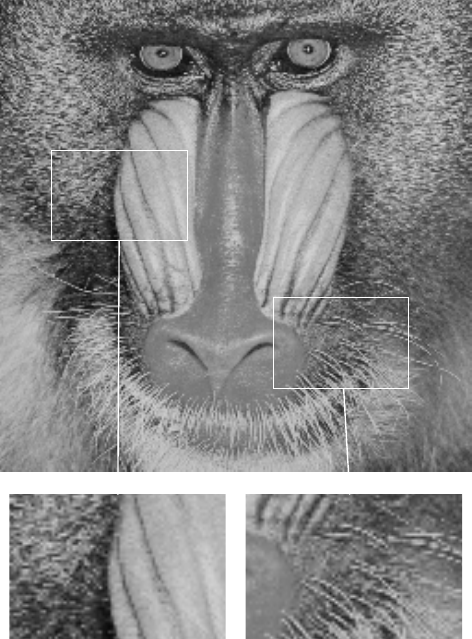}
		\end{subfigure}			
		\begin{subfigure}[t]{0.32\textwidth}
			\centering
			\includegraphics[width=.98\textwidth]{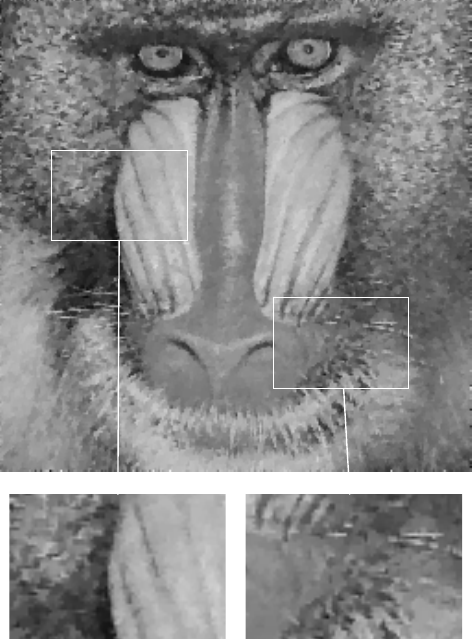}
		\end{subfigure}
		\begin{subfigure}[t]{0.32\textwidth}
			\centering
			\includegraphics[width=.98\textwidth]{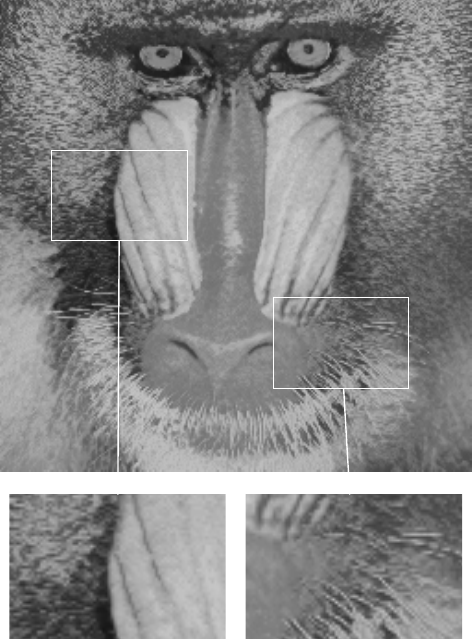}
		\end{subfigure}

		\caption[]{Denoising of an image (left) for $\gamma =5$ using the method proposed in~\cite{MDHY16} (middle) and our N-GMF (right).}
		\label{Fig:examples_gamma_5_2}	
	\end{figure*}

\begin{figure*} [tp]
		\centering
		\begin{subfigure}[t]{0.32\textwidth}
			\centering
			\includegraphics[width=.98\textwidth]{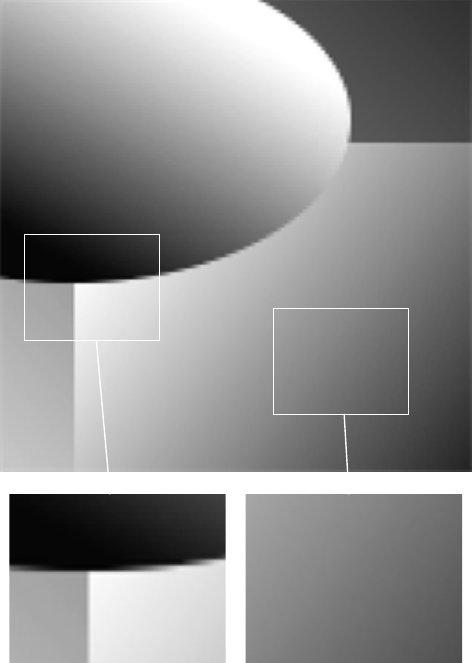}
		\end{subfigure}
		\begin{subfigure}[t]{0.32\textwidth}
			\centering
			\includegraphics[width=.98\textwidth]{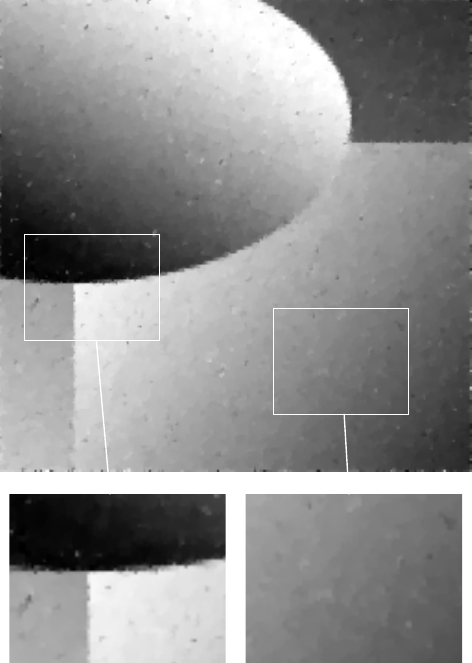}
		\end{subfigure}
		\begin{subfigure}[t]{0.32\textwidth}
			\centering
			\includegraphics[width=.98\textwidth]{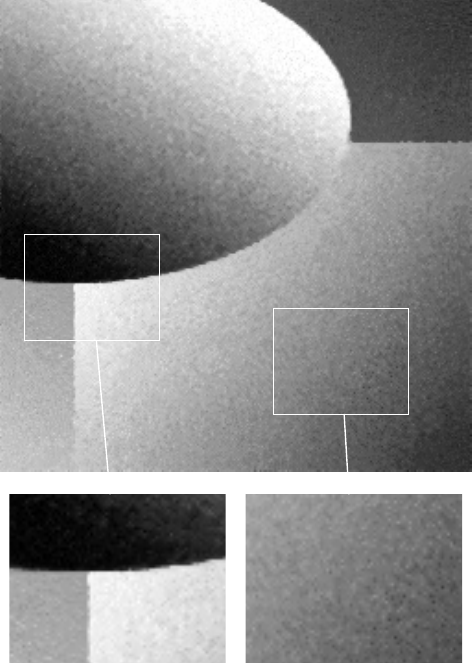}
		\end{subfigure}	
		
		\vspace{0.15cm}	
		
		\begin{subfigure}[t]{0.32\textwidth}
			\centering
			\includegraphics[width=.98\textwidth]{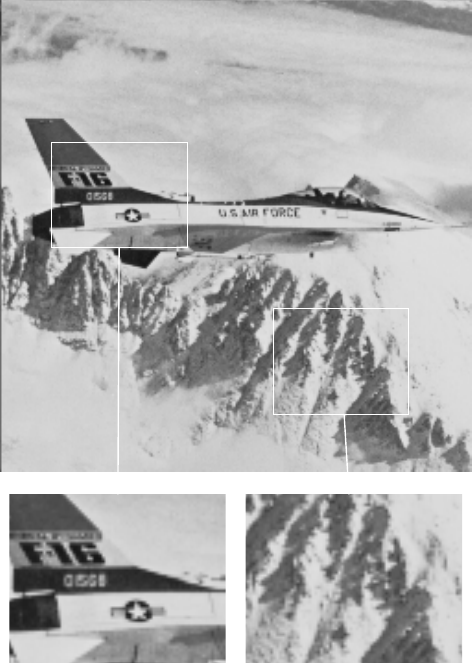}
		\end{subfigure}
		\begin{subfigure}[t]{0.32\textwidth}
			\centering
			\includegraphics[width=.98\textwidth]{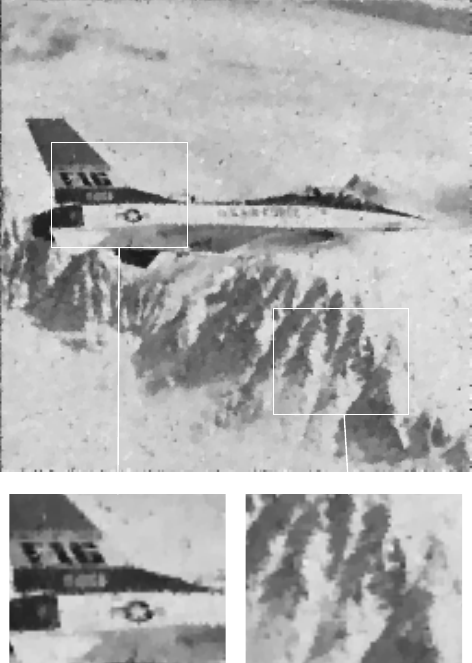}
		\end{subfigure}
		\begin{subfigure}[t]{0.32\textwidth}
			\centering
			\includegraphics[width=.98\textwidth]{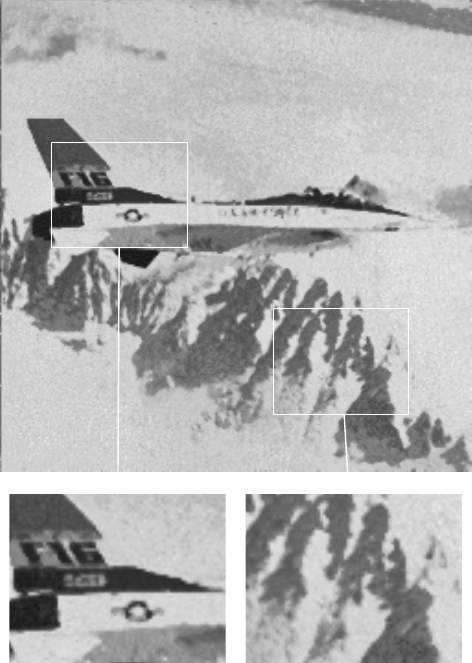}
		\end{subfigure}
		
		\vspace{0.15cm}
		
		\begin{subfigure}[t]{0.32\textwidth}
			\centering
			\includegraphics[width=.98\textwidth]{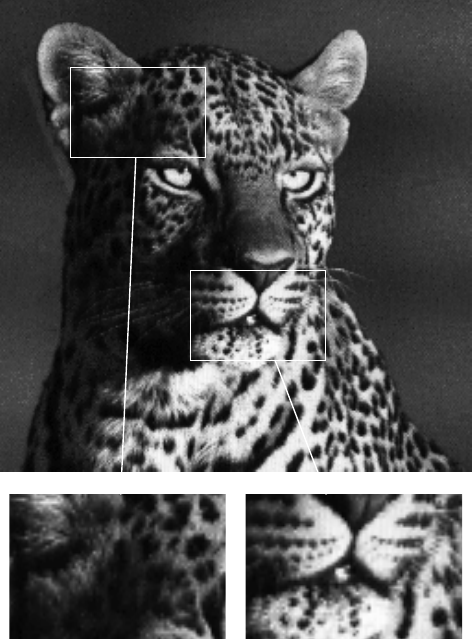}
		\end{subfigure}
		\begin{subfigure}[t]{0.32\textwidth}
			\centering
			\includegraphics[width=.98\textwidth]{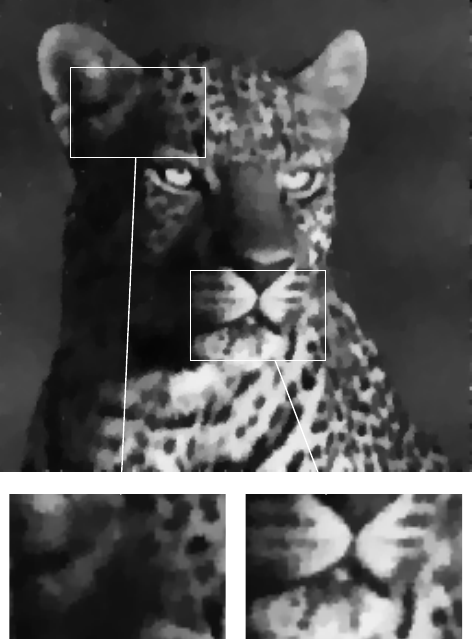}
		\end{subfigure}
		\begin{subfigure}[t]{0.32\textwidth}
			\centering
			\includegraphics[width=.98\textwidth]{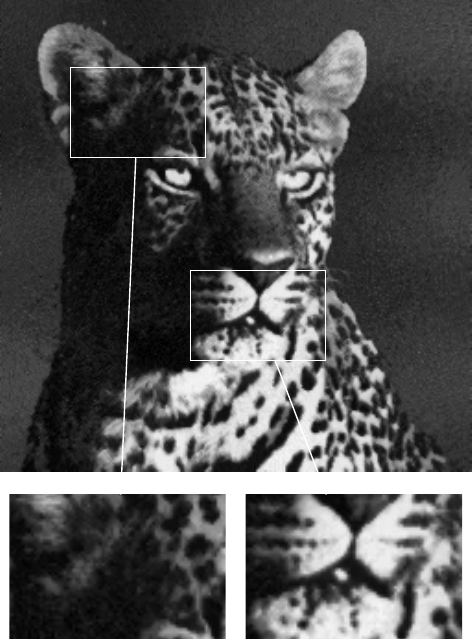}
		\end{subfigure}
		
		\caption[]{Denoising of an image (left) for $\gamma =10$ using the method proposed in~\cite{MDHY16} (middle) and our N-GMF (right).}
		\label{Fig:examples_gamma_10}	
	\end{figure*}

\begin{figure*} [tp]
		\centering

		\begin{subfigure}[t]{0.32\textwidth}
			\centering
			\includegraphics[width=.98\textwidth]{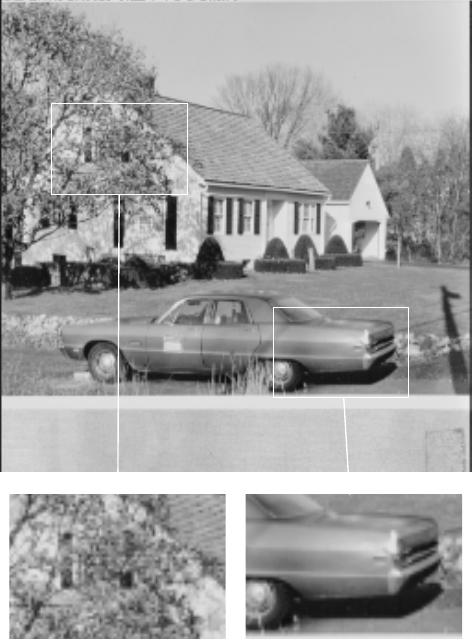}
		\end{subfigure}
		\begin{subfigure}[t]{0.32\textwidth}
			\centering
			\includegraphics[width=.98\textwidth]{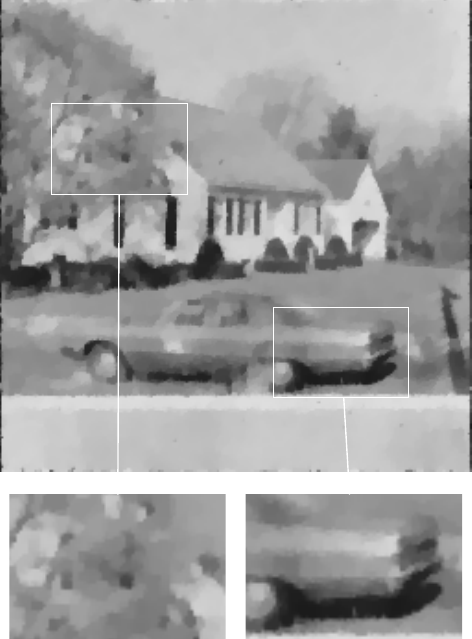}
		\end{subfigure}
		\begin{subfigure}[t]{0.32\textwidth}
			\centering
			\includegraphics[width=.98\textwidth]{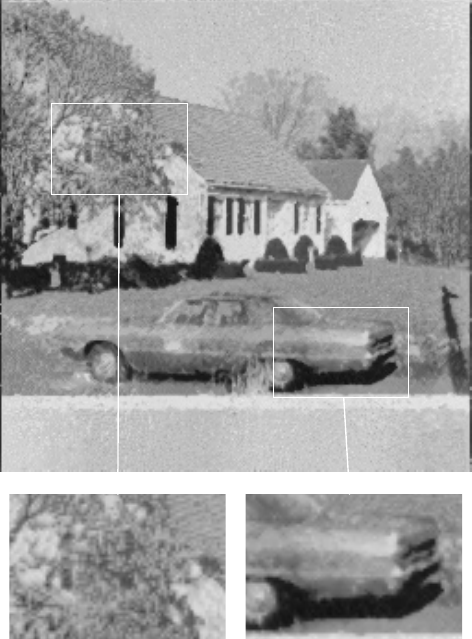}
		\end{subfigure}
		
		\vspace{0.15cm}	
		
		\begin{subfigure}[t]{0.32\textwidth}
			\centering
			\includegraphics[width=.98\textwidth]{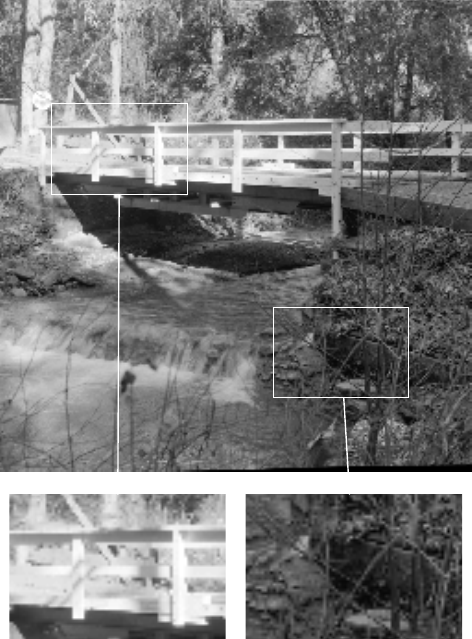}
		\end{subfigure}
		\begin{subfigure}[t]{0.32\textwidth}
			\centering
			\includegraphics[width=.98\textwidth]{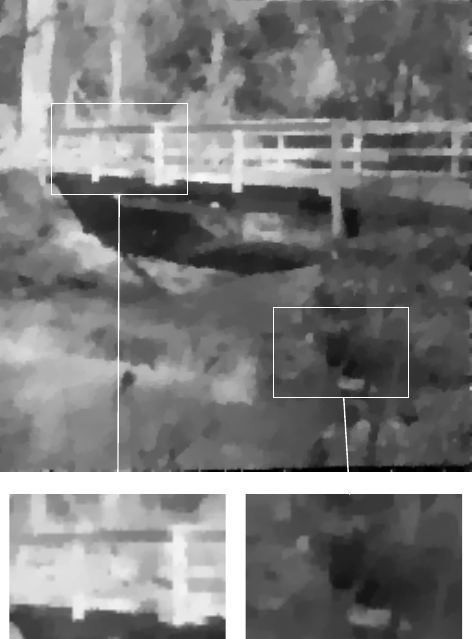}
		\end{subfigure}
		\begin{subfigure}[t]{0.32\textwidth}
			\centering
			\includegraphics[width=.98\textwidth]{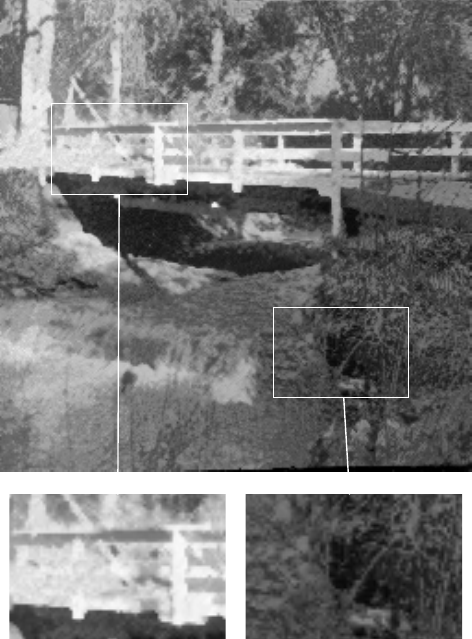}
		\end{subfigure}

		\vspace{0.15cm}	
		
		\begin{subfigure}[t]{0.32\textwidth}
			\centering
			\includegraphics[width=.98\textwidth]{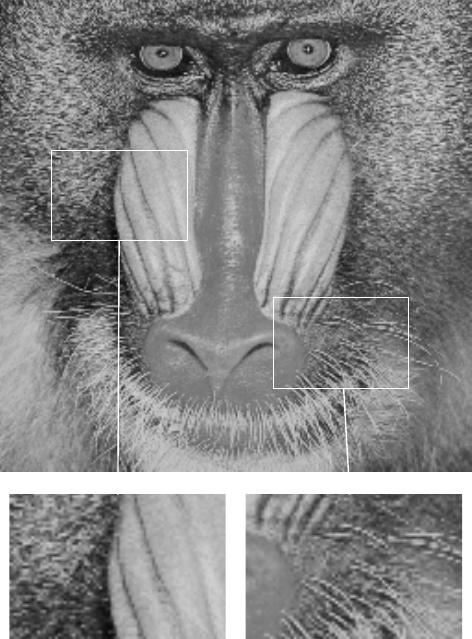}
		\end{subfigure}
		\begin{subfigure}[t]{0.32\textwidth}
			\centering
			\includegraphics[width=.98\textwidth]{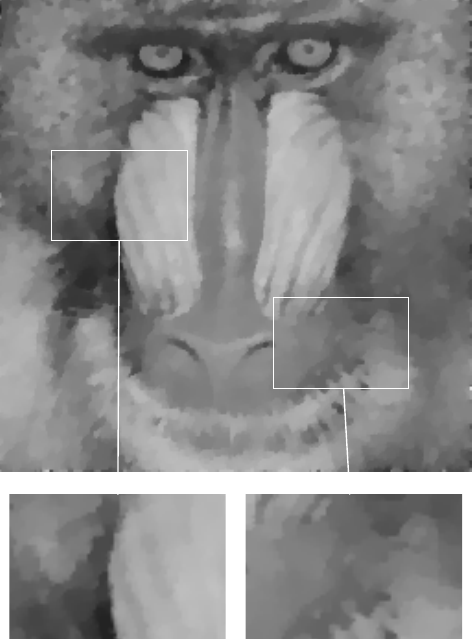}
		\end{subfigure}
		\begin{subfigure}[t]{0.32\textwidth}
			\centering
			\includegraphics[width=.98\textwidth]{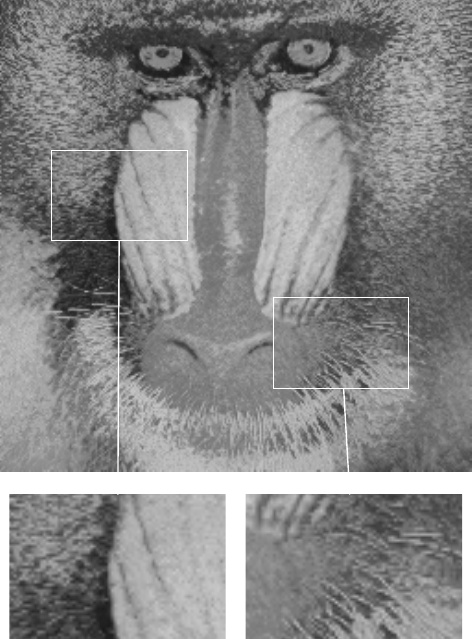}
		\end{subfigure}

		\caption[]{Denoising of an image (left) for $\gamma =10$ using the method proposed in~\cite{MDHY16} (middle) and our N-GMF (right).}
		\label{Fig:examples_gamma_10_2}	
	\end{figure*}

	\begin{table*}[thb]
			\centerline{
				\begin{tabular}{l@{\;}|c@{\;}c@{\;}c@{\;}c}
					\hline
					\hline
					Image & PSNR of~\cite{MDHY16} & PSNR of GNMF & SSIM~\cite{MDHY16}  & SSIM of GNMF\\
					\hline
					\hline
					&&$\gamma=5$ &&\\
					\hline
					test & {\textbf {35.5390}} & 34.6211 &  {\textbf{0.8990} } & 0.8318 \\  
					plane & {\textbf{29.0401}} & 28.4624 &  {\textbf{0.8539}} & 0.8400 \\  
					cameraman & {\textbf{28.9010}} & 28.5065 & 0.8152 & {\textbf{0.8312}} \\  
					boat & {\textbf{29.1804}}  & 28.9941 &  {\textbf{0.8472}}  & 0.8350 \\  
					parrot & {\textbf{29.0874}} & 28.9659 &  {\textbf{0.8607}} & 0.8497 \\  
					leopard & 27.5611 & {\textbf{27.6458}} & 0.8259 & {\textbf{0.8309}} \\  
					peppers & 28.9338 & {\textbf{29.1161}} & 0.8405 & {\textbf{0.8472}} \\  
					house & 27.255  & {\textbf{27.6414}} & 0.8314 & {\textbf{0.8394}} \\  
					bridge & 24.6687 & {\textbf{25.0946}} & 0.7570 & {\textbf{0.7916}} \\  
					montage & 30.9744  & {\textbf{31.7388}} &  {\textbf{0.8913}} & 0.8723 \\  
					baboon & 22.5916 & {\textbf{24.7411}} & 0.6696 & {\textbf{0.7862}} \\  
					barbara & 26.7556 & {\textbf{30.6491}} & 0.8195 & {\textbf{0.8834}} \\  
					\hline
					\hline
					&&$\gamma=10$&& \\
					\hline
					test & 33.5452  & {\textbf{33.9144}} &  {\textbf{0.8651}}  & 0.8410 \\  
					plane & {\textbf{26.9686}}  & 25.4911 &  {\textbf{0.7876}}  & 0.7710 \\  
					cameraman & {\textbf{26.2667}}  & 25.1584 &  {0.6976} & {\textbf{0.7807}} \\  
					boat & {\textbf{26.9668}}  & 25.8286 &  {\textbf{0.7668}}  & 0.7271 \\  
					parrot & 25.9388 & {\textbf{26.1932} } &  {\textbf{0.8113}}  & 0.7876  \\  
					leopard & 23.9575  & {\textbf{24.4084}} &  {{0.7399}} &  {\textbf{0.7505}}  \\  
					peppers & 25.5527  & {\textbf{25.8662}} & 0.6273  & {\textbf{0.6846}} \\  
					house & 24.3855 & {\textbf{24.7098}} & 0.7074  & {\textbf{0.7366}} \\  
					bridge & 22.1702 & {\textbf{22.5982}} & 0.5264  & {\textbf{0.6173}} \\  
					montage & {\textbf{27.7197}} &  25.8032 &  {\textbf{0.8880}}  & 0.8322 \\  
					baboon & 20.6128  & {\textbf{22.0375}} & 0.4030  & {\textbf{0.5976}} \\  
					barbara & 24.0155  & {\textbf{27.9384}} & 0.6845  & {\textbf{0.8121}} \\  				
					\hline				
				\end{tabular}
			}
		\caption{Comparison of  PSNR and SSIM values of the variational method~\cite{MDHY16} and our method. }
		\label{table_psnr}
	\end{table*}
	However, in our opinion, in all except maybe one cases (the \emph{test} image), the visual impression of our restored images is better as they are sharper and provide more details, in particular for a strong noise level $\gamma = 10$. This is illustrated in Figures~\ref{Fig:examples_gamma_5}-\ref{Fig:examples_gamma_10_2}, which show in the left column examples of original images, in the middle the results of~\cite{MDHY16} and in the right the results of our N-GMF. For the \emph{test} image,  the results of~\cite{MDHY16} still contain some corrupted pixels, while our result looks less smooth and slightly grained. Also in other cases, the method of~\cite{MDHY16} does not remove all the outliers, which is in particular visible in constant areas such as the background of the \emph{plane} image or the sky in the \emph{house} image.  In general, our method restores fine structures and textures much better, which can e.g.\ be seen in the whiskers of the \emph{baboon} image, the roof of the \emph{house} or the bushes in the \emph{bridge} image. Furthermore, especially for $\gamma = 10$ the contrast is better preserved with our N-GMF and the results are less blurry, see for instance the \emph{leopard, baboon} or \emph{bridge} image.

	\paragraph{Weighted Generalized Myriad Filtering}
	Finally, we give some results when using  weights in the myriad objective function based on similarity between patches.	Consider a fixed patch $p$ and denote with $q_1,\ldots,q_n$ the patches that are closest to $p$ with respect to the similarity measure $\lambda$. Then, one possibility for choosing weight in the myriad filtering procedure is to use  weights of the form	
	\begin{align}
	w_i \propto \phi_h\bigl(-\log( \lambda(p,q_i))\bigr),\qquad i=1,\ldots,n,\label{weights}
	\end{align}
	followed by normalization such that $\sum\limits_{i=1}^n w_i = 1$. At this point, $\phi_h\colon \R_{\geq 0 } \to \R_{\geq 0 }$ is a continuous kernel function, i.e., $\phi_h$ is nonincreasing, $\phi_h(0)=1$ and $\lim_{t\to \infty} \phi_h(t) = 0$, such as  $\phi_h(t) = \e^{-\frac{t}{h}}$. We chose the latter one, where the parameter $h$ was separately optimized (w.r.t.\ PSNR) for the different noise levels.
	While the difference compared to uniform weights is not visible in the resulting images, PSNR and SSIM values can be further improved, see Table~\ref{Tab:weights}.
	
	\begin{table*}[thb]
			\centerline{
				\begin{tabular}{l@{\;}|c@{\;}c@{\;}c@{\;}c}
					\hline
					\hline
					Image & PSNR  of GNMF & PSNR of wGNMF & SSIM of GNMF & SSIM of wGNMF\\
					\hline
					\hline
					&&$\gamma=5$ &&\\
					\hline
					test &  34.6211 &  \textbf{35.2524}   & 0.8318 & \textbf{0.8341} \\  
					plane &  28.4624 & \textbf{29.0171} & 0.8400 & \textbf{0.8442}\\  
					cameraman &  28.5065 &  \textbf{29.6564}&  {{0.8312}} & \textbf{0.8385}\\  
					boat &   28.9941 &  \textbf{29.4876}  & 0.8350 & \textbf{0.8413}\\  
					parrot &  28.9659 &  \textbf{29.5497}& \textbf{0.8497}& 0.8485\\  
					leopard &  {{27.6458}} & \textbf{27.937}5 & {{0.8309}}& \textbf{0.8322} \\  
					peppers &  {{29.1161}} & \textbf{29.2565}& {{0.8472}} & \textbf{0.8513}\\  
					house &  {{27.6414}} & \textbf{28.1973}& {{0.8394}}&  \textbf{0.8489}\\  
					bridge &  {{25.0946}} & \textbf{25.5402} & {{0.7916}}& \textbf{0.8116}\\  
					montage &  {{31.7388}}& \textbf{32.5221} &  0.8723& \textbf{0.8738}\\  
					baboon &  {{24.7411}} & \textbf{25.0864}& {{0.7862}} & \textbf{0.8029}\\  
					barbara & {{30.6491}} & \textbf{30.9470} & {{0.8834}} & \textbf{0.8842}\\  
					\hline	
					\hline
					&&$\gamma=10$&& \\
					\hline
					test &   {{33.9144}} & \textbf{34.0884} & \textbf{0.8410} & 0.8360\\  
					plane &  25.4911 & \textbf{25.8890} &  \textbf{0.7710} & \textbf{0.7710}\\  
					cameraman &  25.1584 & \textbf{26.6964} &  {{0.7807}}& \textbf{0.7835}\\  
					boat &  25.8286 &  \textbf{26.2730} &   0.7271 & \textbf{0.7362} \\  
					parrot &  26.1932   & \textbf{26.5494} &  \textbf{0.7876} & 0.7854 \\  
					leopard &  {{24.4084}} &\textbf{24.7465} &   {{0.7505}} & \textbf{0.7542}\\  
					peppers &  {{25.8662}}& \textbf{26.0102} & {{0.6846}}& \textbf{0.6945} \\  
					house &  {{24.7098}}& \textbf{25.0779} & {{0.7366}}& \textbf{0.7451}\\  
					bridge &  {{22.5982}}& \textbf{22.8566}  & {{0.6173}} & \textbf{0.6370}\\  
					montage &   25.8032 & \textbf{28.0506}  & \textbf{0.8322} & 0.8265 \\  
					baboon &  {{22.0375}} & \textbf{22.2145} & {{0.5976}} & \textbf{0.6160}\\  
					barbara &  {{27.9384}} & \textbf{28.1885} & {{0.8121}}& \textbf{0.8147}\\  				
					\hline		
				\end{tabular}
			}
		\caption{Comparison of PSNR and SSIM values for GNMF and weighted GNMF. }
		\label{Tab:weights}
	\end{table*}			
	
	\section{Conclusions}\label{sec:conclusions}
	In this work,  we introduced a generalized myriad filtering that can be used 
	to compute the joint ML estimate of the Cauchy distribution. As by-products of our 
	approach we further obtained two algorithms for the estimation of only one parameter 
	if the other one is fixed. Additionally, based on asymptotic estimates  
	we developed a very fast minimization Algorithm~\ref{Alg:myriad_general_fast}. 
	Based on our algorithms we proposed a nonlocal generalized myriad filter for 
	denoising images corrupted by Cauchy noise. It shows an excellent numerical performance, 
	in particular for highly textured images.
	
	There are several directions for future work: 
	On the one hand, in order to handle whole patches 
	and not only single pixels and thereby respecting 
	dependencies between neighboring pixels it might be advantageous to  generalize the myriad filter 
	to a multivariate setting. This is also interesting from a theoretical point of view, since the analysis is 
	in this case much more advanced.
	
	Concerning our denoising approach, fine tuning steps as discussed in~\cite{LCBM12} such as aggregations of patches, 
	the use of an oracle image or a variable patch size to cope better with
	textured and homogeneous image regions may further improve the denoising results. 
	Another question is how to incorporate linear operators (blur, missing pixels) into the image restoration.
	Finally, it would be interesting if other types of impulsive noise can be handled with our filter, 
	for instance impulse noise, Indeed, the myriad filter approaches the mean  for $\gamma\to \infty$ and  thus MFs constitute
	a robust generalization of classical linear mean filters. The other extreme, i.e.\ as $\gamma\to 0$ 
	it converges to the mode-myriad, see~\cite{GA01}.  
	Further it would be interesting to handle a \emph{spatially varying noise level} in our nonlocal method. 
	While our L-GMF can already cope with spatially varying noise, 
	in our N-GMF we would have to adapt the selection 
	of similar patches and the similarity measure, which is a non-trivial~task. 
	
	\subsection*{Acknowledgments} 
	Funding by the German Research Foundation (DFG)  with\-in the Research Training Group 1932,
	project area P3, is gratefully acknowledged. Further, we wish to thank Yiqiu Dong for fruitful discussions and an interesting talk about Cauchy noise removal, as well as the anonymous referees for their careful examination of our manuscript.
	
	\appendix
	\section{Appendix}\label{app:ML_prop}
	This appendix contains the proofs of Section~\ref{sec:ML_prop}.\\
	\textbf{Proof of Theorem~\ref{theo_both}:}
	\begin{proof}
		\begin{enumerate}[1.]
			\item First, we prove that all critical points $(\hat a,\hat{\gamma})$ of $L$ are strict minimizers, by showing that $\nabla^2 L$ is positive definite at all critical points.
			To simplify the notation we set $a_i \coloneqq x_i - a$ and $\hat{a}_i \coloneqq x_i - \hat{a}$, $i=1,\ldots,n$.
			\\[1ex]
			{\textbf{Estimation of $\frac{\partial^2 L}{\partial a^2}$}:} 
			We split the sum in \eqref{der_a2} as follows:
			\begin{align*}
			\sum_{i\colon a_i^2 > \gamma^2} w_i  \frac{a_i^2-\gamma^2 }{\bigl(a_i^2 + \gamma^2\bigr)^2}
			& < \sum_{i\colon a_i^2 > \gamma^2} w_i \frac{a_i^2-\gamma^2 }{\bigl( a_i^2 + \gamma^2\bigr) (\gamma^2 + \gamma^2 )}\\
			&=  \frac{1}{2\gamma^2} \sum_{i\colon a_i^2 > \gamma^2}w_i \frac{a_i^2-\gamma^2 }{  a_i^2 + \gamma^2 }
			\end{align*}
			and similarly, since in this case the summands are negative,
			\begin{equation*}
			\sum_{i\colon a_i^2<\gamma^2} w_i \frac{a_i^2-\gamma^2 }{\bigl(a_i^2 + \gamma^2\bigr)^2} <  \frac{1}{2\gamma^2} \sum_{i\colon a_i^2<\gamma^2}w_i \frac{a_i^2-\gamma^2 }{  a_i^2 + \gamma^2 }.
			\end{equation*}
			Since $n \ge 3$ at least one of these sums exists.
			Thus, for $\hat{\gamma}$ we have
			\begin{align}
			\frac{\partial^2 L}{\partial  a^2}(\hat{a}, \hat{\gamma}) &= - 2 \ \sum_{i=1}^n w_i\frac{\hat{a}_i^2 -\hat\gamma^2}{\bigl(\hat{a}_i^2 + \hat\gamma^2\bigr)^2}\\
			&>- \frac{1}{\hat\gamma^2} \sum_{i=1}^n w_i\frac{\hat{a}_i^2-\hat\gamma^2 }{  \hat{a}_i^2 + \hat\gamma^2 }\\
			& = - \frac{1}{\hat\gamma^2} \sum_{i=1}^n w_i \frac{\hat{a}_i^2+\hat\gamma^2-2\hat\gamma^2 }{  \hat{a}_i^2 + \hat\gamma^2 }\\
			& = - \frac{1}{\hat\gamma^2} + 2 \sum_{i=1}^n w_i \frac{1 }{  \hat{a}_i^2 + \hat\gamma^2 } = 0,
			\end{align}
			where the last equation follows by \eqref{cond_gamma}.
			\\[1ex]
			{\textbf{Estimation of $\det \left( \nabla^2 L(\hat a,\hat\gamma) \right)$}}: Using \eqref{der_a2} - \eqref{der_gamma2} we obtain
			\begin{equation*}
			\frac14 \det \left( \nabla^2 L(\hat a,\hat\gamma) \right) = \left(\sum\limits_{i=1}^n w_i \frac{\hat\gamma^2 - \hat a_i^2}{\bigl(\hat a_i^2 + \hat\gamma^2\bigr)^2}\right)
			\left(\sum\limits_{i=1}^n w_i \frac{\hat a_i^2-\hat\gamma^2}{\bigl(\hat a_i^2 + \hat\gamma^2\bigr)^2}+\frac{1}{2\hat\gamma^2} \right)\\
		 -4\left( \sum\limits_{i=1}^n w_i\frac{\hat\gamma \hat{a}_i}{\bigl(\hat a_i^2 + \hat\gamma^2\bigr)^2}\right)^2.
			\end{equation*}
			With \eqref{cond_gamma} we rewrite the first term as
			\begin{align}\sum_{i=1}^n w_i\frac{\hat a_i^2-\hat \gamma^2 }{  \bigl(\hat a_i^2 + \hat \gamma^2\bigr)^2 }
			&=  	\sum_{i=1}^n w_i\frac{\hat a_i^2+ \hat \gamma^2 -2\hat \gamma^2 }{  \bigl(\hat a_i^2 + \hat \gamma^2\bigr)^2 }\\
			&= \frac{1}{2\hat \gamma^2} - 2\hat \gamma^2 \sum_{i=1}^n w_i\frac{1}{  \bigl(\hat a_i^2 + \hat \gamma^2\bigr)^2 }\\
			& = -\frac{1}{2\hat \gamma^2} + 2 \sum_{i=1}^n w_i\frac{1}{\hat a_i^2 + \hat \gamma^2} - 2 \hat \gamma^2\sum_{i=1}^n w_i\frac{1}{  \bigl(\hat a_i^2 + \hat \gamma^2\bigr)^2 } \\
			& = -\frac{1}{2\hat \gamma^2} \sum_{i=1}^n w_i\frac{  \bigl(\hat a_i^2- \hat \gamma^2\bigr)^2 }{  \bigl(\hat a_i^2 + \hat \gamma^2\bigr)^2 } \label{mixed_2}.
			\end{align}
			With the help of~\eqref{cond_gamma} we can simplify the second term as
			\begin{align} 
			\frac{\partial^2 L}{\partial \gamma^2}(a,\hat \gamma) = 4 \sum_{i=1}^n w_i\frac{a_i^2 }{\bigl(a_i^2 + \hat{\gamma}^2\bigr)^2}
			\end{align}				
			and the third one using \eqref{cond_a} by
			\begin{align}
				\sum_{i=1}^n w_i\frac{\hat a_i }{  \bigl(\hat a_i^2 + \hat \gamma^2\bigr)^2 } 			&= -\frac{1}{2\hat \gamma^2} \sum_{i=1}^n w_i\frac{\hat a_i}{\hat a_i^2 + \hat \gamma^2} + 	\sum_{i=1}^n w_i \frac{\hat a_i }{  \bigl(\hat a_i^2 + \hat \gamma^2\bigr)^2 } \\
			&= -\frac{1}{2\hat \gamma^2}\sum_{i=1}^n w_i\frac{ \hat a_i \bigl(\hat a_i^2 -\hat \gamma^2  \bigr)}{  \bigl(\hat a_i^2 + \hat \gamma^2\bigr)^2 }\label{mixed_1}.
			\end{align}
			Therewith we obtain
			\begin{align}
				\frac14 \det \left( \nabla^2 L(\hat{a},\hat \gamma) \right) = 
			\frac{1}{\hat \gamma^2} \left( 
			\left( \sum_{i=1}^n w_i\frac{  \bigl(\hat a_i^2- \hat \gamma^2\bigr)^2 }{  \bigl(\hat a_i^2 + \hat \gamma^2\bigr)^2 }\right) 
			\left( \sum_{i=1}^n w_i\frac{\hat a_i^2 }{\bigl(\hat a_i^2 + \hat \gamma^2\bigr)^2}\right)- \left(\sum_{i=1}^n w_i\frac{\hat a_i\bigl(\hat a_i^2 -\hat \gamma^2  \bigr)}{  \bigl(\hat a_i^2 + \hat \gamma^2\bigr)^2 }\right)^2  \right)
			\end{align}		
			and by Cauchy-Schwarz' inequality finally
			$$
			\det \left( \nabla^2 L(a,\hat \gamma) \right) >0.
			$$
			Note that we have indeed strict inequality, since otherwise there must exist $\lambda \in \mathbb R$ such that
			$\hat a_i = \lambda (\hat a_i^2 - \hat \gamma^2)$
			for all $i=1,\ldots,n$, which is not possible since $n \ge 3$.
			\item  Next, we show that the first step implies that there is only one critical point.
			For any fixed $a \in (x_1,x_n)$, let $\hat x (a) \coloneqq \hat \gamma (a)^2$ denote the solution of \eqref{cond_gamma} which is unique due to monotonicity, see Lemma \ref{lem_2}.
			Bringing the summands in \eqref{cond_gamma} to the same nominator, 
			we see that $\hat x (a)$ is the unique real zero of a polynomial $P(\cdot,a)$ of degree $n$
			whose coefficients are again polynomials in $a$, say
			\begin{align}
			P(x,a) &=  x^n + p_{n-1} (a) x^{n-1} + \ldots + p_1(a) x + p_0(a).
			\end{align}
			We show that the zero $\hat x(a)$ of $P(\cdot,a)$ is differentiable in $a$. To this end, we
			consider the smooth function $F\colon \mathbb R^n \times \mathbb R \rightarrow \mathbb R$
			given by
			\begin{align*}
			F(c,x)& \coloneqq  x^n + c_{n-1} x^{n-1} + \ldots + c_1 x + c_0,\\
			c &\coloneqq (c_0,\ldots,c_{n-1}).
			\end{align*}
			For an arbitrary fixed $a^* \in (x_1,x_n)$, we define $$c^* \coloneqq (p_0(a^*),\ldots, p_{n-1} (a^*))$$
			and $x^* \coloneqq \hat x (a^*)$.
			Then we have $F(c^*, x^*) = 0$ and, 
			since $x^*$ is a simple zero of $P(\cdot,a^*)$, it holds 
			$
			\frac{\partial}{\partial x} F(c^*,x^*) = P'(x^*,a^*) \not = 0
			$.
			By the implicit function theorem, there exists a continuously differentiable function
			$\varphi\colon \mathbb R^n \rightarrow \mathbb R$ such that 
			$F(c,\varphi(c)) = 0$
			in a neighborhood of $c^*$. 
			Thus, for $c(a) \coloneqq (p_0(a), \ldots, p_{n-1}(a))$ with $a$ in a neighborhood of $a^*$ we have
			$\hat x (a) = \varphi(c(a))$ and
			$$
			\hat x '(a) = \varphi(c(a))' =  \nabla \varphi (c(a)) \cdot \left( p_0'(a) , \ldots,p_{n-1}'(a) \right),
			$$
			which proves the claim.\\		
			Now, the minima of $L$ are given by $(\hat a,\hat \gamma(\hat a))$.
			Assume that there exist two different minimizers $\check a < \tilde a$  of $L$.
			Since they are strict,  $\check a $ and $\tilde a$  are also strict minimizers of the univariate function
			$g(a) \coloneqq L(a,\hat \gamma(a))$. 
			The function  $g(a)$ is continuous, so that there exists a maximizer 
			$\bar a \in (\check a, \tilde a)$
			of $g$ fulfilling
			\begin{equation*}
			0 = g'(\bar a) =  \nabla L(\bar a,\hat \gamma(\bar a)) \cdot \left( (1,\hat \gamma'(\bar a) \right).
			\end{equation*}
			By construction of $\hat \gamma$ we have $\frac{\partial L}{\partial \gamma} (\bar a,\hat \gamma(\bar a)) = 0$.
			This implies  $\frac{\partial L}{\partial a} (\bar a,\hat \gamma(\bar a)) = 0$.
			Consequently, $\nabla L(\bar a,\hat \gamma(\bar a)) = 0$ so that $(\bar a,\hat \gamma(\bar a))$ is critical point of $L$ 
			which is not a strict minimizer. This yields a contradiction and thus we have indeed only one critical point.
			\item Finally, to see the existence of a critical point it remains to show that
			there exists an $a$ such that $S_1(a,\hat \gamma(a)) = 0$.
			By part 2 of the proof $S_1(a,\hat \gamma(a))$ is a continuous function in $a$.
			By the proof of the next Lemma \ref{lem_1} we know that $S_1(a,\hat \gamma(a)) >0$ for $a\le x_1$
			and $S_1(a,\hat \gamma(a)) <0$ for $a\ge x_n$, so that the function has indeed a zero.	
		\end{enumerate}
	\end{proof}
	\textbf{Proof of Lemma~\ref{lem_1}:}
	\begin{proof}
		By \eqref{cond_a}, a critical point of  $L(\cdot,\gamma)$ has to fulfill $s_1(a) = 0$, where
		$s_1 \coloneqq \frac{S_1(\cdot,\gamma)}{\gamma}$.
		All summands in $s_1(a)$ become positive if $a < x_1$ and negative if $a > x_n$.
		Since $n \ge 2$ this implies $s_1(a) > 0$ for $a \le x_1$ and $s_1(a) < 0$ for $a \ge x_n$.
		Hence, the zeros of $s_1$ lie in $(x_1,x_n)$ and there exists at least one zero by
		continuity of $s_1$. Further, 
		\begin{align}
		s_1(a)& =   \frac{P_1(a)}{\prod\limits_{i=1}^n \bigl((x_i-a)^2 + \gamma^2 \bigr)}, \\
		P_1(a) &\coloneqq \sum\limits_{i=1}^n w_i (x_i-a) \prod\limits_{j\neq i} \bigl((x_j-a)^2 + \gamma^2 \bigr),
		\end{align}
		so that the zeros of $s_1$ are the real roots of the nontrivial polynomial $P_1$ of degree $2n-1$, which  are at most $2n-1$.
	\end{proof}
	\textbf{Proof of Lemma~\ref{Lemma:no_saddle_point}:}
	\begin{proof}
		From the previous proof we know that the zeros of $\frac{\partial L}{\partial a}(\cdot,\gamma)$ coincide with those of~$P_1$.
		By \eqref{der_a2} we have 
			\begin{align}
			\frac{\partial^2 L}{\partial a^2}(a,\gamma)
			&=  \frac{-2P_2(a)}{ \prod\limits_{i=1}^n \bigl((x_i-a)^2 + \gamma^2 \bigr)^2}, \\
			P_2(a) 
			&\coloneqq \sum\limits_{i=1}^n w_i \bigl((x_i-a)^2 - \gamma^2\bigr) \prod\limits_{j\neq i} \bigl((x_j-a)^2 + \gamma^2 \bigr)^2,
			\end{align}	
		so that the zeros of $\frac{\partial^2 L}{\partial a^2}(\cdot,\gamma)$ are those of the polynomial $P_2$.
		The coefficients of the polynomials $P_i$, $i=1,2$  are polynomials in $x_1,\ldots,x_n$. 
		Now, \eqref{saddle_point} states that $P_1$ and $P_2$ have a common root,
		which implies that the resultant $\operatorname{Res}(P_1,P_2)$ of $P_1$ and $P_2$  equals zero, see, e.g.~\cite{GKZ08}. 
		The resultant is defined as the determinant of the associated Sylvester matrix, so it is a polynomial expression in the  coefficients of $P_1$ and $P_2$ as well, 
		i.e., a polynomial in $x_1,\ldots,x_n$.
		Since the  set of roots of an $n$-variate polynomial is a set of measure zero in $\R^n$, this finishes the proof.  
	\end{proof}
	\textbf{Proof of Lemma~\ref{lem_2}:}
	\begin{proof}
		\begin{enumerate}[1.]
			\item 	By \eqref{cond_gamma} the critical points of $L(a,\cdot)$ have to fulfill
			$s_0(\gamma^2) = \frac{1}{2}$, where
			$s_0 (\gamma^2) \coloneqq S_0(a,\gamma)$.
			The continuous function $s_0$ is strictly increasing in $\gamma^2$.
			Since $s_0(0) = 0$ and $\lim_{\gamma \rightarrow \infty} s_0(\gamma) = 1$, 
			we conclude that
			$s_0(\gamma^2) = \frac{1}{2}$ has a unique solution $\hat \gamma^2$. 
			Moreover, in~\eqref{der_gamma2}, we obtain 
			\begin{align} \label{der_gamma2_crit}
			\frac{\partial^2 L}{\partial \gamma^2}(a,\hat \gamma) = 4 \sum_{i=1}^n w_i\frac{(x_i-a)^2 }{\bigl((x_i-a)^2 + \hat{\gamma}^2\bigr)^2}>0
			\end{align}
			so that $\hat \gamma$ is a minimizer.
			\item 	Concerning the range of $\gamma$ it follows from~\eqref{cond_gamma} that $\hat \gamma^2 \in \left( \min_i (x_i-a)^2,\max_i (x_i-a)^2 \right)$
			which gives together with the fact $a\in (x_1,x_n)$ the upper bound for $\hat \gamma$.
			To see the lower bound, assume that $\hat \gamma^2 \leq d^2 \epsilon^2$ and distinguish two cases:
			\begin{enumerate}[i)]
				\item First, let $a$ be one of the sample points, say $a=x_i$. Then, since $s_0$ is strictly increasing and $(x_j-a)^2 \ge d^2$ for $j\neq i$, it holds
				\begin{align}
				S_0(a,\gamma) &< w_i + \sum_{j \not = i} w_j \frac{  d^2 \epsilon^2}{d^2 + d^2 \epsilon^2} \\
				&= w_i + (1-w_i)\frac{\epsilon^2}{1+\epsilon^2} \\
				& \leq w_{\text{max}} + (1-w_{\text{max}})\frac{\epsilon^2}{1+\epsilon^2}\\ &= \frac{1}{1+\epsilon^2} \frac{1}{2}\leq\frac12, 
				\end{align}
				which is in contradiction to \eqref{cond_gamma}.
				\item 		Next, let $a \in (x_i,x_{i+1})$. Similarly we obtain in this case
					\begin{align}
						S_0(a,\gamma)	 < w_i \frac{d^2 \epsilon^2}{(a-x_i)^2 + d^2 \epsilon^2} + w_{i+1}  \frac{d^2 \epsilon^2}{(a-x_{i+1})^2 + d^2 \epsilon^2} + \frac{\epsilon^2}{1+\epsilon^2} (1-w_i - w_{i+1}).
					\end{align}
				If the weights are not equal, say $w_i > w_{i+1}$, then the right-hand side becomes largest for $a = x_i$ and we are in the previous case i).
				For $w_i = w_{i+1} = w$ we get
					\begin{align}
					S_0(a,\gamma)				< \epsilon^2 w \left( \frac{1}{(\frac{a-x_i}{d})^2 + \epsilon^2} + \frac{1}{(\frac{a-x_{i+1}}{d})^2 + \epsilon^2}\right)
					& \quad + \frac{\epsilon^2}{1+\epsilon^2}(1-2w)
					\end{align}
				and by replacing $d$ by $x_{i+1} -x_i$ and 
				denoting by $z \in \left[0,\frac12\right]$ the distance of $a$ to the midpoint of the normalized interval,
					\begin{align}
					S_0(a,\gamma) 
					&< \epsilon^2 w \left( \frac{1}{ (\frac12 + z)^2 + \epsilon^2}
					+ \frac{1}{(\frac12 - z)^2 + \epsilon^2}\right) + \frac{\epsilon^2}{1+\epsilon^2}(1-2w)\\
					&= \epsilon^2 w \left( \frac{2( \frac14 + z^2 + \epsilon^2)}{(\frac14 - z^2)^2 + 2 \epsilon^2 (\frac14 + z^2) + \epsilon^4} \right) + \frac{\epsilon^2}{1+\epsilon^2}(1-2w)\\
					&=
					2 \epsilon^2 w \left( \frac{ \frac14 + z^2 + \epsilon^2}{( \frac14 + z^2 + \epsilon^2)^2 - z^2} \right)+ \frac{\epsilon^2}{1+\epsilon^2}(1-2w).\label{estimate_gamma}
					\end{align}	
				Now, $\frac{ \frac14 + z^2 + \epsilon^2}{ \left( \frac14 + z^2 + \epsilon^2\right)^2 - z^2}$ becomes largest iff 
				\begin{equation*}
				\frac{\left( \frac14 + z^2 + \epsilon^2\right)^2 - z^2}{ \frac14 + z^2 + \epsilon^2} 
				= \left(\frac14 + z^2 + \epsilon^2\right) - \frac{z^2}{ \frac14 + z^2 + \epsilon^2}
				\end{equation*}
				becomes smallest. Substituting $y \coloneqq z^2 + \frac14 \in \left[\frac14,\frac12\right]$, we obtain the function
				\begin{align*}
				f(y)& \coloneqq  y + \epsilon^2 - \frac{y - \frac14}{y + \epsilon^2} \\
				&= - 1 +  y + \epsilon^2 + \frac{\frac{1}{4}+\epsilon^2}{y + \epsilon^2},
				\end{align*}
				whose derivatives are given by
				\begin{align*}
				f'(y) = 1-\frac{\epsilon^2 + \frac{1}{4}}{(y + \epsilon^2)^2},\qquad f''(y) = 2 \frac{\epsilon^2 + \frac{1}{4}}{(y + \epsilon^2)^2}.
				\end{align*}
				Setting the derivative to zero results in the positive solution $y = -\epsilon^2 + \sqrt{\epsilon^2 + \frac{1}{4}}$, which is the global minimum on 
				$\left[\frac14,\frac12\right]$ since $f$ is convex. Resubstituting and plugging it in~\eqref{estimate_gamma} yields 
					\begin{align*}
					S_0(a,\gamma)&< 2 \epsilon^2 w \frac{1}{2\sqrt{\epsilon^2 + \frac{1}{4}}+1}+ \frac{\epsilon^2}{1+\epsilon^2}(1-2w) \\
					&\leq  w \epsilon^2+ \frac{\epsilon^2}{1+\epsilon^2}(1-2w)\\
					& =  \underbrace{\epsilon^2\left(\frac{\epsilon^2 -1}{\epsilon^2+1} \right)}_{<0}w		+ \frac{\epsilon^2}{1+\epsilon^2}\\
					& \leq 	\frac{\epsilon^2}{1+\epsilon^2}\leq \frac{1}{3},
					\end{align*}		
				since $\epsilon^2 \in \left(0,\frac{1}{2}\right)$.
			\end{enumerate}		
		\end{enumerate}	
	\end{proof}
	\section{Appendix}\label{app:alg}
	This appendix contains the proofs of Section~\ref{sec:alg}.\\
	\textbf{Proof of Theorem~\ref{Theo:myriad_general}:}
	\begin{proof}
		\begin{enumerate}[1.]
			\item We show that the objective function $L(a_r,\gamma_r)$ decreases for increasing $r$.
			By concavity of the logarithm we have
			\begin{align}
			L(a_{r+1},\gamma_{r+1}) - L(a_r,\gamma_r) &= \sum_{i=1}^n w_i \log \left( \frac{(x_i-a_{r+1})^2 + \gamma_{r+1}^2}{(x_i-a_{r})^2 + \gamma_{r}^2} \frac{\gamma_{r}}{\gamma_{r+1}}\right)\\
			&\le \log \Biggl( \underbrace{\sum_{i=1}^n w_i\frac{(x_i-a_{r+1})^2 + \gamma_{r+1}^2}{(x_i-a_{r})^2 + \gamma_{r}^2} \frac{\gamma_{r}}{\gamma_{r+1}} }_{\Upsilon}\Biggr),
			\end{align}
			so that it suffices to show that $\Upsilon \leq 1$. 
			Setting $S_{0r} \coloneqq S_0(a_r,\gamma_r)$ and $S_{1r} \coloneqq S_1(a_r,\gamma_r)$ we obtain with Algorithm~\ref{Alg:myriad_general}
				\begin{align}
				\Upsilon& = \sqrt{\frac{S_{0r}}{1-S_{0r}}} \sum_{i=1}^n w_i\frac{(x_i-a_r + a_r- a_{r+1})^2 
					+ \frac{1-S_{0r}}{S_{0r}} \gamma_{r}^2}{(x_i-a_{r})^2 + \gamma_{r}^2} \\
				&= \sqrt{\frac{S_{0r}}{1-S_{0r}}} 
				\Bigg( \underbrace{\sum_{i=1}^n w_i\frac{(x_i-a_r )^2 }{(x_i-a_{r})^2 + \gamma_{r}^2}}_{1-S_{0r}}\\
				&
				\quad +2 \underbrace{(a_r - a_{r+1} )}_{-\gamma_r \tfrac{S_{1r}}{S_{0r}}} 
				\underbrace{\sum_{i=1}^n w_i\frac{x_i-a_r  }{(x_i-a_{r})^2 + \gamma_{r}^2}}_{\tfrac{S_{1r}}{\gamma_r}} 
				\\
				&\quad
				+  \sum_{i=1}^n w_i\frac{\overbrace{(a_r-a_{r+1})^2}^{\gamma_r^2\tfrac{S_{1r}^2}{S_{0r}^2}}  }{(x_i-a_{r})^2 + \gamma_{r}^2}\\
				&\quad	+
				\frac{1-S_{0r}}{S_{0r}} \underbrace{\sum_{i=1}^n w_i \frac{\gamma_r^2 }{(x_i-a_{r})^2 + \gamma_{r}^2}}_{S_{0r}} \Bigg)\\
				&= 2 \sqrt{S_{0r} (1-S_{0r})} - 2 \frac{S_{1r}^2}{\sqrt{S_{0r} (1-S_{0r}) }}+\sqrt{\frac{S_{0r}}{1-S_{0r}}} \frac{S_{1r}^2}{S_{0r}^2} 
				S_{0r}\\
				&= 2 \sqrt{S_{0r} (1-S_{0r})} - \frac{S_{1r}^2}{\sqrt{S_{0r} (1-S_{0r})}}.
				\end{align}
			The function 
			\begin{equation*}
			f\colon(0,1) \rightarrow \mathbb R,\quad f(z) \coloneqq 2\sqrt{z(1-z)} - \frac{\alpha^2}{\sqrt{z(1-z)}}
			\end{equation*}
			attains its global maximum in $z = \frac12$, where
			$f(\frac12) = 1 - 2\alpha^2 \leq 1$. Consequently,
			$
			\Upsilon \leq 1
			$
			with equality if and only if $S_{1r} = 0$ and $S_{0r} = \frac12$, that is, $(a_{r+1},\gamma_{r+1}) = (a_r,\gamma_r)$.
			\item By \eqref{cond_gamma} and \eqref{cond_a} we know that $(a_{r+1},\gamma_{r+1}) = (a_r,\gamma_r)$ 
			is a fixed point of $(a_{r+1},\gamma_{r+1}) \coloneqq T(a_r,\gamma_r)$
			in Algorithm \ref{Alg:myriad_general} if and only if it is the minimizer of $L$.
			Let $(a_{r+1},\gamma_{r+1}) \not = (a_r,\gamma_r)$ for all $r \in \mathbb N_0$.
			The sequence $\{(a_r,\gamma_r)\}_{r\in \N}$ is bounded: for $a_r$ we have by \eqref{update_a_convex} that $a_r$ is always a convex combination of the $x_i$ so that
			$a_r \in (x_1,x_n)$; for $\gamma_r$ this follows from Lemma~\ref{lem_2} and Theorem~\ref{Theo:convergence_gamma}, that is shown later on.
			Together with part 1 of the proof we see that $L_r \coloneqq L(a_r,\gamma_r)$ is a strictly decreasing, bounded sequence of numbers
			which must converge to some number~$\hat L$.
			Further, $\{ (a_r,\gamma_r)\}_{r\in \N}$ contains a convergent subsequence $\{ (a_{r_j},\gamma_{r_j})\}_{j\in \N}$
			which converges to some $(\hat a,\hat \gamma)$. By the continuity of $L$ and  $T$  we obtain
			\begin{align}
			L (\hat a,\hat \gamma) 
			&= \lim_{j\rightarrow \infty} L(a_{r_j}, \gamma_{r_j}) = \lim_{j\rightarrow \infty} L_{r_j}=\lim_{j\rightarrow \infty} L_{r_j+1} \\
			&= \lim_{j\rightarrow \infty} L(a_{r_j+1}, \gamma_{r_j+1}) \\
			&	= \lim_{j\rightarrow \infty} L\left(T(a_{r_j}, \gamma_{r_j}) \right)
			= L\left(T(\hat a,\hat \gamma) \right).
			\end{align}
			However, this implies $(\hat a,\hat \gamma) = T(\hat a,\hat \gamma)$ so that $(\hat a,\hat \gamma)$ is a fixed point of $T$ and
			consequently the minimizer.
			Since the minimizer is unique, the whole sequence $\{(a_r,\gamma_r)\}_{r\in \N}$ converges to $(\hat a,\hat \gamma)$
			and we are done.
		\end{enumerate}
	\end{proof}
	\textbf{Proof of Theorem~\ref{Theo:descent}:}
	\begin{proof}
		We follow the lines of the proof of Theorem~\ref{Theo:myriad_general}. Recall, that
		\begin{equation}
		Q(a) = L(a,\gamma) + \log(\gamma) = \sum_{i=1}^n w_i \log \left((x_i-a)^2 + \gamma^2 \right).	
		\end{equation}
		\begin{enumerate}
			\item First, we show that the objective function $Q(a_r)$ decreases for increasing $r$.
			By concavity of the logarithm we have
			\begin{align}
			Q(a_{r+1})-Q(a_r) & = 
			L(a_{r+1},\gamma) - L(a_r,\gamma) \\
			&= \sum_{i=1}^n w_i \log \left( \frac{(x_i-a_{r+1})^2 + \gamma^2}{(x_i-a_{r})^2 + \gamma^2} \right)\\
			&\leq \log \Biggl( \underbrace{\sum_{i=1}^n w_i\frac{(x_i-a_{r+1})^2 + \gamma^2}{(x_i-a_{r})^2 + \gamma^2}  }_{\Upsilon}\Biggr),
			\end{align}
			and it suffices to show that $\Upsilon \leq 1$.
			Setting $S_{0r} \coloneqq S_0(a_r,\gamma)$ and $S_{1r} \coloneqq S_1(a_r,\gamma)$ (note that $\gamma$ is fixed here) we obtain with Algorithm~\ref{Alg:myriad_a}
			\begin{align*}
			\Upsilon & =	\sum_{i=1}^n w_i\frac{(x_i-a_{r+1})^2 + \gamma^2}{(x_i-a_{r})^2 + \gamma^2}\\
			& = 1-S_{0r} - \frac{S_{1r}^2}{S_{0r}} + S_{0r}= 1 - \frac{S_{1r}^2}{S_{0r} }\leq 1
			\end{align*}
			with equality if and only if $S_{1r} = 0$, i.e.\ $a_{r+1}= a_r$, in which case $\hat a \coloneqq a_r$ is a critical point of $Q$.
			\item If $a_{r+1}\not = a_r$ for all $r \in \mathbb N_0$, 
			the sequence $Q_r \coloneqq Q( a_r)$ is strictly decreasing
			and bounded   below by $\log (\gamma^2)$, so that $Q_r \to \hat Q$ as $r \to \infty$.
			Further, since $Q$ is continuous and coercive, the sequence $\{a_r\}_{r\in \N}$ is bounded. 
			Consequently, it contains a convergent subsequence $\{a_{r_j} \}_{j\in \N}$
			which converges to some $\hat a$ and by continuity of $Q$ we have $Q(\hat a) = \hat Q$.
			By continuity of $Q$ and the operator $T_1$ given by $a_{r+1} \coloneqq T_1(a_r)$ 
			in Algorithm~\ref{Alg:myriad_a} it follows
			\begin{align} 
			Q(\hat a)  
			&=  \lim_{j \rightarrow \infty} Q(a_{r_j}) = \lim_{j \rightarrow \infty} Q_{r_j} 
			= \lim_{j \rightarrow \infty} Q_{r_j+1}\\
			& =  \lim_{j \rightarrow \infty} Q(a_{r_j+1})
			= \lim_{j \rightarrow \infty} Q(T_1(a_{r_j}) = Q(T_1(\hat a)).
			\end{align}
			By the first part of the proof this implies that $\hat a$ is a fixed point of $T_1$ and thus a critical point of $Q$.
			\item Observing that $\frac{S_{1r}^2}{S_{0r}} = \frac{S_0r}{\gamma^2}(a_{r+1}-a_r)^2$ and $- \log(1-y) \ge y$, $y \in (0,1)$, we have
			\begin{align}	
			Q(a_r) - Q(a_{r+1})& \geq - \log \left(1- \frac{S_{0r}}{\gamma^2}(a_{r+1}-a_r)^2\right)\\
			&\geq \frac{S_{0r}}{\gamma^2} (a_r-a_{r+1})^2.		
			\end{align}
			Since $1+|x-y|^2 <2 (1+|x|^2)(1+|y|^2)$  and $a_r\in [x_1,x_n]$ we estimate
				\begin{align*}
				\frac{S_0r}{\gamma^2} &= \sum_{i=1}^n w_i \frac{1 }{(x_i-a_r)+\gamma^2} \geq \frac{\gamma^2}{2\bigl( a_r^2+\gamma^2 \bigr)} \sum_{i=1}^n w_i \frac{1 }{x_i^2+ \gamma^2} \\
				&\geq \frac{\gamma^2}{2\bigl(\max\{ x_1^2,x_n^2\} +\gamma^2 \bigr)} \sum_{i=1}^n w_i\frac{1 }{x_i^2+\gamma^2} \eqqcolon \tau_0 > 0,
				\end{align*}	
			which results in
			\begin{equation*}
			Q(a_r) - Q(a_{r+1}) \ge \tau_0 (a_r-a_{r+1})^2.
			\end{equation*}
			Since by the second part of the proof $\lim_{r\rightarrow \infty} Q(a_r) - Q(a_{r+1}) = 0$ we also have
			$\lim_{r\rightarrow \infty} |a_r-a_{r+1}| = 0$.	
			\item  Assume now that there exists a subsequence $\{a_{r_l} \}_{l\in \N}$ which converges to some $a^\ast \neq \hat a$. 
			Since the set of critical points is finite, there exists $\varepsilon > 0$ such that $|\hat a - a^\ast| \ge \varepsilon$.
			On the other hand we have by the third part of the proof for $l,j$ large enough that
			$\varepsilon > |a_{r_l} - a_{r_j}|$. For $l,j \rightarrow \infty$ this leads to a contradiction.
		\end{enumerate}
	\end{proof}
	\textbf{Proof of Theorem~\ref{Theorem:convergence_local_min}:}
	\begin{proof}
		By Theorem~\ref{Theo:descent} we know that
		$\hat a = \lim_{r\to \infty} a_r$ exists 
		and $\hat a = T(\hat a)$ is a stationary points of $Q$ fulfilling $Q'(\hat a)  = 0$.
		We distinguish the following cases:
		\begin{align}
		\mbox{Case  I :} \; \exists\, r_0\in \N\colon & a_{r_0 + 1} =  a_{r_0},\\ 
		\mbox{Case  II:} \; \forall\, r  \in \N\colon & a_{r+1}\neq   a_r.
		\end{align}
		We show that Case I occurs with probability zero and that in Case II, the probability of $\hat a$ being a local minimum is one.
		By~\eqref{update_a_convex} we get
		\begin{align}
		a_{r+1} &=  T(a_r) 
		= \frac{	\sum\limits_{i=1}^n w_i x_i \prod\limits_{j\neq i}\left((x_j-a_r)^2 + \gamma^2\right)}
		{	\sum\limits_{i=1}^n \prod\limits_{j\neq i}\left( (x_j-a_r)^2 + \gamma^2\right)}.
		\end{align}
		Rearranging yields 
		\begin{align*}
		0 & = a_{r+1}	
		\underbrace{ \sum\limits_{i=1}^n \prod\limits_{j\neq i} \left((x_j-a_r)^2 + \gamma^2\right)}_{p_1(a_r)} \\
		&\quad- \underbrace{\sum\limits_{i=1}^n w_i x_i \prod\limits_{j\neq i}\left((x_j-a_r)^2 + \gamma^2\right)}_{p_2(a_r)}\\
		& = a_{r+1} p_1(a_r) - p_2(a_r) ,
		\end{align*}
		where $p_1$ and $p_2$ are polynomials. 
		This polynomial equation in 	$a_r$ has only finitely many solutions  (up to $2(n-1)$). 
		Recursively, backtracking each of the possible values for $a_r$ in a similar way, we end up with at most $2^{r+1} (n-1)^{r+1}$ 
		starting points $a_0 \in (x_1,x_n)$ that can lead to the point $a_{r+1}$ after exactly $r+1$ iterations.
		
		\textbf{Case I}: As seen above, there are only finitely many starting points $a_0$ 
		for which the sequence $\{ a_r\}_{r\in \N}$ 
		reaches a fixed point after exactly $r_0$ steps. 
		Since the set of natural numbers $\N$ is countable and countable unions of finite sets are countable, 
		the set of starting points leading to Case I is  countable and has consequently Lebesgue measure zero.
		
		\textbf{Case II}: Since $Q$ is smooth, there might occur the following cases for the critical point~$\hat a$:
		a) $Q''(\hat a) <0$ (local maximum),
		b) $Q''(\hat a) =0$ (saddle point),
		c) $Q''(\hat a) >0$ (local minimum).
		Indeed case a) cannot happen since we have seen in the proof of Theorem \ref{Theo:descent} that
		$\{Q_r\}_{r\in \N}$ is decreasing.
		Addressing case b), according to Lemma~\ref{Lemma:no_saddle_point}  
		the  function $Q$ has with probability one only minima and maxima, but no saddle points.
		Since cases a) and b) occur each with probability zero, case c) occurs with probability one.
		This finishes the proof.
	\end{proof}
	\textbf{Proof of Theorem~\ref{Theo:convergence_gamma}:}
	\begin{proof} 
		\begin{enumerate}
			\item  First, we show property \eqref{mono}.
			From~\eqref{cond_gamma} we see immediately
			\begin{equation*}
			S_0(a,\gamma)
			\begin{cases}
			< \frac12& \text{if } \gamma < \hat \gamma,\\
			= \frac12& \text{if } \gamma = \hat \gamma,\\
			> \frac12& \text{if } \gamma > \hat \gamma,
			\end{cases}	
			\end{equation*}
			so that $\gamma_r < \hat \gamma$ implies
			$\gamma_{r+1} > \gamma_r$ and $\gamma_r > \hat \gamma$  results in $\gamma_{r+1} < \gamma_r$.
			To see that the iterates cannot skip $\hat \gamma$ we consider the quotient
			\begin{align}
			\frac{\gamma_{r+1}^2}{\hat \gamma^2} = \frac{\gamma_r^2}{\hat \gamma^2}\frac{1-S_0(a,\gamma_r)}{S_0(a,\gamma_r)}
			= \frac {\sum\limits_{i=1}^n w_i \frac{a_i^2}{\alpha a_i^2+ \hat\gamma^2} }{\sum\limits_{i=1}^n w_i \frac{\hat \gamma^2}{\alpha a_i^2 + \hat \gamma^2 }},
			\end{align}
			where $\alpha \coloneqq \left( \frac{\hat \gamma }{\gamma_r} \right)^2$.	We have to show that $\alpha < 1$ implies $\frac{\gamma_{r+1}^2}{\hat \gamma^2} > 1$ and conversely,
			$\alpha > 1$ implies $\frac{\gamma_{r+1}^2}{\hat \gamma^2} < 1$.
			Alternatively, we can prove that the function
			\begin{equation*}
			f(\alpha) = \sum_{i=1}^n w_i \frac{\gamma^2-a_i^2 }{\alpha a_i^2 + \gamma^2}, \
			\end{equation*}
			fulfills 
			\begin{equation}\label{claim}
			f(\alpha) \left\{
			\begin{array}{ll}
			< 0 & \mathrm{if} \; \alpha \in (0,1), \\
			> 0 & \mathrm{if} \; \alpha \in (1,+\infty).
			\end{array}
			\right.
			\end{equation}
			We have $f(1) = 0$ and the derivatives of $f$ are given by
			\begin{align}
			f'(\alpha) &= \sum_{i=1}^nw_i \frac{a_i^2(a_i^2 -\gamma^2)}{(\alpha a_i^2 + \gamma^2)^2},\\
			f''(\alpha) &=2 \sum_{i=1}^n w_i \frac{a_i^4(\gamma^2-a_i^2)}{(\alpha a_i^2 + \gamma^2)^3}.
			\end{align}
			For $f'$ we estimate similarly as in the proof of Theorem~\ref{theo_both},
				\begin{align}
				\sum_{i\colon a_i^2 
					> \gamma^2} w_i \frac{a_i^2(a_i^2 -\gamma^2)}{(\alpha a_i^2 + \gamma^2)^2} 
				& > 	\sum_{i\colon a_i^2 
					> \gamma^2}^n w_i \frac{a_i^2(a_i^2 -\gamma^2)}{(\alpha a_i^2 + a_i^2)(\alpha a_i^2 + \gamma^2)} \\
				&=\frac{1}{\alpha+1}\sum_{i\colon a_i^2 > \gamma^2}^n w_i \frac{a_i^2 - \gamma^2 }{\alpha a_i^2 + \gamma^2}  
				\end{align}	
			and analogously  for the negative summands, so that in summary
				\begin{align}\label{estimate_fp}
				f'(\alpha) & > \frac{1}{\alpha+1}\sum_{i=1}^n w_i \frac{a_i^2 -\gamma^2}{\alpha a_i^2 + \gamma^2}
				= -\frac{1}{\alpha+1}f(\alpha).
				\end{align}
			Analogously, we obtain for $f''$ 
		
				\begin{align}
				f''(\alpha) & < \frac{2}{\alpha+1}\sum_{i=1}^n w_i \frac{a_i^2(\gamma^2-a_i^2)}{(\alpha a_i^2 + \gamma^2)^2}
				= -\frac{2}{\alpha+1}f'(\alpha)\label{estimate_fpp}.
				\end{align}
		
			From~\eqref{estimate_fp} it follows $f'(1) > \frac12 f(1) = 0$ and therewith further $f(\alpha) <0 $ for $\alpha\in (0,1)$. 
			Consider the case $\alpha > 1$.
			By continuity  of $f$, we have $f(\alpha)>0$ for $\alpha$ sufficiently close to 1. 
			Since
			$\lim_{\alpha\to \infty} f(\alpha) = 0$ we conclude that $f'$ has at least one root for $\alpha > 1$.
			On the other hand, $f'$ has at most one root, since according to~\eqref{estimate_fpp} any root of $f'$ is a local maximum of $f$.
			Thus, $f'$ has exactly one root, so that $f$ has exactly one critical point (a local maximum). 
			Since $\lim_{\alpha\to \infty} f(\alpha) = 0$ this implies $f(\alpha)>0$ for all $\alpha >1$. 
			\item Since $L(a,\cdot)$ is continuous and has only one critical point, it follows immediately from \eqref{mono} that
			$L(a,\gamma_r) \geq L(a,\gamma_{r+1}) = L(a,T_2 (\gamma_r))$ 
			with equality if and only if $\gamma_r = \gamma_{r+1} = \hat \gamma$.
			In the latter case we are done, so assume that $\gamma_r \not = \gamma_{r+1}$ for all $r \in \mathbb N_0$.
			By part~1  of the proof, the sequence $\{\gamma_r\}_r$ is monotone and  bounded, so it converges to some~$\gamma^*$. By continuity of $L(a,\cdot)$ and $T_2$ we get
			\begin{align*}
			L(a,\gamma^*) &= \lim_{r \rightarrow \infty} L(a,\gamma_r) = \lim_{r \rightarrow \infty} L(a,\gamma_{r+1})\\
			&= \lim_{r \rightarrow \infty} L(a,T_2(\gamma_r)) = L(a,T_2(\gamma^*)),
			\end{align*}
			which is only possible it $\gamma^* = T_2(\gamma^*)$, i.e., if $\gamma^* = \hat \gamma$.	
		\end{enumerate}
	\end{proof}
	
	\section{Appendix}\label{app:asymp}
	This appendix contains the proofs of Section~\ref{sec:asymp}.\\
	\textbf{Proof of Lemma~\ref{Theo:expected_values}:}
	\begin{proof}
		Consider the functions 
		$g(x) \coloneqq \frac{1}{1 + x^2}$ and $h(x)  \coloneqq \frac{x}{1 + x^2}$.
		Both functions are measurable and since
		\begin{align*}
		\int_{-\infty}^\infty  g(x) p(x)\dx &\leq \int_{-\infty}^\infty  p(x) \dx=1 ,\\
		\int_{-\infty}^\infty |h(x) p(x)|\dx  &\leq \int_{-\infty}^\infty \frac{1}{2}  p(x) \dx=\frac{1}{2},
		\end{align*}
		where $p$ denotes the density function of $C(a,\gamma)$, the expected values $\E(Y),\E(Z)$ exist.\\	
		For $g$ and $a\neq 0$, we compute
		\begin{align*}
		&\int_{-\infty}^\infty g(x) p(x)\dx  = \int_{-\infty}^\infty \frac{1}{1+x^2} \frac{1}{\pi \gamma}\frac{\gamma^2}{(x-a)^2 + \gamma^2} \dx\\
		& {=  \left[\frac{ a\gamma \log\left(\frac{(x-a)^2 + \gamma^2}{x^2 + 1}\right) + \gamma(a^2 + \gamma^2-1)\arctan(x) 
				+ (a^2 -\gamma^2 + 1)\arctan\left(\tfrac{x-a}{\gamma}\right)}{\pi\bigl(a^4 + 2a^2(\gamma^2 + 1) + (\gamma^2-1)^2\bigr) } \right]_{-\infty}^{\infty}}.
		\end{align*}
		Since $\lim_{x\to \pm \infty}\log\left(\frac{(x-a)^2 + \gamma^2}{x^2 + 1}\right) = 0$ and $\lim_{x\to \pm \infty} \arctan(x) = \pm \frac{\pi}{2}$, 
		we obtain the first equation in \eqref{E1}.
		For $a=0$ we have 
		\begin{align}
		\int_{-\infty}^\infty g(x) p(x)\dx& = \int_{-\infty}^\infty \frac{1}{1+x^2} \frac{1}{\pi \gamma}\frac{\gamma^2}{x^2 + \gamma^2} \dx\\
		& =\left[\frac{\gamma^2 \arctan(x) - \gamma\arctan\left(\tfrac{x}{\gamma}\right) }{\pi \gamma(\gamma^2-1)}\right]_{-\infty}^\infty,
		\end{align}
		which results in the second equation in \eqref{E1}.
		Similarly, for $h$ and $ (a,\gamma) \not = (0,1)$,
		\begin{align}
		&\int_{-\infty}^\infty h(x) p(x) \dx = \int_{-\infty}^\infty \frac{x}{1+x^2} \frac{1}{\pi \gamma}\frac{\gamma^2}{(x-a)^2 + \gamma^2} \dx\\
		& =\left[\frac{\gamma(a^2 + \gamma^2 - 1) \log\left(\frac{x^2 + 1}{(x-a)^2 + \gamma^2}\right) 
			+ 2a(a^2 + \gamma^2 + 1)\arctan\left(\tfrac{x-a}{\gamma}\right)-4a\gamma\arctan(x)}{2\pi\bigl(a^4 + 2a^2(\gamma^2 + 1) + (\gamma^2-1)^2\bigr) }\right]_{-\infty}^\infty,
		\end{align}
		so that the first equation in \eqref{E2} follows.
		Finally, for $(a,\gamma)  = (0,1)$  it holds
		\begin{align}
		\int_{-\infty}^\infty h(x) p(x)\dx &= \int_{-\infty}^\infty \frac{x}{1+x^2} \frac{1}{\pi }\frac{1}{x^2 + 1} \dx\\
		& = \left[-\frac{1}{2\pi(1 + x^2)}\right]_{-\infty}^\infty,	
		\end{align}
		and consequently $E(Z) = 0$.	
	\end{proof}
	\textbf{Proof of Corollary~\ref{Coro:mean_Y_ir}:}
	\begin{proof}
		By Proposition~\ref{Prop:properties_Cauchy} we have $X_r\sim C\left(\tfrac{a-a_r}{\gamma_r},\tfrac{\gamma}{\gamma_r}\right)$.
		Setting
		$\tilde{a} = \frac{a-a_r}{\gamma_r}$ and $\tilde{\gamma} = \frac{\gamma}{\gamma_r}$ and applying the results of Lemma~\ref{Theo:expected_values} we obtain
		\begin{align}
		\E\bigl(g(X_r)\bigr)  
		&= 	\frac{\tilde{\gamma}(\tilde{a}^2 + \tilde{\gamma}^2-1) + \tilde{a}^2-\tilde{\gamma}^2 + 1}{(\tilde{a}^2 + \tilde{\gamma}^2 + 1)^2 - 4\tilde{\gamma}^2}  \\
		& \textstyle= 	\frac{\frac{\gamma}{\gamma_r}\left[\left(\frac{a-a_r}{\gamma_r}\right)^2 + \left(\frac{\gamma}{\gamma_r}\right)^2-1\right]+\left(\frac{a-a_r}{\gamma_r}\right)^2 
			- \left(\frac{\gamma}{\gamma_r}\right)^2+1}{\left[\left(\frac{a-a_r}{\gamma_r}\right)^2 + \left(\frac{\gamma}{\gamma_r}\right)^2 + 1\right]^2 - 4\left(\frac{\gamma}{\gamma_r}\right)^2}\\
		& = \frac{\gamma_r(\gamma + \gamma_r) }{(a-a_r)^2 + (\gamma + \gamma_r)^2}.
		\end{align}
		Similarly, we compute
		\begin{align*}
		\E\bigl(h(X_r)\bigr)  
		&= 	\frac{\tilde{a}(\tilde{a}^2 + \tilde{\gamma}^2+1-2\tilde{\gamma})}{(\tilde{a}^2 + \tilde{\gamma}^2 + 1)^2 - 4\tilde{\gamma}^2}  \\
		& = 	\frac{\frac{a-a_r}{\gamma_r}\left[\left(\frac{a-a_r}{\gamma_r}\right)^2 
			+ \left(\frac{\gamma}{\gamma_r}\right)^2+1-2\frac{\gamma}{\gamma_r}\right]}{\left[\left(\frac{a-a_r}{\gamma_r}\right)^2 
			+ \left(\frac{\gamma}{\gamma_r}\right)^2 + 1\right]^2 - 4\left(\frac{\gamma}{\gamma_r}\right)^2}\\
		& = \frac{\gamma_r(a-a_r) }{(a-a_r)^2 + (\gamma + \gamma_r)^2}.
		\end{align*}
	\end{proof}
	\textbf{Proof of Theorem~\ref{conv_komisch}:}
	\begin{proof}
		i) Since $\gamma,\tilde{\gamma}_0>0$, it follows inductively from ~\eqref{gamma_rek_mean} that $\tilde{\gamma}_r>0$. Further, if $\tilde{\gamma}_r<\gamma$, then
		\begin{equation*}
		\tilde{\gamma}^2_{r+1}-\tilde{\gamma}^2_r = \tilde{\gamma}_r(\gamma-\tilde{\gamma}_r) + \tilde{\gamma}_r \frac{(a-\tilde{a}_r)^2}{\gamma  + \tilde{\gamma}_r}>0.
		\end{equation*}
		On the other hand, if $\tilde{\gamma}_r>\gamma$, then it holds
		\begin{equation*}
		\tilde{\gamma}^2_{r+1}-\gamma^2 = \gamma(\tilde{\gamma}_r - \gamma) + \tilde{\gamma}_r \frac{(a-\tilde{a}_r)^2}{\gamma  + \tilde{\gamma}_r}>0.	
		\end{equation*}
		Thus, to summarize it holds $\tilde{\gamma}^2_{r+1}\geq \min\{\tilde{\gamma}^2_r,\gamma^2\}$ and inductively we have $\tilde{\gamma}^2_{r+1}\geq \min\{\tilde{\gamma}^2_0,\gamma^2\}$. 
		\\[1ex]		
		ii) Since $\gamma,\tilde{\gamma}_r>0$, this is a direct consequence of~\eqref{a_rek_mean}.
		\\[1ex]		
		iii)
		Let $q = \max\left\{\frac{1}{2},\frac{\gamma}{\gamma + \tilde{\gamma}_0}\right\}$, so that $\frac{1}{2} \leq q < 1$. For the sequence $\{\tilde{a}_r\}_{r\in \N}$ we estimate
		\begin{align}
		|\tilde{a}_{r+1} - a| 
		& = \left|\tilde{a}_r + \tilde{\gamma}_r \frac{a-\tilde{a}_r}{\gamma + \tilde{\gamma}_r} - a  \right|\\
		& = \left| \left(1-\frac{\tilde{\gamma}_r}{\gamma + \tilde{\gamma}_r}\right) (\tilde{a}_r - a) \right|
		= \frac{\gamma}{\gamma + \tilde{\gamma}_r} |\tilde{a}_r - a|\\
		&\leq q  |\tilde{a}_r - a|
		\leq \ldots \le
		q^{r+1}|\tilde{a}_0 - a| \overset{r\to \infty}{\to} 0.
		\end{align}
		Similarly, we obtain for the sequence $\{\tilde{\gamma}_r\}_{r\in \N}$ ,
			\begin{align}
			&	|\tilde{\gamma}^2_{r+1} - \gamma^2|\\
			& = \left|\tilde{\gamma}_r \left(\gamma + \frac{(a-\tilde{a}_r)^2}{\gamma + \tilde{\gamma}_r}\right) - \gamma^2 \right|\\
			&= \left|\gamma(\tilde{\gamma}_r-\gamma) +\frac{\tilde{\gamma}_r}{\gamma + \tilde{\gamma}_r}(a-\tilde{a}_r)^2  \right| \\
			& = \left|\frac{\gamma}{\gamma + \tilde{\gamma}_r}(\tilde{\gamma}_r^2 - \gamma^2)+\frac{\tilde{\gamma}_r}{\gamma + \tilde{\gamma}_r}(a-\tilde{a}_r)^2  \right|\\
			& \leq \frac{\gamma}{\gamma + \tilde{\gamma}_r} \left|\tilde{\gamma}_r^2 - \gamma^2\right| + \frac{\tilde{\gamma}_r}{\gamma + \tilde{\gamma}_r}(a-\tilde{a}_r)^2 \\
			&\leq q \left|\tilde{\gamma}_r^2 - \gamma^2\right| + (a-\tilde{a}_r)^2
			\leq q \left|\tilde{\gamma}_r^2 - \gamma^2\right| + q^{2r}(a-\tilde{a}_0)^2\\
			& \leq q\bigl(q \left|\tilde{\gamma}_r^2 - \gamma^2\right| + q^{2(r-1)}(a-\tilde{a}_0)^2 \bigr)+ q^{2r}(a-\tilde{a}_0)^2\\
			& \;\;\vdots\\
			& \leq q^{r+1} \left|\tilde{\gamma}_0^2 - \gamma^2\right| + \left(\sum\limits_{i=0}^r q^{r+i}\right) (a-\tilde{a}_0)^2\\
			& = q^{r+1} \left|\tilde{\gamma}_0^2 - \gamma^2\right| + q^r \frac{1-q^{r+1}}{1-q}(\tilde{a}_0 - a)^2\overset{r\to \infty}{\to} 0.	
			\end{align}	
	\end{proof}
	\textbf{Proof of Theorem~\ref{Theo:alg_fast}:}
	\begin{proof}
		By strict concavity of the logarithm function and since $w_i > 0$ 
		we have	
		\begin{equation*}
		L(a_{r+1},\gamma_{r+1}) - L(a_{r},\gamma_{r})
		\leq \log \Biggl(\underbrace{\sum_{i=1}^n w_i\frac{(x_i-a_{r+1})^2 + \gamma_{r+1}^2}{(x_i-a_{r})^2 
				+ \gamma_{r}^2} \frac{\gamma_{r}}{\gamma_{r+1}}}_{\Upsilon}\Biggr),
		\end{equation*}	
		with equality if and only if $(a_r,\gamma_r) = (a_{r+1},\gamma_{r+1})$.
		From Algorithm \ref{Alg:myriad_general_fast} we obtain similarly as in the proof of~\ref{Theo:myriad_general},
			\begin{align}
			\Upsilon &= \frac{S_{0r}^2+ S_{1r}^2}{S_{0r} (1-S_{0r}) - S_{1r}^2}
			\left(
			(1-S_{0r}) - 2 \gamma_r \frac{S_{1r}}{S_{0r}^2+ S_{1r}^2} \frac{S_{1r}}{\gamma_r}  \right. \\
			& \left.  + \frac{S_{0r}S_{1r}^2}{(S_{0r}^2+ S_{1r}^2)^2}
			+ \frac{S_{0r} \left(S_{0r} (1-S_{0r}) - S_{1r}^2 \right)^2}{(S_{0r}^2+ S_{1r}^2)^2}\right)\\
			&=\frac{(S_{0r}^2+ S_{1r}^2)^2 (1-S_{0r}) - 2S_{1r}^2(S_{0r}^2+ S_{1r}^2) + S_{0r} S_{1r}^2 + S_{0r}\left(S_{0r} (1-S_{0r}) - S_{1r}^2 \right)^2 }
				{\left(S_{0r} (1-S_{0r}) - S_{1r}^2 \right)(S_{0r}^2+ S_{1r}^2)}\\
			&=1.
			\end{align}
		Thus, $L(a_{r+1},\gamma_{r+1}) \leq L(a_{r},\gamma_{r})$ with equality if and only if $(a_r,\gamma_r) = (a_{r+1},\gamma_{r+1})$.
		The convergence result follows as in part 2 of the proof of Theorem \ref{Theo:myriad_general}.
	\end{proof}
	
	\section{Appendix}\label{app:myr}
	This appendix contains the proof of Section~\ref{sec:myr}.\\
	\textbf{Proof of Lemma~\ref{Lem:Cauchy_similarity}:}
	\begin{proof}
		Under $\HH_0$ (i.e., $n=2$), the ML estimate is not unique, but one easily verifies using~\eqref{cond_a} that 
		\begin{align*}
		\hat{\theta}&=\frac{1}{2}(x_1 + y_1)\in\argmax_{\theta \in \HH_0}\bigl\{\LL(\theta|x_1,y_1)\bigr\}\\
		&=\argmax_{\theta \in \Theta} \bigl\{\LL(\theta|x_1)\LL(\theta|y_1)\bigr\},
		\end{align*}
		and therewith
		\begin{align*}
		\sup_{\theta\in \HH_0} \LL(\theta|x_1,y_1) = \LL(\hat{\theta}_0|x_1)\LL(\hat{\theta}_0|y_1)
		& = \frac{1}{\pi \gamma } \frac{\gamma^2 }{\left(x_1-\frac{x_1+ y_1}{2}\right)^2 + \gamma^2}  \frac{1}{\pi \gamma } \frac{\gamma^2 }{\left(y_1-\frac{x_1+ y_1}{2}\right)^2 + \gamma^2} \\
		& = \frac{1}{\pi^2 \gamma^2 }\frac{1 }{\left(\left(\frac{x_1- y_1}{2\gamma}\right)^2 + 1\right)^2}.
		\end{align*}
		Note that although the ML estimate $\hat{\theta} $ is not unique, the value of the log-likelihood function does not change when using another solution.\\	
		Under $\HH_1$ (i.e., $n=1$), the ML estimate simply reads as
		\begin{equation*}
		\hat{\theta}_i = \argmax_{\theta \in \Theta}\LL(\theta|x_i)=x_i,\qquad i=1,2,
		\end{equation*}
		resulting in 
		\begin{align*}
		\sup_{\theta\in \Theta} \LL(\theta|x_i)= \frac{1}{\pi \gamma } \frac{\gamma^2 }{\left(x_i-x_i\right)^2 + \gamma^2}=\frac{1}{\pi \gamma },\qquad i=1,2.
		\end{align*}	
		Therewith, the LR statistic becomes
		\begin{align*}
		\lambda(x_1,y_1)& = 
		\frac{\sup\limits_{\theta\in \Theta}\bigl\{\LL(\theta|x_1)\LL(\theta|y_1)\bigr\}}{\sup\limits_{\theta\in \Theta}\bigl\{\LL(\theta|x_1)\bigr\}\sup\limits_{\theta\in \Theta}\bigl\{\LL(\theta|y_1)\bigr\}}\\
		&=\frac{\frac{1}{\pi^2 \gamma^2 }\frac{1 }{\left(\left(\frac{x_1- y_1}{2\gamma}\right)^2 
				+ 1\right)^2}}{\frac{1}{\pi \gamma } \frac{1}{\pi \gamma }  }\\
		&=\left(\left(\frac{x_1- y_1}{2\gamma}\right)^2 + 1\right)^{-2}.
		\end{align*}
	\end{proof}
	
	\bibliographystyle{abbrv}
	\bibliography{Cauchy_myriad}

\end{document}